\documentclass[10pt,reqno,twoside]{amsart}

\usepackage{amsfonts}
\usepackage{paralist}
  \usepackage{graphics}
  \usepackage{epsfig}
 \usepackage[colorlinks=true]{hyperref}
\usepackage{hyperref}

\numberwithin{equation}{section}
\usepackage{latexsym}
\usepackage{amsmath}
\usepackage{amssymb}
\usepackage{mathrsfs}
\usepackage{graphicx}
\usepackage{ifthen}
\usepackage[T1]{fontenc}
\usepackage[latin1]{inputenc}
\usepackage{color}
\newtheorem{Satz}{Theorem}[section]

\newtheorem{Def}[Satz]{Definition}

\newcommand{\oo}{\overline \Omega}

\newcommand{\po}{\partial\Omega}

\newcommand{\dist}{\text {dist}}

\newcommand{\diag}{\text {diag}}

\newcommand{\supp}{\text {supp}}

\newcommand{\sign}{\text {sign}}

\begin{document}

\title[Boundary concentration phenomena]
{
Boundary concentration phenomena for an anisotropic  Neumann problem in $\mathbb{R}^2$
}

\author[Yibin Zhang]{Yibin Zhang}
\address{College of Sciences, Nanjing Agricultural University, Nanjing 210095, China}
\email{yibin10201029@njau.edu.cn}



\begin{abstract}
Given  a  smooth bounded domain  $\Omega$ in $\mathbb{R}^2$,
we study the following anisotropic  Neumann problem
$$
\begin{cases}
-\nabla(a(x)\nabla u)+a(x)u=\lambda a(x) u^{p-1}e^{u^p},\,\,\,\,
u>0\,\,\,\,\,
\textrm{in}\,\,\,\,\,
\Omega,\\[2mm]
\frac{\partial u}{\partial\nu}=0\,\,
\qquad\quad\qquad\qquad\qquad\qquad\qquad
\ \ \ \ \,\qquad\quad\,
\textrm{on}\,\,\,
\partial\Omega,
\end{cases}
$$
where
$\lambda>0$ is a small parameter,
$0<p<2$, $a(x)$ is a positive smooth function
over $\overline{\Omega}$ and $\nu$
denotes the outer unit normal vector to $\partial\Omega$.
Under  suitable assumptions on anisotropic
coefficient $a(x)$, we construct
solutions $u_\lambda$ of this problem with arbitrarily many
mixed interior and boundary bubbles
which concentrate at  totally
different
strict local
maximum or minimal boundary points  of $a(x)$ restricted to $\partial\Omega$, or accumulate to
the same strict local
maximum boundary  point of $a(x)$ over $\overline{\Omega}$ as $\lambda\rightarrow0$.
Furthermore, for these bubbling solutions $u_\lambda$ we compute the delicate expansion of the corresponding energy
$$
\aligned
\beta_\lambda
=\,\frac{\lambda p}{2}
\left(
\int_{\Omega}
a(x)\big(e^{u_\lambda^p}-1\big)dx
\right)^{\frac{2-p}{p}}
\left(
\int_{\Omega}
a(x)u_\lambda^p
e^{u_\lambda^p}dx
\right)^{\frac{2(p-1)}{p}}
\endaligned
$$
and  exhibit the sharp difference of tendency of energy $\beta_\lambda$ between $0<p<1$ and $1<p<2$.
\\
%\noindent{\it Keywords:}\,\,
%Boundary concentration phenomena;
%Anisotropic  elliptic Neumann problem;
%Lyapunov-Schmidt reduction procedure.
%
%
%\noindent{\it Keywords:}\,\,
%Boundary concentration phenomena;
%Lyapunov-Schmidt reduction method;
%Anisotropic Neumann problem.
\end{abstract}

\maketitle

\section{Introduction}
This paper is concerned with  the analysis of solutions to
the  anisotropic  Neumann problem
\begin{equation}\label{1.1}
\left\{\aligned
&-\nabla(a(x)\nabla u)+a(x)u=\lambda a(x) u^{p-1}e^{u^p},\,\,\,\,
u>0\,\,\,\,\,
\textrm{in}\,\,\,\,\,
\Omega,\\[1mm]
&\frac{\partial u}{\partial\nu}=0\,\,
\qquad\quad\qquad\qquad\qquad\qquad\qquad
\ \ \ \ \,\qquad\quad\,
\textrm{on}\,\,\,
\partial\Omega,
\endaligned\right.
\end{equation}
where  $\Omega$ is a smooth  bounded   domain in $\mathbb{R}^2$,
$\lambda>0$ is a small parameter,
$0<p<2$, $a(x)$ is a positive smooth function
over $\overline{\Omega}$ and $\nu$
denotes the outer unit normal vector to $\po$.
Any solution of this problem is exactly a critical point of
 the functional
\begin{equation}\label{1.2}
\aligned
J_\lambda(u)=\frac12\int_{\Omega}a(x)(|\nabla u|^2+u^2)-\frac{\lambda}{p}
\int_{\Omega}a(x)e^{u_{+}^p},\,\,\ \ \,\,u\in H^1(\Omega),
\endaligned
\end{equation}
where $u_{+}=\max\{u,0\}$ and the positivity of
every critical point follows from the maximum principle.

%which is well defined because the critical Moser-Trudinger
%inequality implies the validity of the Sobolev-Orlicz  compact
%subcritical embedding
%$$
%\aligned
%\sup_{u\in H^1(\Omega)\setminus\{0\}}
%\int_{\Omega}a(x)\exp\left(\frac{u^p}{\|u\|_{H^1(\Omega)}^p}\right)<+\infty.
%\endaligned
%$$

We are interested in the existence of solutions of equation (\ref{1.1})
that exhibit  the boundary concentration phenomenon  as the parameter $\lambda$
tends to zero.
This work is strongly  stimulated  by
some extensive research involving
the isotropic case $a(x)\equiv1$ in
equation (\ref{1.1}):
\begin{equation}\label{1.3}
\left\{\aligned
&-\Delta
\upsilon+\upsilon=\lambda  \upsilon^{p-1}e^{\upsilon^p},\,\,\,\ \,
\upsilon>0\,\,\,\,\,
\textrm{in}\,\,\,\,\,
\mathcal{D},\\[1mm]
&\frac{\partial \upsilon}{\partial\nu}=0
\,\ \ \,
\qquad\quad\qquad\qquad\qquad\quad
\textrm{on}\,\,\,
\partial\mathcal{D},
\endaligned\right.
\end{equation}
where $\mathcal{D}$ is a smooth  bounded domain in $\mathbb{R}^N$ with $N\geq2$.
In the case of  $p=1$, this scalar  equation
 is equivalent to an elliptic system representing
the  stationary Keller-Segel chemotaxis system with linear sensitivity:
\begin{equation}\label{1.6}
\left\{\begin{array}{ll}
\Delta \psi-\nabla(\psi\nabla \upsilon)=0,\,\,\,\textrm{in}\quad\,\,\,\,\mathcal{D},\\
\Delta \upsilon-\upsilon+\psi=0,\,\,\,\,\quad\,\textrm{in}\quad\,\,\,\,\mathcal{D},\\
\upsilon,\,\psi>0,\quad\quad\quad\quad\quad\textrm{in}\quad\,\,\,\,\mathcal{D},\\
\frac{\partial }{\partial \nu}\upsilon=\frac{\partial }{\partial \nu}\psi=0,\,\,\,\,\quad\,\textrm{on}\quad\,\,\partial\mathcal{D},
\end{array}\right.
\end{equation}
because the first equation in system (\ref{1.6}) implies
$$
\aligned
\int_{\mathcal{D}}\psi|\nabla(\log \psi-\upsilon)|^2=0,
\endaligned
$$
and hence $\psi=\lambda e^{\upsilon}$ for some positive constant $\lambda$.
Steady states of  system (\ref{1.6}), namely its solutions,
are of basic importance
for a better understanding of global dynamics to the following  Keller-Segel system
with $\tau\geq0$:
\begin{equation}\label{1.9}
\left\{\begin{array}{ll}
\psi_t=\Delta \psi-\nabla(\psi\nabla \upsilon),\,
\,\,\,\,\,\textrm{in}\quad\,\,\,\,\mathcal{D},\\
\tau \upsilon_t=\Delta \upsilon-\upsilon+\psi,\,\,\,\,\quad\,\textrm{in}\quad\,\,\,\,\mathcal{D},\\
\upsilon,\,\psi>0,\quad\quad\quad\quad
\quad\quad\textrm{in}\quad\,\,\,\,\mathcal{D},\\
\frac{\partial }{\partial \nu}\upsilon=\frac{\partial }{\partial \nu}\psi=0,\,\,\,\,
\quad\quad\,\textrm{on}\quad\,\,\partial\mathcal{D},
\end{array}\right.
\end{equation}
which describes chemotactic
feature of cellular slime molds   sensitive to the gradient of a chemical substance secreted by themselves
(see \cite{KS}).
The one-dimensional form of
system (\ref{1.6}) was first studied by Schaaf \cite{S}.
In higher dimensions $N\geq2$
Biler \cite{B} established the existence of
non-constant radially symmetric  solution  to  (\ref{1.6})
when the domain $\mathcal{D}$ is a ball.
In the general two-dimensional case, Wang-Wei \cite{WW}, independently of  Senba-Suzuki
\cite{SS}, proved that for any $\mu\in(0,1/|\mathcal{D}|+\mu_1)\setminus\{4\pi m|m=1,2,\ldots\}$
(where $\mu_1$ denotes the first positive eigenvalue of $-\Delta$ with Neumann boundary condition),
system (\ref{1.6}) has a non-constant solution such that $\int_{\mathcal{D}}\psi=\mu|\mathcal{D}|$.
Meanwhile, if space dimension is $N=2$,
it is known that
as infinite time blow-up solutions of the
parabolic-elliptic system (\ref{1.9}) from chemotaxis,  steady states of  (\ref{1.6})
produce a
significant  concentration phenomenon  in mathematical biology
referred as  `chemotactic collapse', namely
the blow-up for the quantity $\psi$  in (\ref{1.6}) takes place as a finite sum
of Dirac measure at points with masses    equal to
$8\pi$ or $4\pi$, respectively, depending on whether the  blow-up
points lie inside the domain or on the boundary.
By analyzing the asymptotic behavior of families of solutions to equation $(\ref{1.3})|_{p=1}$
Senba-Suzuki \cite{SS,SS2} exhibited this phenomenon  for
the term $\lambda e^\upsilon$ in $(\ref{1.3})|_{p=1}$ with   positive, uniformly
bounded mass  $\lambda\int_\mathcal{D} e^\upsilon$
as $\lambda$ tends to zero. More precisely,  if $\upsilon_\lambda$ is a family of solutions of $(\ref{1.3})$
under $p=1$ and $N=2$,
 such that
$$
\aligned
\lim_{\lambda\rightarrow0}\lambda\int_{\mathcal{D}} e^{\upsilon_\lambda}=L>0,
\endaligned
$$
then there exist non-negative integers $k$, $l$ with $k+l\geq1$
for which $L=4\pi(k+2l)$. Moreover, once $\lambda$ tends to zero,
this family of solutions concentrate at $l$ different points $\xi_1,\ldots,\xi_l$
inside the domain $\mathcal{D}$ and $k$ different points $\xi_{l+1},\ldots,\xi_{k+l}$
on the boundary  $\partial\mathcal{D}$.
In particular, far away from these concentration points the asymptotic
profile of $\upsilon_{\lambda}$ is uniformly described as
$$
\aligned
\upsilon_{\lambda}(x)\rightarrow
\sum\limits_{i=1}^{l}
8\pi G(x,\xi_i)+\sum\limits_{i=l+1}^{k+l}
4\pi G(x,\xi_i).
\endaligned
$$
In addition, these concentration points or blow-up points $\xi=(\xi_1,\ldots,\xi_{k+l})$
are nothing but  critical points
of a functional
\begin{equation}\label{1.4}
\aligned
\varphi_{k+l}(\xi)=\varphi_{k+l}(\xi_1,\ldots,\xi_{k+l})
=
\sum_{i=1}^{k+l}c_i^2H(\xi_i,\xi_i)+\sum\limits_{i\neq j}^{k+l}
c_ic_j G(\xi_i,\xi_j),
\endaligned
\end{equation}
where
$c_i=8\pi$ for $i=1,\ldots,l$, but $c_i=4\pi$
for $i=l+1,\ldots,k+l$, $G(x,y)$  denotes
the
Green's function  of the   problem
$$
\left\{\aligned
&-\Delta_xG(x,y)+G(x,y)=\delta_y(x),\,\,\,\,\,\,\,
x\in\Omega,\\
&\frac{\partial G}{\partial\nu_x}(x,y)=0,
\qquad\,\,
\qquad\qquad\qquad\,\,
x\in\partial\Omega,
\endaligned\right.
$$
and $H(x,y)$  its regular part  defined as
$$
\aligned
H(x,y)=\left\{\aligned
&G(x,y)+\frac{1}{2\pi}\log|x-y|,\,\quad\,y\in\Omega,\\[1mm]
&G(x,y)+\frac1{\pi}\log|x-y|,\,\ \quad\,y\in\po.
\endaligned\right.
\endaligned
$$
Reciprocally, in the spirit
of the  Lyapunov-Schmidt finite-dimensional reduction method
del Pino-Wei \cite{DW}  constructed a family of mixed interior and boundary
bubbling solutions for equation $(\ref{1.3})|_{p=1,N=2}$ with exactly  the asymptotic profile above.
Successively, when  $N=2$ and $p$ is between $0$ and $2$,  Deng \cite{D} used a  reductional argument to build solutions for
equation (\ref{1.3})
with bubbling profiles at points inside $\mathcal{D}$ and on the boundary
$\partial\mathcal{D}$, which recovered the result in \cite{DW} when $p=1$.
In general, such bubbling solutions are called
{\it solutions concentrating on $0$-dimensional sets with uniformly bounded mass}.

Clearly, a natural question is to ask whether there exists a family of solutions
of  equation (\ref{1.3})
concentrating on higher dimensional subsets of $\overline{\mathcal{D}}$ with or without
uniformly bounded mass as the parameter $\lambda$ tends to zero.
The first result in this  direction was obtained by
Pistoia-Vaira \cite{PV} in the case that
$p=1$ and the domain  $\mathcal{D}$ is a ball with dimension $N\geq2$.
Based on a fixed-point argument,  they constructed
 a family of uniformly unbounded mass radial solutions
of  $(\ref{1.3})|_{p=1}$ in the ball, which
blow up on the entire boundary  and hence produce the boundary concentration layer.
Following closely the techniques of  $(\ref{1.3})|_{p=1}$  in \cite{PV}, Bonheure-Casteras-Noris \cite{BCN}
constructed a family of boundary layer solutions  in the annulus
blowing up simultaneously along both boundaries,
a family of  internal  layer solutions in the unit ball
blowing up on an interior sphere,
and a family of   solutions in the unit ball
with  an internal layer and a boundary layer
blowing up simultaneously on an interior sphere and the boundary.
Recently,  when the domain $\mathcal{D}$ is a unit disk
(corresponding to a unit ball with dimension $N=2$),
Bonheure-Casteras-Rom\'{a}n \cite{BCR} have successfully constructed
a family of uniformly unbounded mass radial solutions of  $(\ref{1.3})|_{p=1}$
which concentrate at the origin and blow up on the entire boundary.
Additionally, some bifurcation analyses of radial solutions to $(\ref{1.3})|_{p=1}$
in a ball with dimension $N\geq2$ were also performed by Bonheure et al. in \cite{BCF,BCN1}.
As for a general smooth two-dimensional  domain $\mathcal{D}$,
it is very worth mentioning that inspired by the novel result in \cite{PV},
del Pino-Pistoia-Varia \cite{DPV}   applied an infinite-dimensional form of Lyapunov-Schmidt reduction
to establish the existence of  a family of solutions  $\upsilon_\lambda$
for equation $(\ref{1.3})|_{p=1,N=2}$ with
unbounded mass $\lambda\int_\mathcal{D} e^{\upsilon_\lambda}$,  which exhibit a sharp boundary layer
and blow up along
the entire  $\partial\mathcal{D}$  as  $\lambda$ tends to zero but remains
suitably away from a sequence of critical small values where certain resonance phenomenon occurs.
Finally, when the domain  $\mathcal{D}$
has  suitable rational symmetries in higher dimensions $N\geq3$,
 Agudelo-Pistoia \cite{AP}  constructed  several families  of  layered  solutions  $\upsilon_\lambda$  of the stationary Keller-Segel
chemotaxis equation $(\ref{1.3})|_{p=1}$  with uniformly bounded mass
$\lambda\int_\mathcal{D} e^{\upsilon_\lambda}$, which only exhibit
{\it three } different types of chemoattractant concentration  along
 suitable $(N-2)$-dimensional minimal submanifolds of the boundary.

Problem (\ref{1.1}) is seemingly  similar to equation (\ref{1.3}).
Our original motivation   in equation (\ref{1.3}) is based on  the fact that except
for $p=1$,
nothing  is known about the existence or the boundary concentration phenomenon
for  solutions  of  equation  (\ref{1.3})
in higher dimensions $N\geq3$.
For this aim  our idea is to  consider partially axially symmetric solutions of  equation (\ref{1.3})
when  the domain $\mathcal{D}$ has some rotational symmetries, which implies that
problem (\ref{1.1}) can be viewed as a special case of equation (\ref{1.3}) in higher dimensions $N\geq3$.
Indeed,
take $n\in\{1,2\}$ as a fixed integer.
Let $\Omega$  be  a smooth bounded domain in $\mathbb{R}^2$ such that
\begin{equation*}\label{}
\aligned
\overline{\Omega}\subset
\{
(x_1,x_n,x')\in\mathbb{R}^n\times\mathbb{R}^{2-n}|\,
\,\,x_i>0,\,\,\,i=1,n
\}.
\endaligned
\end{equation*}
Fix $k_1,k_n\in\mathbb{N}$ with $k_1+k_n=N-2\geq1$
and set
\begin{equation*}\label{}
\aligned
\mathcal{D}:=
\big\{
(y_1,y_n,x')\in\mathbb{R}^{k_1+1}\times\mathbb{R}^{k_n+1}\times\mathbb{R}^{2-n}|\,
\,(|y_1|,|y_n|,x')\in\Omega
\big\}.
\endaligned
\end{equation*}
Then $\mathcal{D}$ is a smooth bounded domain in $\mathbb{R}^N$
which is  invariant under the action of the group $\Upsilon:=\mathcal{O}(k_1+1)\times\mathcal{O}(k_n+1)$
on $R^N$ given by
$$
\aligned
(g_1,g_n)(y_1,y_n,x'):=
(g_1y_1,g_ny_n,x').
\endaligned
$$
Note that $\mathcal{O}(k_i+1)$ is the group of linear isometries of $\mathbb{R}^{k_i+1}$,
$\mathbb{S}^{k_i}$ is the unit sphere in $\mathbb{R}^{k_i+1}$ and $|\mathbb{S}^{k_i}|$  represents  its area.
If  we seek $\Upsilon$-invariant solutions of equation (\ref{1.3}),
i.e. solutions $\upsilon$ of the form
\begin{equation*}\label{1.10}
\aligned
\upsilon(y_1,y_n,x')=u(|y_1|,|y_n|,x'),
\endaligned
\end{equation*}
a direct calculus shows that equation (\ref{1.3}) is transformed to
\begin{equation}\label{1.11}
\left\{\aligned
&-\Delta u-\sum_{i=1}^n\frac{k_i}{\,x_i\,}\frac{\partial u}{\partial x_i}+u=\lambda u^{p-1}e^{u^p},
\,\,\,\,\,
u>0\,\,\,\,\,\,
\textrm{in}\,\,\,\,\,
\Omega,\\
&\frac{\partial u}{\partial\nu}=0\,\,\,\,\,
\ \,\ \,\qquad\qquad\qquad\qquad\,\,\,\,
\,\ \,\,\quad\qquad\qquad\,\,\
\textrm{on}\,\,\,\,
\partial\Omega.
\endaligned\right.
\end{equation}
Clearly,   if we take anisotropic coefficient
\begin{equation}\label{1.12}
\aligned
a(x)=a(x_1,x_n,x'):=x_1^{k_1}\cdot x_n^{k_n},
\endaligned
\end{equation}
then equation (\ref{1.11}) (or (\ref{1.3})) can be rewritten as problem (\ref{1.1}).
Moreover, the  solution $\upsilon$ of equation (\ref{1.3})
always has the same energy
$\beta_\lambda$ as the corresponding solution $u$ of problem (\ref{1.1}), namely
\begin{equation*}
\aligned
\beta_\lambda:=\frac{\lambda p}{2}\frac{1}{|\mathbb{S}^{k_1}|}
\frac{1}{|\mathbb{S}^{k_n}|}
\left(
\int_{\mathcal{D}}
\big(e^{\upsilon^p}-1\big)dy
\right)^{\frac{2-p}{p}}
\left(
\int_{\mathcal{D}}
\upsilon^p
e^{\upsilon^p}dy
\right)^{\frac{2(p-1)}{p}}
=
\frac{\lambda p}{2}
\left(
\int_{\Omega}
a(x)\big(e^{u^p}-1\big)dx
\right)^{\frac{2-p}{p}}
\left(
\int_{\Omega}
a(x)u^p
e^{u^p}dx
\right)^{\frac{2(p-1)}{p}}.
\endaligned
\end{equation*}
Hence  by considering rotational symmetry of $\mathcal{D}$,
a fruitful  approach for seeking layered  solutions
of equation (\ref{1.3}) with concentration
along some $(N-2)$-dimensional minimal submanifolds of $\overline{\mathcal{D}}$
diffeomorphic to
$\mathbb{S}^{k_1}\times\mathbb{S}^{k_n}$
is to reduce it to produce  pointwise  blow-up solutions of the
anisotropic problem (\ref{1.1}) in the domain
$\Omega$ of dimension $2$. This approach, together
with some Lyapunov-Schmidt finite-dimensional reduction arguments, has recently been taken
to construct multi-layer positive
solutions of equation (\ref{1.3})  concentrating
along some $(N-2)$-dimensional minimal submanifolds of $\partial\mathcal{D}$,
which can be  found in \cite{AP}  only for the case $p=1$.

In this paper, our goal  is to obtain the existence of
boundary separated or clustered  layer positive solutions
for equation (\ref{1.3}) in the higher-dimensional domain with some
rotational symmetries,
by constructing bubbling solutions  for  the anisotropic planar problem (\ref{1.1})
with  simple or non-simple boundary  concentration points
when $p$ is between $0$ and $2$.
We try to use a new reductional argument to investigate the effect of anisotropic coefficient
$a(x)$  on the existence
of boundary concentrating solutions to problem (\ref{1.1}).
As a result,
with the help of some suitable assumptions on anisotropic
coefficient $a(x)$ we prove that there exists a
family of positive solutions $u_\lambda$ of problem (\ref{1.1})
with an arbitrary number of
mixed interior and boundary bubbles
which concentrate at  totally
different strict local
maximum or minimal boundary
points  of $a(x)$ restricted to $\po$, or accumulate to
the same strict local maximum boundary
 point of $a(x)$ over $\oo$ as  $\lambda$ tends to zero.
Furthermore, for these bubbling solutions $u_\lambda$ we compute the delicate expansion of the corresponding energy
\begin{equation}\label{1.26}
\aligned
\beta_\lambda
=\,
\frac{\lambda p}{2}
\left(
\int_{\Omega}
a(x)\big(e^{u_\lambda^p}-1\big)dx
\right)^{\frac{2-p}{p}}
\left(
\int_{\Omega}
a(x)u_\lambda^p
e^{u_\lambda^p}dx
\right)^{\frac{2(p-1)}{p}}
\endaligned
\end{equation}
and  exhibit the sharp difference of tendency of $\beta_\lambda$ between $0<p<1$ and $1<p<2$.
In particular, we recover and improve those results for $p=1$
  in \cite{AP} and for $a(x)\equiv1$ in \cite{D}.

Before  precisely stating our  results, let us start with some notations.
Let $\gamma$ and $\varepsilon$ be positive parameters  chosen  by the relations
\begin{equation}\label{1.5}
\aligned
p\lambda\gamma^{2(p-1)}\varepsilon^2e^{\gamma^p}=1,
\endaligned
\end{equation}
and
\begin{equation}\label{2.7}
\aligned
p\gamma^{p}=-4\log\varepsilon,
\endaligned
\end{equation}
where
$\lambda\rightarrow0$ if and only if $\gamma\rightarrow+\infty$
and  $\varepsilon\rightarrow0$,
especially  $\lambda=\varepsilon^2$
and $\lambda^2e^{\gamma}=1$ if
$p=1$.
Let
$$
\aligned
\Delta_au=\frac1{a(x)}\nabla(a(x)\nabla u)=\Delta u+\nabla\log a(x)\nabla u,
\endaligned
$$
and $G_a(x,y)$ be the anisotropic Green's function associated to the Neumann  equation
\begin{equation}\label{1.7}
\left\{\aligned
&-\Delta_aG_a(x,y)+G_a(x,y)=\delta_y(x),\,\,\,\,\,\,\,
x\in\Omega,\\
&\frac{\partial G_a}{\partial\nu_x}(x,y)=0,
\qquad\,\,
\qquad\qquad\qquad\ \,\,\,
x\in\partial\Omega,
\endaligned\right.
\end{equation}
for every $y\in\oo$.
The regular part of $G_a(x,y)$ is defined depending on whether $y$ lies inside
the domain or on its boundary as
\begin{equation}\label{1.8}
\aligned
H_a(x,y)=\left\{\aligned
&G_a(x,y)+\frac{1}{2\pi}\log|x-y|,\,\quad\,y\in\Omega,\\
&G_a(x,y)+\frac1{\pi}\log|x-y|,\,\ \quad\,y\in\po.
\endaligned\right.
\endaligned
\end{equation}
In this way, $y\in\oo\mapsto H_a(\cdot,y)\in C\big(\Omega, C^{\alpha}(\overline{\Omega})\big)
\cap C\big(\partial\Omega, C^{\alpha}(\overline{\Omega})\big)$ and
$H_a(x,y)\in C^\alpha\big(\oo\times\Omega\big)\cap C^\alpha\big(\oo\times\po\big)
\cap C^1\big(\oo\times\Omega\setminus\{x=y\}\big)
\cap C^1\big(\oo\times\po\setminus\{x=y\}\big)$
for any $\alpha\in(0,1)$,
and the corresponding Robin's function $y\in\oo\mapsto H_a(y,y)$ belongs to $C^1(\Omega)\cap C^1(\partial\Omega)$
(see \cite{AP}). Moreover, by the maximum principle, for any $y\in\overline{\Omega}$, $G_a(\cdot,y)>0$
over $\overline{\Omega}$.

Our  first result concerns the existence of
solutions of problem (\ref{1.1})
whose mixed interior and boundary bubbles are uniformly far away from each other
and interior bubbles
lie in the domain with  distance
to the boundary uniformly approaching zero.

\vspace {1mm}
\vspace {1mm}
\vspace {1mm}
\vspace {1mm}

\noindent{\bf Theorem 1.1.}\,\,\,{\it
Let $k$, $l$  be any non-negative integers with $k+l\geq1$,
$0<p<2$ and assume that there exist
$k+l$ different points  $\xi^*_1,\ldots,\xi^*_{k+l}\in\partial\Omega$
such that
each $\xi_i^*$ is either a strict local
maximum or a strict local minimum point of $a(x)$ restricted to $\po$  and
satisfies  for all $i=1,\ldots,l$,
$\partial_{\nu}a(\xi_i^*):=\langle\nabla a(\xi_i^*),\,\nu(\xi_i^*)\rangle>0$.
Then for any sufficiently small $\lambda$,
there exists  a
family of  positive solutions $u_\lambda$
for problem {\upshape (\ref{1.1})}
with $k$ different boundary bubbles and $l$ different interior bubbles
located at distance $O\left(1/|\log\varepsilon|\right)$ from $\po$
such that, as $\lambda\rightarrow0$,
%$$
%\aligned
%\lim\limits_{\lambda\rightarrow0}\varepsilon^{\frac{2(2-p)}{p}}\int_{\Omega}a(x)e^{(u_\lambda)^p}=
%\sum_{i=1}^{k+l}c_ia(\xi_i^*),
%\endaligned
%$$
\begin{eqnarray}\label{1.19}
p\gamma^{p-1}
u_{\lambda}(x)\rightarrow
\sum\limits_{i=1}^{k+l}
c_i G_a(x,\xi_i^*)
\,\quad\textrm{in}\,\,\,\,
C_{loc}\big(\overline{\Omega}\setminus\{\xi^*_1,\ldots,\xi^*_{k+l}\}\big),
\end{eqnarray}
and
\begin{equation}\label{1.14}
\aligned
p\gamma^{p-1}\lambda u_\lambda^{p-1}e^{u_\lambda^p}
\rightharpoonup
\sum_{i=1}^{k+l}c_i\delta_{\xi_i^*}
\qquad
\textrm{weakly in the sense of measure in}
\,\,\,
\overline{\Omega},
\endaligned
\end{equation}
and
%$$
%\aligned
%u_\lambda(x)=\frac1{p\gamma^{p-1}}
%\sum\limits_{i=1}^{k+l}\left[\,\log
%\frac1{(\varepsilon^2\mu_i^2+|x-\xi_i^\varepsilon|^2)^2}
%+c_iH_a(x,\xi_i^\varepsilon)
%+o(1)\,\right],
%\endaligned
%$$
\begin{eqnarray}\label{1.18}
u_\lambda(x)=\frac1{p\gamma^{p-1}}
\sum\limits_{i=1}^{k+l}\left\{\,
\log
\frac1{\,((\varepsilon\mu^{\varepsilon}_i)^2+|x-\xi_i^\varepsilon|^2)^2\,}
+c_i H_a(x,\xi_i^\varepsilon)
+\sum_{j=1}^4\left(\frac{p-1}{p}\right)^j\frac{1}{\gamma^{jp}}
\left[\omega^j_{\mu^{\varepsilon}_i}\left(\frac{x-\xi_i^\varepsilon}{\varepsilon}\right)
\right.\right.
\nonumber\\[1mm]
\left.\left.
+D^j_{\mu^{\varepsilon}_i}\log(\varepsilon\mu^{\varepsilon}_i)
-
\frac{1}{4} D^j_{\mu^{\varepsilon}_i}
c_i H_a(x,\xi_i^\varepsilon)
\right]
+O\left(\frac{1}{|\log\varepsilon|^3}\right)\right\},
\quad\qquad\qquad\qquad
\qquad\qquad\qquad\qquad
\end{eqnarray}
uniformly in $\oo$, where
$\gamma$ and $\varepsilon$  are defined in  {\upshape(\ref{1.5})}-{\upshape(\ref{2.7})},
$c_i=8\pi$ for $i=1,\ldots,l$, but $c_i=4\pi$
for $i=l+1,\ldots,k+l$, $\omega^j_{\mu^{\varepsilon}_i}$ and
$D^j_{\mu^{\varepsilon}_i}$ are defined in {\upshape(\ref{2.9})} and {\upshape(\ref{2.14})}, respectively,
the parameter  $\mu_i^\varepsilon$ satisfies
$$
\aligned
\frac1C\leq\mu_i^\varepsilon\leq|\log\varepsilon|^C,
\endaligned
$$
for some $C>0$ independent of $\varepsilon$, and
 $(\xi^\varepsilon_1,\ldots,\xi^\varepsilon_{k+l})\in\Omega^l\times(\po)^{k}$ satisfies
$$
\aligned
\xi^\varepsilon_i\rightarrow\xi^*_i
\,\quad\,\textrm{for all}\,\,\,i,
\,\qquad\quad\,
\textrm{and}
\,\qquad\quad\,
\dist(\xi^\varepsilon_i,\po)=O\left(1/|\log\varepsilon|\right)
\quad\,\,\forall\,\,i=1,\ldots,l.
\endaligned
$$
Furthermore, the corresponding energy defined in  {\upshape(\ref{1.26})} satisfies
%\begin{equation}\label{}
%\aligned
%\beta_\lambda
%=
%\frac{\lambda p}{2}
%\left(
%\int_{\Omega}
%a(x)\big(e^{u_\lambda^p}-1\big)dx
%\right)^{\frac{2-p}{p}}
%\left(
%\int_{\Omega}
%a(x)u_\lambda^p
%e^{u_\lambda^p}dx
%\right)^{\frac{2(p-1)}{p}}
%=
%\frac{1}{2}
%\left(\sum\limits_{i=1}^{k+l}
%c_i a(\xi_i^\varepsilon)\right)
%\left[
%1
%+
%O\left(\frac{\log^2|\log\varepsilon|}{|\log\varepsilon|^2}\right)
%\right],
%\endaligned
%\end{equation}
%\begin{eqnarray}\label{1.15}
%\beta_\lambda
%=\,
%\frac{\lambda p}{2}
%\left(
%\int_{\Omega}
%a(x)\big(e^{u_\lambda^p}-1\big)dx
%\right)^{\frac{2-p}{p}}
%\left(
%\int_{\Omega}
%a(x)u_\lambda^p
%e^{u_\lambda^p}dx
%\right)^{\frac{2(p-1)}{p}}
%\ \,\ \qquad\qquad\qquad
%\nonumber\\[1mm]
%=\,
%\frac{1}{2}
%\left(\sum\limits_{i=1}^{k+l}
%c_i a(\xi_i^\varepsilon)\right)
%\left[
%1
%+
%O\left(\frac{\log^2|\log\varepsilon|}{|\log\varepsilon|^2}\right)
%\right]
%\,
%\rightarrow
%\,
%\frac{1}{2}
%\sum\limits_{i=1}^{k+l}
%c_i a(\xi_i^*)
%\qquad
%\textrm{as}\,\,\,
%\lambda\rightarrow0,
%\end{eqnarray}
\begin{eqnarray}\label{1.15}
\beta_\lambda
=\,
\frac{1}{2}
\left(\sum\limits_{i=1}^{k+l}
c_i a(\xi_i^\varepsilon)\right)
\left[
1
+
O\left(\frac{\log^2|\log\varepsilon|}{|\log\varepsilon|^2}\right)
\right]
\,
\rightarrow
\,
\frac{1}{2}
\sum\limits_{i=1}^{k+l}
c_i a(\xi_i^*)
\qquad
\textrm{as}\,\,\,
\lambda\rightarrow0,
\end{eqnarray}
and if  $0<p<1$,
\begin{equation}\label{1.16}
\aligned
\beta_\lambda\leq
\frac{1}{2}
\left(\sum\limits_{i=1}^{k+l}
c_i a(\xi_i^\varepsilon)\right)
\left\{1+
\frac{4(p-1)}{\,p^2\gamma^{2p}\,}
\left[1+
O\left(\frac{\log^D|\log\varepsilon|}{|\log\varepsilon|}\right)
\right]\right\}
<\,\frac{1}{2}
\sum\limits_{i=1}^{k+l}
c_i a(\xi_i^\varepsilon),
\endaligned
\end{equation}
but if $1<p<2$,
\begin{equation}\label{1.17}
\aligned
\beta_\lambda\geq
\frac{1}{2}
\left(\sum\limits_{i=1}^{k+l}
c_i a(\xi_i^\varepsilon)\right)
\left\{1+
\frac{4(p-1)}{\,p^2\gamma^{2p}\,}
\left[1+
O\left(\frac{\log^D|\log\varepsilon|}{|\log\varepsilon|}\right)
\right]\right\}
>\frac{1}{2}
\sum\limits_{i=1}^{k+l}
c_i a(\xi_i^\varepsilon),
\endaligned
\end{equation}
for some $D>1$, independent of $\varepsilon$.
}

\vspace{1mm}
\vspace{1mm}
\vspace{1mm}
\vspace{1mm}

Our  next  result
concerns the existence of solutions of problem (\ref{1.1})
with mixed interior and boundary bubbles which  accumulate to
the same  boundary point.

\vspace{1mm}
\vspace{1mm}
\vspace{1mm}
\vspace{1mm}

\noindent{\bf Theorem 1.2.}\,\,\,{\it
Let $k$, $l$ be any non-negative integers with $k+l\geq1$,
$0<p<2$
and assume that $\xi_*\in\po$ is a strict local maximum point of $a(x)$
over $\overline{\Omega}$
and satisfies $\partial_{\nu}a(\xi_*):=\langle\nabla a(\xi_*),\,\nu(\xi_*)\rangle=0$.
Then for any sufficiently small $\lambda$,
there exists a
family of positive solutions $u_\lambda$ for
problem {\upshape (\ref{1.1})}
with $k$ different boundary bubbles and $l$ different interior bubbles
which accumulate to $\xi_*$ as $\lambda\rightarrow0$. More precisely, as $\lambda\rightarrow0$,
\begin{eqnarray}\label{1.20}
p\gamma^{p-1}
u_{\lambda}(x)\rightarrow
4\pi(k+2l)G_a(x,\xi_*)
\,\quad\textrm{in}\,\,\,\,
C_{loc}\big(\overline{\Omega}\setminus\{\xi_*\}\big),
\end{eqnarray}
and
\begin{equation}\label{1.21}
\aligned
p\gamma^{p-1}\lambda u_\lambda^{p-1}e^{u_\lambda^p}
\rightharpoonup
4\pi(k+2l)\delta_{\xi_*}
\qquad
\textrm{weakly in the sense of measure in}
\,\,\,
\overline{\Omega},
\endaligned
\end{equation}
and
%$$
%\aligned
%u_\lambda(x)=\frac1{p\gamma^{p-1}}
%\sum\limits_{i=1}^{k+l}\left[\,\log
%\frac1{(\varepsilon^2\mu_i^2+|x-\xi_i^\varepsilon|^2)^2}
%+c_iH_a(x,\xi_i^\varepsilon)
%+o(1)\,\right],
%\endaligned
%$$
\begin{eqnarray}\label{1.22}
u_\lambda(x)=\frac1{p\gamma^{p-1}}
\sum\limits_{i=1}^{k+l}\left\{\,
\log
\frac1{((\varepsilon\mu^{\varepsilon}_i)^2+|x-\xi_i^\varepsilon|^2)^2}
+c_i H_a(x,\xi_i^\varepsilon)
+\sum_{j=1}^4\left(\frac{p-1}{p}\right)^j\frac{1}{\gamma^{jp}}
\left[\omega^j_{\mu^{\varepsilon}_i}\left(\frac{x-\xi_i^\varepsilon}{\varepsilon}\right)
\right.\right.
\nonumber\\[1mm]
\left.\left.
+\,D^j_{\mu^{\varepsilon}_i}\log(\varepsilon\mu^{\varepsilon}_i)
- \frac{1}{4} D^j_{\mu^{\varepsilon}_i}
c_i H_a(x,\xi_i^\varepsilon)
\right]
+O\left(\frac{1}{|\log\varepsilon|^3}\right)\right\},
\quad\qquad\qquad\qquad\qquad
\qquad\qquad\quad\ \,
\end{eqnarray}
uniformly in $\oo$, where
$\gamma$ and $\varepsilon$  are defined in  {\upshape(\ref{1.5})}-{\upshape(\ref{2.7})},
$c_i=8\pi$ for $i=1,\ldots,l$, but $c_i=4\pi$
for $i=l+1,\ldots,k+l$, $\omega^j_{\mu^{\varepsilon}_i}$ and
$D^j_{\mu^{\varepsilon}_i}$ are defined in {\upshape(\ref{2.9})} and {\upshape(\ref{2.14})}, respectively,
the parameter  $\mu_i^\varepsilon$ satisfies
$$
\aligned
\frac1C\leq\mu_i^\varepsilon\leq|\log\varepsilon|^C,
\endaligned
$$
for some $C>0$ independent of $\varepsilon$, and
 $(\xi^\varepsilon_1,\ldots,\xi^\varepsilon_{k+l})\in\Omega^l\times(\po)^{k}$ satisfies
$$
\aligned
\xi^\varepsilon_i\rightarrow\xi_*
\,\quad\,\forall\,\,i,
\qquad
|\xi^\varepsilon_i-\xi^\varepsilon_j|>\frac{1}{\,|\log\varepsilon|^{2(m^2+1)}}
\quad\forall\,\,i\neq j,
\qquad
\textrm{and}
\qquad
\dist(\xi^\varepsilon_i,\po)>\frac{1}{\,|\log\varepsilon|^{2(m^2+1)}}
\quad\forall\,\,i=1,\ldots,l.
\endaligned
$$
Furthermore, the corresponding energy defined in  {\upshape(\ref{1.26})} satisfies
\begin{eqnarray}\label{1.23}
\beta_\lambda
=\,
\frac{1}{2}
\left(\sum\limits_{i=1}^{k+l}
c_i a(\xi_i^\varepsilon)\right)
\left[
1
+
O\left(\frac{\log^2|\log\varepsilon|}{|\log\varepsilon|^2}\right)
\right]
\,
\rightarrow
\,
2\pi(k+2l)a(\xi_*)
\qquad
\textrm{as}\,\,\,
\lambda\rightarrow0,
\end{eqnarray}
and if  $0<p<1$,
\begin{equation}\label{1.24}
\aligned
\beta_\lambda\leq
\frac{1}{2}
\left(\sum\limits_{i=1}^{k+l}
c_i a(\xi_i^\varepsilon)\right)
\left\{1+
\frac{4(p-1)}{\,p^2\gamma^{2p}\,}
\left[1+
O\left(\frac{\log^D|\log\varepsilon|}{|\log\varepsilon|}\right)
\right]\right\}
<\,\frac{1}{2}
\sum\limits_{i=1}^{k+l}
c_i a(\xi_i^\varepsilon)
<\,
2\pi(k+2l)a(\xi_*),
\endaligned
\end{equation}
but if $1<p<2$,
\begin{equation}\label{1.25}
\aligned
\beta_\lambda\geq
\frac{1}{2}
\left(\sum\limits_{i=1}^{k+l}
c_i a(\xi_i^\varepsilon)\right)
\left\{1+
\frac{4(p-1)}{\,p^2\gamma^{2p}\,}
\left[1+
O\left(\frac{\log^D|\log\varepsilon|}{|\log\varepsilon|}\right)
\right]\right\}
>\frac{1}{2}
\sum\limits_{i=1}^{k+l}
c_i a(\xi_i^\varepsilon),
\endaligned
\end{equation}
for some $D>1$, independent of $\varepsilon$.
}

\vspace{1mm}
\vspace{1mm}
\vspace{1mm}
\vspace{1mm}

From Theorem 1.2 we derive  that  if $0<p<2$ and  the domain  $\mathcal{D}$
has some rational symmetries in higher dimensions
$N\geq3$ such that the corresponding anisotropic coefficient
$a(x)$ given by (\ref{1.12}) satisfies the assumptions in Theorem 1.2, then
equation (\ref{1.3})  has a family of positive solutions with
arbitrarily many   mixed interior and boundary layers which
collapse to the same  $(N-2)$-dimensional minimal submanifold
of  $\partial\mathcal{D}$ as $\lambda$ tends to zero. Meanwhile,
 we  observe  that the assumptions  in
Theorem $1.2$  contain the following two cases:\\
\indent (C1)
$\xi_*\in\po$ is a strict
local maximum point of $a(x)$ restricted to $\po$;\\
\indent (C2)
$\xi_*\in\po$ is a strict local maximum
point of $a(x)$ restricted in $\Omega$
and satisfies $\partial_{\nu}a(\xi_*)=\langle\nabla a(\xi_*),\,\nu(\xi_*)\rangle=0$.\\
Arguing as in the proof of Theorem $1.2$, we readily prove that if
(C1) holds,  then problem (\ref{1.1}) has positive solutions with
arbitrarily many
boundary bubbles which  accumulate to $\xi_*$ along $\po$; while
if (C2) holds,  then problem (\ref{1.1}) has positive solutions with
arbitrarily many
interior bubbles which  accumulate to $\xi_*$ along the neighborhood  near the
inner normal direction of $\po$. As  for   the latter  case,  our result seems to close  some \textit{gap} which
was left open in the  literature \cite{AP} regarding
such type of chemoattractant concentration
from the stationary  Keller-Segel system with linear
chemotactical
sensitivity function, namely  involving the existence of
solutions of  equation $(\ref{1.3})|_{p=1,N\geq3}$ with
an arbitrary number of interior layers
which simultaneously
accumulate along a suitable $(N-2)$-dimensional minimal submanifold
of $\partial\mathcal{D}$  as $\lambda$ tends to zero.
Moreover,
 this type of boundary concentration phenomenon has also appeared in \cite{Z}
 for  positive solutions of
a higher-dimensional   elliptic Neumann problem with large exponent.
Finally, it is necessary to point out  that
radial  solutions  of  equation $(\ref{1.3})|_{p=1}$ with concentration on an arbitrary number
of  internal  spheres  were  built by Bonheure-Casteras-Noris \cite{BCN1} when the domain  $\mathcal{D}$
is a  ball with dimension $N\geq2$, but
a remarkable fact is that, in opposition to our result
or an analogous  one given by  Malchiodi-Ni-Wei \cite{MNW} for
a singularly  perturbed elliptic Neumann  problem on a ball, the layers of those  solutions do not accumulate
to the boundary of $\mathcal{D}$ as $\lambda$ tends to zero.

The proof of  our main results  relies on a
well known  Lyapunov-Schmidt finite-dimensional reduction procedure in which the fundamental problem is how to attain  bounded invertibility of the linearized operator
$\mathcal{L}$ defined in (\ref{2.30}).
The same strategy
 has been used by Deng \cite{D} to construct mixed interior and boundary separated
 bubbling solutions for the
 Neumann problem
\begin{equation}\label{1.13}
\left\{\aligned
&-\Delta u+u=\lambda  u^{p-1}e^{u^p},\,\,\,\,
u>0\,\,\ \,\,\,
\textrm{in}\,\,\,\,\,
\Omega,\\[1mm]
&\frac{\partial u}{\partial\nu}=0\,\,
\qquad\quad\qquad\qquad
\ \ \ \ \,\qquad\,
\textrm{on}\,\,\,
\partial\Omega,
\endaligned\right.
\end{equation}
where  $\Omega\subset\mathbb{R}^2$ is a smooth  bounded   domain,
$\lambda>0$ is a small parameter,
$0<p<2$ and $\nu$
denotes the outer unit normal vector to $\po$.
However,  in contrast with
the result in \cite{D},
the location of mixed interior and boundary
 separated or clustered
bubbles in our present work
is not been characterized as an absolute minimum point of $\varphi_{k+l}(\xi)$ defined in (\ref{1.4}),
but some totally different
strict local  maximum or minimum boundary points of $a(x)$ restricted to $\po$,
or a certain strict local maximum boundary point of $a(x)$ over $\oo$,
which needs us to investigate deeply the effect of anisotropic coefficient
$a(x)$  on the existence of boundary bubbling solutions
in the configuration space $\mathcal{O}_\varepsilon$ of
 separated and clustered concentration points defined by (\ref{2.4}).
For this aim,
in order to make the error  $E_{\xi'}$ of the approximation
of the scaling equation (\ref{2.26})  still small enough under the
occurrence of the accumulation of concentration points, unlike \cite{D} we
build the initial approximate form of up to fourth  order for the solution of problem (\ref{1.1})
that not only has the original cell $\omega_{\varepsilon,\mu_i,\xi_i}$, but also has the ingredients
 $\omega^j_{\mu_{i}}$,
$j=1,2,3,4$ (see (\ref{2.8})). Furthermore,   from the choice of  the concentration parameter $\mu_i$
in system (\ref{2.22}) we can easily construct a good approximation for the solution of the Neumann problem (\ref{1.1})
near each concentration point $\xi_i$. As a consequence, we estimate  the error term $E_{\xi'}$ as $O\left(1/|\log\varepsilon|^4\right)$
under a $L^{\infty}$-norm  $\|\cdot\|_{*}$ with the new weighted function (\ref{2.32})
and compute the Taylor expansion (\ref{3.2}) for the potential $W_{\xi'}$ of equation (\ref{2.26})
(or (\ref{2.29})) near  each scaling point $\xi'_i$ which has a more precise  form than the corresponding  expansion
of $W_{\xi'}$ in Page 739 of \cite{D} because expansion (\ref{3.2}) of  $W_{\xi'}$ does not only look like
 $e^{\omega_{\mu_i}}$, but  also  relates to  the ingredient  $\omega^1_{\mu_{i}}$.
Here, when  we  thoroughly check the validity of the analysis of bounded invertibility of the linearized operation  $\mathcal{L}$
stated in the whole proof of Proposition 4.1 of \cite{D}, it is unfortunate  that
by recalling the corresponding  expansion
of $W_{\xi'}$ in Page 739 of \cite{D}
we don't  obtain the final equality in Page 753 of \cite{D} so that we never conclude bounded invertibility of  $\mathcal{L}$
stated in Proposition 4.1 of \cite{D}.
Notice that $\omega^1_{\mu_{i}}$ has also appeared in expansion (\ref{3.2}) of  $W_{\xi'}$.
So, very differently from the analysis of bounded invertibility of the linearized operation  $\mathcal{L}$
stated
in the proof of Proposition 4.1 of \cite{D}, in this article we should try to apply  a
 new finite-dimensional reduction technique regarding the key ingredient $\omega^1_{\mu_{i}}$
to  obtain bounded invertibility of  $\mathcal{L}$
although we have to confront the essential difficulty from the appearance of   anisotropic coefficient $a(x)$
and mixed interior and boundary separated  or   clustered bubbling points
$\xi_i$, $i=1,\ldots,k+l$.
This is the delicate description during we carry out the whole finite-dimensional  reduction procedure to construct boundary
separated and clustered bubbling  solutions
of problem (\ref{1.1}). Please the readers refer to the the whole proof of Proposition 3.2 in this article.
Finally, let us emphasize  that our reductional arguments regarding  $\omega^1_{\mu_{i}}$ actually
fill  the giant \textit{gap } from the analysis of bounded invertibility of the linearized operator
$\mathcal{L}$ in the proof of Proposition 4.1 of \cite{D}
because  these arguments
will be also suitable for those easy cases  in \cite{D} that anisotropic coefficient $a(x)$ is vanishing, i.e. $a(x)\equiv1$,
 mixed interior and boundary bubbling points $\xi_i$, $i=1,\ldots,k+l$
are uniformly separated from each other
and
interior bubbling points $\xi_i$, $i=1,\ldots,l$
are uniformly far away from
the boundary $\po$.

The proof of expansions
(\ref{1.15})-(\ref{1.17}) and (\ref{1.23})-(\ref{1.25}) of  energy $\beta_\lambda$ of  $k+l$-bubble solutions
described by Theorems 1.1-1.2
relies on the
vanishing identity (\ref{8.18}) of first order and
the identity (\ref{9.10}) of
second order given by \cite{Z1}.
For  $p\in(0,2)$, the second order expansions  (\ref{1.15}) and (\ref{1.23})
of energy $\beta_\lambda$
of  these $k+l$-bubble
solutions can be proven by the second order  Taylor expansion and
the vanishing identity
(\ref{8.18}) of first order. Furthermore,
for  $k+l=1$, with the help of  the identity
(\ref{9.10}) of second order we always observe
%\begin{eqnarray}\label{6.14a}
%\beta_\lambda
%=4\pi
%\left\{1+
%\frac{4(p-1)}{\,p^2\gamma^{2p}\,}
%\left[1+
%O\left(\frac{1}{|\log\varepsilon|}\right)
%\right]\right\}.
%\end{eqnarray}
\begin{equation*}\label{}
\aligned
\beta_\lambda
=
\frac{1}{2}
c_1 a(\xi_1^\varepsilon)
\left\{1+
\frac{4}{\,p^2\gamma^{2p}\,}
\left[p-1+
O\left(\frac{\log^D|\log\varepsilon|}{|\log\varepsilon|}\right)
\right]\right\}
\,
\rightarrow
\,
\frac{1}{2}
c_1 a(\xi_1^*)
\qquad
\textrm{as}\,\,\,
\lambda\rightarrow0,
\endaligned
\end{equation*}
where
$\xi_1^\varepsilon\rightarrow\xi_1^*\in\po$
if $\xi_1^*$ satisfies the assumptions in Theorem 1.1 or 1.2,
$c_1=8\pi$ if $\xi_1^\varepsilon\in\Omega$, but
$c_1=4\pi$ if  $\xi_1^\varepsilon\in\po$.
For general  $k+l\geq1$, by H\"{o}lder's inequality for vectors in $\mathbb{R}_{+}^m$
and  two key  identities   (\ref{8.18})  and (\ref{9.10}) of  first order and second order
 we conclude expansions
(\ref{1.16})-(\ref{1.17}) and (\ref{1.24})-(\ref{1.25})
and hence observe the sharp difference of tendency of energy $\beta_\lambda$ between $0<p<1$ and $1<p<2$.
Here, it is necessary to point out that when anisotropic coefficient $a(x)$ is vanishing, i.e. $a(x)\equiv1$,
 for any $k+l$-bubble solution $u_\lambda$ of problem (\ref{1.13})
that is constructed  by \cite{D} and
concentrates at $l$ different interior points $\xi_1,\ldots,\xi_l$
and $k$ different boundary points $\xi_{l+1},\ldots,\xi_{k+l}$,
we can apply
analogous arguments in
(\ref{1.15})-(\ref{1.17})
%the
%vanishing identity (\ref{8.18}) of first order and
%the identity (\ref{9.10}) of
%second order
to compute the second order Taylor expansion of the corresponding energy
\begin{equation*}
\aligned
\beta_\lambda
=
\frac{\lambda p}{2}
\left(
\int_{\Omega}
\big(e^{u_\lambda^p}-1\big)dx
\right)^{\frac{2-p}{p}}
\left(
\int_{\Omega}
u_\lambda^p
e^{u_\lambda^p}dx
\right)^{\frac{2(p-1)}{p}}.
\endaligned
\end{equation*}
More precisely,  for $k+l=1$,
\begin{eqnarray}\label{1.27}
\beta_\lambda
=2\pi(k+2l)
\left\{1+
\frac{4}{\,p^2\gamma^{2p}\,}
\left[
p-1
+
O\left(\frac{1}{|\log\varepsilon|}\right)
\right]\right\}
\,
\rightarrow
\,
2\pi(k+2l)
\qquad
\textrm{as}\,\,\,
\lambda\rightarrow0,
\end{eqnarray}
where $k=0$ and $l=1$ if $\xi_1\in\Omega$, but
 $k=1$ and $l=0$ if  $\xi_1\in\po$. For general  $k+l\geq1$,
\begin{equation}\label{1.28}
\aligned
\beta_\lambda
=
2\pi(k+2l)
+
O\left(\frac{1}{|\log\varepsilon|^2}\right)
\,
\rightarrow
\,
2\pi(k+2l)
\qquad
\textrm{as}\,\,\,
\lambda\rightarrow0,
\endaligned
\end{equation}
and if $0<p<1$,
\begin{eqnarray}\label{1.29}
\beta_\lambda\leq
2\pi(k+2l)
\left\{1+
\frac{4(p-1)}{\,p^2\gamma^{2p}\,}
\left[1+
O\left(\frac{1}{|\log\varepsilon|}\right)
\right]\right\}
<2\pi(k+2l)
\qquad
\textrm{and}
\qquad
\beta_\lambda\uparrow 2\pi(k+2l),
\end{eqnarray}
%\begin{equation}\label{1.7a}
%\aligned
%\beta_\lambda
%=
%\frac{\lambda p}{2}
%\left(
%\int_{\Omega}
%\big(e^{|u_\lambda|^p}-1\big)dx
%\right)^{\frac{2-p}{p}}
%\left(
%\int_{\Omega}
%|u_\lambda|^p
%e^{|u_\lambda|^p}dx
%\right)^{\frac{2(p-1)}{p}}
%\,\big\uparrow
%\, 4\pi m,
%\endaligned
%\end{equation}
but if $1<p<2$,
\begin{eqnarray}\label{1.30}
\beta_\lambda\geq
2\pi(k+2l)
\left\{1+
\frac{4(p-1)}{\,p^2\gamma^{2p}\,}
\left[1+
O\left(\frac{1}{|\log\varepsilon|}\right)
\right]\right\}
>2\pi(k+2l)
\qquad
\textrm{and}
\qquad
\beta_\lambda\downarrow 2\pi(k+2l).
\end{eqnarray}
In addition, given  $0<p<2$ and  $\beta>0$, we consider
the functional
\begin{equation}\label{1.31}
\aligned
I_{p,\beta}(u)=\frac{2-p}2\left(
\frac{p\|u\|_{H^1(\Omega)}^2}{2\beta}
\right)^{\frac{p}{2-p}}
-\log
\int_{\Omega}\left(e^{u_{+}^p}-1
\right)dx,\,\,\ \ \,\,u\in H^1(\Omega),
\endaligned
\end{equation}
whose critical points are exactly these
 solutions of the mean-field type equation
\begin{equation}\label{1.32}
\left\{\aligned
&-\Delta u+u=\lambda  u^{p-1}e^{u^p},\,\,\,
\quad\quad\quad
u>0\,\,\ \,
\quad\quad
\textrm{in}\,\,\,\,\,\,
\Omega,\\[1mm]
&\frac{\partial u}{\partial\nu}=0\,\,
\qquad\quad\qquad\qquad\quad\quad\quad\quad\quad
\ \ \ \qquad\,
\textrm{on}\,\,\,\,
\partial\Omega,
\\[1mm]
&
\frac{\lambda p}{2}
\left(
\int_{\Omega}
\big(e^{u^p}-1\big)dx
\right)^{\frac{2-p}{p}}
\left(
\int_{\Omega}
u^p
e^{u^p}dx
\right)^{\frac{2(p-1)}{p}}=\beta,
\endaligned\right.
\end{equation}
where the positivity of  solutions arises from the maximum principle and $\lambda>0$
is given by the last relation of the energy $\beta$.
Very recently, when the smooth bounded planar domain $\Omega$ is replaced by a closed Riemann surface $(\Sigma,g_0)$, it has been proven in \cite{DMMT}
that for every $1<p<2$ and $\beta>0$, or for $p=1$ and $0<\beta\not\in4\pi\mathbb{N}$,
the functional (\ref{1.31})  always has a family of positive critical points $u_{p,\beta}$   in $H^1(\Sigma)$ whose pointwise blow-up
means $\lambda\rightarrow0+$ and can occur only for
$\beta\downarrow4\pi\mathbb{N^{\star}}$ or for $p\rightarrow1$ and $\beta\rightarrow 4\pi\mathbb{N^{\star}}$.
Notice that equation (\ref{1.32}) seems to have the same form as  problem (\ref{1.13}) except that
$\lambda>0$ is not a free parameter but
a number $\lambda=\lambda_\beta$
defined by  the relation of the energy parameter $\beta$. These, together with our
expansions (\ref{1.27})-(\ref{1.30}), seem to imply that
%for any smooth bounded planar domain $\Omega$ and
for
 any non-negative integers $k$,  $l$ with $k+l\geq1$,
 if the free parameter
$\beta\rightarrow 2\pi(k+2l)$
for $p=1$,
$\beta \uparrow 2\pi(k+2l)$
for $p\in(0,1)$, but $\beta \downarrow 2\pi(k+2l)$
for $p\in(1,2)$,
 then,
the functional (\ref{1.31})    will always have a  family of
positive bubbling critical points $u_\beta$ in $H^1(\Omega)$
that concentrates at  $l$ different interior points
and $k$ different boundary points  and $\lambda=\lambda_\beta\rightarrow0+$.

This article is organized as follows:
In Section $2$ we provide a good
approximation for the solution of problem (\ref{1.1}) and
estimate the scaling error created by this approximation.
Then we rewrite  problem (\ref{1.1}) in terms of a linear operator
$\mathcal{L}$ for which
a solvability theory
is performed  through solving a linearized  problem in Section $3$.
In Section $4$ we
solve  an auxiliary nonlinear problem.
In Section $5$ we reduce the problem of finding bubbling solutions
of (\ref{1.1}) to that of finding a critical point of a finite-dimensional function.
In section $6$ we give an asymptotic expansion of the energy functional associated to the approximate
solution. In the final section   we provide the detailed proof
of Theorems $1.1$-$1.2$.

Notation:
In this paper the letters $C$ and $D$ will always denote
a  universal positive  constant independent of $\lambda$ and $\varepsilon$,
which could be changed from one line to another.
The symbol $o(t)$ (respectively $O(t)$) will denote a quantity for which
$\frac{o(t)}{|t|}$ tends to zero
(respectively, $\frac{O(t)}{|t|}$ stays bounded )
as parameter $t$ goes to zero.
Moreover, we will use the notation
$o(1)$ (respectively $O(1)$)
to stand for a quantity which tends to zero
(respectively, which remains uniformly bounded) as $\lambda$ or $\varepsilon$
tends to zero.

\vspace{1mm}
\vspace{1mm}

\section{An approximation for the solution}
The original cells for the construction of an approximate solution of problem (\ref{1.1})
are based on the four-parameter family of functions
\begin{equation}\label{2.1}
\aligned
\omega_{\varepsilon,\mu,\xi}(z)=\log\frac{8\mu^2}{(\varepsilon^2\mu^2+|z-\xi|^2)^2},\quad
\,\,\varepsilon>0,\,\,\,\,\mu>0,\,\,\,\,\xi\in\mathbb{R}^2,
\endaligned
\end{equation}
which exactly solve
\begin{equation}\label{2.2}
\aligned
-\Delta
\omega=\varepsilon^2e^\omega\,\ \,\,\,\,
\textrm{in}\,\,\,\,
\mathbb{R}^2,
\qquad\quad\qquad
\int_{\mathbb{R}^2}\varepsilon^2e^\omega=8\pi.
\endaligned
\end{equation}
Set
\begin{equation}\label{2.3}
\aligned
\omega_{\mu}(z)=\omega_{1,\mu,(0,0)}(|z|)\equiv\log\frac{8\mu^2}{(\mu^2+|z|^2)^2}.
\endaligned
\end{equation}
The configuration space for $m$ concentration points $\xi=(\xi_1,\ldots,\xi_m)$ we
try to look for  is the following
\begin{eqnarray}\label{2.4}
\mathcal{O}_\varepsilon:=\left\{\,\xi=(\xi_1,\ldots,\xi_m)\in\Omega^l\times(\po)^{m-l}
\left|\,
\min_{i,j=1,\ldots,m,\,i\neq j}|\xi_i-\xi_j|>\frac{1}{|\log\varepsilon|^\kappa},
\,\,\,\,\,\,\,
\min_{1\leq i\leq l}\dist(\xi_i,\po)>\frac{1}{|\log\varepsilon|^\kappa}
\right.\right\},
\end{eqnarray}
where $l\in\{0,\ldots,m\}$,
$\varepsilon$ is sufficiently small and uniquely defined by $\lambda$ in  (\ref{1.5})-(\ref{2.7}),
and $\kappa$ is given by
\begin{eqnarray}\label{2.5}
\kappa=2(m^2+1).
\end{eqnarray}
Let  $m\in\mathbb{N}^*$ and
$\xi=(\xi_1,\ldots,\xi_m)\in\mathcal{O}_\varepsilon$ be fixed.
For numbers $\mu_i$,
$i=1,\ldots,m$, yet to be determined, but we always assume
\begin{equation}\label{2.6}
\aligned
1/C\leq\mu_i\leq |\log\varepsilon|^{C},
\qquad i=1,\ldots,m,
\endaligned
\end{equation}
for some $C>0$, independent of $\varepsilon$.
For each $i=1,\ldots,m$, we define
\begin{equation}\label{2.8}
\aligned
U_i(x)=\frac{1}{p\gamma^{p-1}}\left[
\omega_{\varepsilon,\mu_i,\xi_i}(x)
+\sum_{j=1}^4\left(\frac{p-1}{p}\right)^j\frac{1}{\gamma^{jp}}\omega^j_{\mu_{i}}\left(\frac{x-\xi_i}{\varepsilon}\right)
\right].
\endaligned
\end{equation}
%\begin{equation}\label{2.8}
%\aligned
%U_i(x)=\frac{1}{p\gamma^{p-1}}\left[
%\omega_{\varepsilon,\mu_i,\xi_i}(x)
%+\frac{p-1}p\frac{1}{\gamma^{p}}\omega^1_{\mu_{i}}\left(\frac{x-\xi_i}{\varepsilon}\right)
%+\left(\frac{p-1}{p}\right)^2\frac{1}{\gamma^{2p}}\omega^2_{\mu_{i}}\left(\frac{x-\xi_i}{\varepsilon}\right)
%+\left(\frac{p-1}{p}\right)^3\frac{1}{\gamma^{3p}}\omega^3_{\mu_{i}}\left(\frac{x-\xi_i}{\varepsilon}\right)
%\right].
%\endaligned
%\end{equation}
Here,    $\omega^j_{\mu_i}$, $j=1,2,3,4$,
 are radial solutions of
\begin{equation}\label{2.9}
\aligned
\Delta\omega^j_{\mu_i}+e^{\omega_{\mu_i}(|z|)}\omega^j_{\mu_i}=e^{\omega_{\mu_i}(|z|)}f^j_{\mu_i}
\qquad
\textrm{in}\,\,\ \,\,\mathbb{R}^2,
\endaligned
\end{equation}
with
\begin{equation}\label{2.10}
\aligned
f^1_{\mu_i}=-\left[
\omega_{\mu_i}+\frac{1}{2}(\omega_{\mu_i})^2
\right],
\endaligned
\end{equation}
and
\begin{equation}\label{2.11}
\aligned
f^2_{\mu_i}=-\left\{
\left[
\omega^{1}_{\mu_i}+\frac{p-2}{2(p-1)}(\omega_{\mu_i})^2\right]
+\omega_{\mu_i}\left[\omega^{1}_{\mu_i}+\frac{1}{2}(\omega_{\mu_i})^2\right]
+\omega_{\mu_i}\omega^{1}_{\mu_i}+\frac{p-2}{6(p-1)}(\omega_{\mu_i})^3
+\frac12\left[\omega^{1}_{\mu_i}+\frac{1}{2}(\omega_{\mu_i})^2\right]^2
\right\},
\endaligned
\end{equation}
and
\begin{eqnarray}\label{2.12}
f^3_{\mu_i}=-\left\{
\left[
\omega^{2}_{\mu_i}+\frac{p-2}{p-1}\omega_{\mu_i}\omega^1_{\mu_i}
+\frac{(p-2)(p-3)}{6(p-1)^2}(\omega_{\mu_i})^3
\right]
+
\left[
\omega^{1}_{\mu_i}+\frac{p-2}{2(p-1)}(\omega_{\mu_i})^2\right]
\left[\omega^{1}_{\mu_i}+\frac{1}{2}(\omega_{\mu_i})^2\right]
\right.
\qquad\qquad\quad\,
\nonumber\\
\left.
+\omega_{\mu_i}\left[
\omega^{2}_{\mu_i}+\omega_{\mu_i}\omega^{1}_{\mu_i}+\frac{p-2}{6(p-1)}(\omega_{\mu_i})^3
+\frac12\left(\omega^{1}_{\mu_i}+\frac{1}{2}(\omega_{\mu_i})^2\right)^2\right]
+\frac12(\omega^{1}_{\mu_i})^2
+
\omega_{\mu_i}\omega^{2}_{\mu_i}
+\frac{p-2}{2(p-1)}(\omega_{\mu_i})^2\omega^{1}_{\mu_i}
\right.
\nonumber\\
\left.
+\frac{(p-2)(p-3)}{24(p-1)^2}(\omega_{\mu_i})^4
+
\left[\omega^{1}_{\mu_i}+\frac{1}{2}(\omega_{\mu_i})^2\right]\left[
\omega^{2}_{\mu_i}+\omega_{\mu_i}\omega^{1}_{\mu_i}+\frac{p-2}{6(p-1)}(\omega_{\mu_i})^3\right]
+\frac16\left[\omega^{1}_{\mu_i}+\frac{1}{2}(\omega_{\mu_i})^2\right]^3
\right\},
\quad\,\,\,\,
\end{eqnarray}
and
\begin{eqnarray}\label{2.121}
f^4_{\mu_i}=-\left\{
\left[
\omega^{3}_{\mu_i}+\frac{p-2}{2(p-1)}(\omega^1_{\mu_i})^2
+\frac{p-2}{p-1}\omega_{\mu_i}\omega^2_{\mu_i}
+\frac{(p-2)(p-3)}{2(p-1)^2}(\omega_{\mu_i})^2\omega^1_{\mu_i}
+\frac{(p-2)(p-3)(p-4)}{24(p-1)^3}(\omega_{\mu_i})^4
\right]
\right.
\qquad\qquad\ \,
\nonumber\\
+\left[
\omega^{2}_{\mu_i}+\frac{p-2}{p-1}\omega_{\mu_i}\omega^1_{\mu_i}
+\frac{(p-2)(p-3)}{6(p-1)^2}(\omega_{\mu_i})^3
\right]
\left[\omega^{1}_{\mu_i}+\frac{1}{2}(\omega_{\mu_i})^2\right]
+\frac1{24}\left[\omega^{1}_{\mu_i}+\frac{1}{2}(\omega_{\mu_i})^2\right]^4
\qquad\qquad\qquad\qquad\quad\,\ \,\,\,
\nonumber\\
+\left[
\omega^{1}_{\mu_i}+\frac{p-2}{2(p-1)}(\omega_{\mu_i})^2\right]
\left[
\omega^{2}_{\mu_i}+\omega_{\mu_i}\omega^{1}_{\mu_i}+\frac{p-2}{6(p-1)}(\omega_{\mu_i})^3
+\frac12\left(\omega^{1}_{\mu_i}+\frac{1}{2}(\omega_{\mu_i})^2\right)^2
\right]
\qquad\qquad\qquad\qquad\qquad\qquad
\nonumber\\
+\,
\omega_{\mu_i}
\left[
\omega^{3}_{\mu_i}
+\frac12(\omega^{1}_{\mu_i})^2
+
\omega_{\mu_i}\omega^{2}_{\mu_i}
+\frac{p-2}{2(p-1)}(\omega_{\mu_i})^2\omega^{1}_{\mu_i}
+\frac{(p-2)(p-3)}{24(p-1)^2}(\omega_{\mu_i})^4
\right.
\qquad\qquad\qquad\qquad\qquad\qquad\qquad\quad\,
\nonumber\\
\left.
+
\left(\omega^{1}_{\mu_i}+\frac{1}{2}(\omega_{\mu_i})^2\right)
\left(
\omega^{2}_{\mu_i}+\omega_{\mu_i}\omega^{1}_{\mu_i}+\frac{p-2}{6(p-1)}(\omega_{\mu_i})^3
\right)
+\frac16\left(\omega^{1}_{\mu_i}+\frac{1}{2}(\omega_{\mu_i})^2\right)^3
\right]
\qquad\qquad\qquad\qquad\qquad\qquad\quad\,\,\,\,
\nonumber\\
+\left[
\omega^1_{\mu_i}\omega^{2}_{\mu_i}
+\frac{p-2}{2(p-1)}(\omega_{\mu_i})^2\omega^{2}_{\mu_i}
+\frac{p-2}{2(p-1)}\omega_{\mu_i}(\omega^{1}_{\mu_i})^2
+\frac{(p-2)(p-3)}{6(p-1)^2}(\omega_{\mu_i})^3\omega^{1}_{\mu_i}
+\frac{(p-2)(p-3)(p-4)}{120(p-1)^3}(\omega_{\mu_i})^5
\right]\nonumber\\
+\left[\omega^{1}_{\mu_i}+\frac{1}{2}(\omega_{\mu_i})^2\right]
\left[
\omega^{3}_{\mu_i}
+\frac12(\omega^{1}_{\mu_i})^2
+
\omega_{\mu_i}\omega^{2}_{\mu_i}
+\frac{p-2}{2(p-1)}(\omega_{\mu_i})^2\omega^{1}_{\mu_i}
+\frac{(p-2)(p-3)}{24(p-1)^2}(\omega_{\mu_i})^4
\right]
\qquad\qquad\qquad\qquad\quad
\nonumber\\
\left.
+\frac12\left[
\omega^{2}_{\mu_i}+\omega_{\mu_i}\omega^{1}_{\mu_i}+\frac{p-2}{6(p-1)}(\omega_{\mu_i})^3\right]^2
+\frac12\left[\omega^{1}_{\mu_i}+\frac{1}{2}(\omega_{\mu_i})^2\right]^2
\left[
\omega^{2}_{\mu_i}+\omega_{\mu_i}\omega^{1}_{\mu_i}+\frac{p-2}{6(p-1)}(\omega_{\mu_i})^3\right]
\right\},
\,\,\,\,\,
\textrm{for}
\,\,p\neq1.
\quad
\end{eqnarray}
According to   \cite{CI,EMP,MM}, it  readily follows that for any $j=1, 2, 3, 4$,
\begin{equation}\label{2.13}
\aligned
\omega^j_{\mu_i}(r)=\frac{D^j_{\mu_i}}{2}\log\left(1+\frac{r^2}{\mu_i^2}\right)+O\left(\frac{\mu_i}{\mu_i+r}\right),
\,\,\quad\,\,
\partial_{r}\omega^j_{\mu_i}(r)=\frac{D^j_{\mu_i} r}{\mu_i^2+r^2}+O\left(\frac{\mu_i}{\mu_i^2+r^2}\right)
\quad\,\,\,\textrm{as}\,\,\,r\rightarrow+\infty,\,\,\,r=|z|,
\endaligned
\end{equation}
where
\begin{equation}\label{2.14}
\aligned
D^j_{\mu_i}=\frac{1}{2\pi}\int_{\mathbb{R}^2}
\Delta\big[\omega^j_{\mu_i}(\mu_iy)\big]dy
\qquad\quad
\textrm{and}
\qquad\quad
D^j_{\mu_i}=8\int_{0}^{+\infty}t\frac{t^2-1}{(t^2+1)^3}f^j_{\mu_i}(\mu_i t)dt.
\endaligned
\end{equation}
Obviously, for every $j=1,2,3,4$, the coefficient $D^j_{\mu_i}$ has at most polynomial growth with respect to $\log\mu_i$.
Moreover,
\begin{equation}\label{2.15}
\aligned
D^1_{\mu_i}=4\log8-8-8\log\mu_i.
\endaligned
\end{equation}
Now we construct the approximate solution of problem (\ref{1.1}) by
\begin{equation}\label{2.16}
\aligned
U_\xi(x):=\sum_{i=1}^mPU_i(x)=\sum_{i=1}^m\big[U_i(x)+H_i(x)\big],
\endaligned
\end{equation}
where $H_i$  is a correction term defined as the solution of
\begin{equation}\label{2.17}
\aligned
\left\{\aligned
&-\Delta_a
H_i+H_i=\nabla\log a(x)\nabla U_i-U_i\,\,\,\,\,\,
\textrm {in}\,\,\,\,\,\,\Omega,\\[1mm]
&\frac{\partial H_i}{\partial \nu}=-\frac{\partial U_i}{\partial\nu}\,
\,\qquad\qquad\qquad\qquad\qquad\,\,\,
\textrm{on}\,\,\,\,\po.
\endaligned\right.
\endaligned
\end{equation}
To state the asymptotic behavior of each correction term $H_i$ in terms of $\varepsilon$,
$\mu_i$ and $\xi_i$, we first  use the convention
\begin{equation}\label{2.18}
\aligned
c_i=\left\{
\aligned
&8\pi,
\qquad\textrm{if}\quad\,\xi_i\in\Omega,
\\[1mm]
&4\pi,
\qquad\textrm{if}\quad \xi_i\in\po.
\endaligned
\right.
\endaligned
\end{equation}

\noindent{\bf Lemma 2.1.}\,\,{\it
For any
$\xi=(\xi_1,\ldots,\xi_m)\in\mathcal{O}_\varepsilon$
and  $\alpha\in(0,1)$,   then we have
\begin{eqnarray}\label{2.19}
H_i(x)=\frac{1}{p\gamma^{p-1}}\left\{
\left[
1-
\frac{1}{4}
\sum_{j=1}^4\left(\frac{p-1}{p}\right)^j\frac{D^j_{\mu_i}}{\gamma^{jp}}
\right]
c_i H_a(x,\xi_i)-\log(8\mu_i^2)
+\left[\sum_{j=1}^4\left(\frac{p-1}{p}\right)^j\frac{D^j_{\mu_i}}{\gamma^{jp}}
\right]
\log(\varepsilon\mu_i)
+O\left((\varepsilon\mu_i)^{\alpha/2}\right)
\right\}
\end{eqnarray}
uniformly in $\oo$,
where
$H_a$ is the regular part of the anisotropic Green's function defined in
{\upshape(\ref{1.8})}.
}

\begin{proof}
Inserting  (\ref{2.1}), (\ref{2.8}) and (\ref{2.13})
into (\ref{2.17}), we have that for any
 $\beta\in(0,1)$,
 $$
\aligned
\left\{\aligned&
-\Delta_a
H_i+H_i=\frac{1}{p\gamma^{p-1}}\left\{
\left[-4
+\sum_{j=1}^4\left(\frac{p-1}{p}\right)^j\frac{D^j_{\mu_i}}{\gamma^{jp}}
\right]\left[
\frac{(x-\xi_i)\cdot\nabla\log a(x)}{\varepsilon^2\mu_i^2+|x-\xi_i|^2}
-\frac{1}{2}\log\big(\varepsilon^2\mu_i^2+|x-\xi_i|^2\big)\right]
-\log(8\mu_i^2)
\right.\\[1mm]
&\left.\,\quad\qquad
+\left[
\sum_{j=1}^4\left(\frac{p-1}{p}\right)^j\frac{D^j_{\mu_i}}{\gamma^{jp}}
\right]\log(\varepsilon\mu_i)
+\frac{p-1}{p}\frac{1}{\gamma^p}
O_{\large L^{\infty}\big(\Omega\setminus B_{(\varepsilon\mu_i)^{\beta/2}}(\xi_i)\big)}
\left(\frac{\varepsilon\mu_i}{\varepsilon\mu_i+|x-\xi_i|}
+\frac{\varepsilon\mu_i}{\varepsilon^2\mu_i^2+|x-\xi_i|^2}\right)
\right.\\[1mm]
&\left.\,\quad\qquad
+\left[\sum_{j=1}^4\left(\frac{p-1}{p}\right)^j\frac{D^j_{\mu_i}}{\gamma^{jp}}
\right]
O_{\large L^{\infty}\big(\Omega\bigcap B_{(\varepsilon\mu_i)^{\beta/2}}(\xi_i)\big)}\left(\frac{|(x-\xi_i)\cdot\nabla\log a(x)|}{\varepsilon^2\mu_i^2+|x-\xi_i|^2}
+\log\frac{\varepsilon^2\mu_i^2+|x-\xi_i|^2}{\varepsilon^2\mu_i^2}
\right)
\right\}
\quad\,\,
\textrm{in}\,\,\,\ \ \,\Omega,\\[1mm]
&\frac{\partial H_i}{\partial \nu}=-\frac{1}{p\gamma^{p-1}}\left\{
\left[-4
+\sum_{j=1}^4\left(\frac{p-1}{p}\right)^j\frac{D^j_{\mu_i}}{\gamma^{jp}}
\right]
\frac{(x-\xi_i)\cdot\nu(x)}{\varepsilon^2\mu_i^2+|x-\xi_i|^2}
+\frac{p-1}{p}\frac{1}{\gamma^p}
O_{\large L^{\infty}\big(\partial\Omega\setminus B_{(\varepsilon\mu_i)^{\beta/2}}(\xi_i)\big)}\left(\frac{\varepsilon\mu_i}{\varepsilon^2\mu_i^2+|x-\xi_i|^2}\right)
\right.\\[1mm]
&\left.\,
\quad\qquad
+\left[\sum_{j=1}^4\left(\frac{p-1}{p}\right)^j\frac{D^j_{\mu_i}}{\gamma^{jp}}
\right]
O_{\large L^{\infty}\big(\partial\Omega\bigcap B_{(\varepsilon\mu_i)^{\beta/2}}(\xi_i)\big)}\left(\frac{|(x-\xi_i)\cdot\nu(x)|}{\varepsilon^2\mu_i^2+|x-\xi_i|^2}\right)
\right\}
\,\qquad\qquad\qquad\qquad\qquad\qquad\,
\,\textrm{on}\,\,\ \ \,\po.
\endaligned\right.
\endaligned
$$
Using (\ref{1.7})-(\ref{1.8}) we get that
the regular part of Green's function, $H_a(x,\xi_i)$,  satisfies
\begin{equation*}\label{}
\left\{\aligned
&-\Delta_a
H_a(x,\xi_i)+H_a(x,\xi_i)=
\frac{4}{c_i}\log|x-\xi_i|-
\frac{4}{c_i}\frac{(x-\xi_i)\cdot\nabla\log a(x)}{|x-\xi_i|^2}
\quad
\,\,\textrm{in}\,\,\,\,\,\Omega,\\[1mm]
&\frac{\partial H_a(x,\xi_i)}{\partial \nu}
=\frac{4}{c_i}\frac{(x-\xi_i)\cdot\nu(x)}{|x-\xi_i|^2}
\qquad\qquad\qquad\qquad
\qquad\qquad\qquad\qquad\quad
\textrm{on}\,\,\,\po.
\endaligned\right.
\end{equation*}
Set
$$
\aligned
Z_i(x)=p\gamma^{p-1}H_i(x)
-\left[
1-
\frac{1}{4}
\sum_{j=1}^4\left(\frac{p-1}{p}\right)^j\frac{D^j_{\mu_i}}{\gamma^{jp}}
\right]c_iH_a(x,\xi_i)
+\log(8\mu_i^2)
-\left[\sum_{j=1}^4\left(\frac{p-1}{p}\right)^j\frac{D^j_{\mu_i}}{\gamma^{jp}}
\right]\log(\varepsilon\mu_i).
\endaligned
$$
Then
$$
\aligned
\left\{\aligned&
-\Delta_a
Z_i+Z_i=
\left[-4
+
\sum_{j=1}^4\left(\frac{p-1}{p}\right)^j\frac{D^j_{\mu_i}}{\gamma^{jp}}
\right]\left[\,\frac{1}{2}\log\left(
\frac{|x-\xi_i|^2}{\varepsilon^2\mu_i^2+|x-\xi_i|^2}\right)
-\frac{(x-\xi_i)\cdot\nabla\log a(x)}{|x-\xi_i|^2}\cdot\frac{\varepsilon^2\mu_i^2}{\varepsilon^2\mu_i^2+|x-\xi_i|^2}
\right]\\[1mm]
&\quad\qquad
+\left[\sum_{j=1}^4\left(\frac{p-1}{p}\right)^j\frac{D^j_{\mu_i}}{\gamma^{jp}}
\right]
O_{\large L^{\infty}\big(\Omega\bigcap B_{(\varepsilon\mu_i)^{\beta/2}}(\xi_i)\big)}\left(\frac{|(x-\xi_i)\cdot\nabla\log a(x)|}{\varepsilon^2\mu_i^2+|x-\xi_i|^2}
+\log\frac{\varepsilon^2\mu_i^2+|x-\xi_i|^2}{\varepsilon^2\mu_i^2}
\right)\\[1mm]
&\quad\qquad
+\frac{p-1}{p}\frac{1}{\gamma^p}
O_{\large L^{\infty}\big(\Omega\setminus B_{(\varepsilon\mu_i)^{\beta/2}}(\xi_i)\big)}
\left(\frac{\varepsilon\mu_i}{\varepsilon\mu_i+|x-\xi_i|}
+\frac{\varepsilon\mu_i}{\varepsilon^2\mu_i^2+|x-\xi_i|^2}\right)
\quad\qquad\qquad\qquad\qquad\qquad\,\qquad
\textrm{in}\,\,\ \,\,\ \,\,\Omega,\\[1mm]
&\frac{\partial Z_i}{\partial \nu}=
\left[-4
+
\sum_{j=1}^4\left(\frac{p-1}{p}\right)^j\frac{D^j_{\mu_i}}{\gamma^{jp}}
\right]
\frac{(x-\xi_i)\cdot\nu(x)}{|x-\xi_i|^2}
\frac{\varepsilon^2\mu_i^2}{\varepsilon^2\mu_i^2+|x-\xi_i|^2}
+\frac{p-1}{p}\frac{1}{\gamma^p}
O_{\large L^{\infty}\big(\partial\Omega\setminus B_{(\varepsilon\mu_i)^{\beta/2}}(\xi_i)\big)}\left(\frac{\varepsilon\mu_i}{\varepsilon^2\mu_i^2+|x-\xi_i|^2}\right)
\\[1mm]
&\qquad\quad
+\left[\sum_{j=1}^4\left(\frac{p-1}{p}\right)^j\frac{D^j_{\mu_i}}{\gamma^{jp}}
\right]
O_{\large L^{\infty}\big(\partial\Omega\bigcap B_{(\varepsilon\mu_i)^{\beta/2}}(\xi_i)\big)}\left(\frac{|(x-\xi_i)\cdot\nu(x)|}{\varepsilon^2\mu_i^2+|x-\xi_i|^2}\right)
\qquad\qquad\qquad\qquad\qquad\qquad\quad\,\,\,\textrm{on}\,\ \ \,\,\,\po.
\endaligned\right.
\endaligned
$$
Direct computations show that there exists a constant $C>0$ such that
for any  $q>1$,
$$
\aligned
\left\|
\log
\left(\frac{|x-\xi_i|^2}{\varepsilon^2\mu_i^2+|x-\xi_i|^2}\right)
\right\|_{L^q(\Omega)}\leq C(\varepsilon\mu_i)^{\frac{2}{q}},
\endaligned
$$
and
$$
\aligned
\left\|
\frac{\varepsilon\mu_i}{\varepsilon\mu_i+|x-\xi_i|}
+\frac{\varepsilon\mu_i}{\varepsilon^2\mu_i^2+|x-\xi_i|^2}
\right\|_{L^q\big(\Omega\setminus B_{(\varepsilon\mu_i)^{\beta/2}}(\xi_i)\big)}
\leq C(\varepsilon\mu_i)^{1-\beta+\frac{\beta}{q}},
\endaligned
$$
and for  any $1<q<2$,
$$
\aligned
\left\|
\frac{|(x-\xi_i)\cdot\nabla\log a(x)|}{\varepsilon^2\mu_i^2+|x-\xi_i|^2}
+\log\frac{\varepsilon^2\mu_i^2+|x-\xi_i|^2}{\varepsilon^2\mu_i^2}
\right\|_{L^q\big(\Omega\bigcap B_{(\varepsilon\mu_i)^{\beta/2}}(\xi_i)\big)}
\leq C
(\varepsilon\mu_i)^{\beta(\frac{1}{q}-\frac12)},
\endaligned
$$
and
$$
\aligned
\left\|
\frac{(x-\xi_i)\cdot\nabla\log a(x)}{|x-\xi_i|^2}\cdot\frac{\varepsilon^2\mu_i^2}{\varepsilon^2\mu_i^2+|x-\xi_i|^2}
\right\|_{L^q(\Omega)}\leq C(\varepsilon\mu_i)^{\frac{2}{q}-1}.
\endaligned
$$
Hence for any $\xi_i\in\oo$ and any  $1<q<2$,
\begin{equation*}\label{}
\aligned
\big\|-\Delta_a
Z_i+Z_i\big\|_{L^q(\Omega)}\leq C(\varepsilon\mu_i)^{\beta(\frac{1}{q}-\frac{1}{2})}.
\endaligned
\end{equation*}
On the other hand,
if $\xi_i\in\po$, from  the fact that
$|(x-\xi_i)\cdot\nu(x)|\leq C|x-\xi_i|^2$ for any $x\in\po$ (see \cite{AP})
we can compute  that for any $q>1$,
$$
\aligned
\left\|
\frac{(x-\xi_i)\cdot\nu(x)}{|x-\xi_i|^2}
\cdot
\frac{\varepsilon^2\mu_i^2}{\varepsilon^2\mu_i^2+|x-\xi_i|^2}
\right\|_{L^q(\partial\Omega)}
\leq
 C(\varepsilon\mu_i)^{\frac{1}{q}},
\endaligned
$$
$$
\aligned
\left\|
\frac{(x-\xi_i)\cdot\nu(x)}{\,\varepsilon^2\mu_i^2+|x-\xi_i|^2\,}
\right
\|_{L^q\big(\partial\Omega\bigcap B_{(\varepsilon\mu_i)^{\beta/2}}(\xi_i)\big)}
\leq C(\varepsilon\mu_i)^{\frac{\beta}{2q}},
\endaligned
$$
$$
\aligned
\left\|
\frac{\varepsilon\mu_i}{\,\varepsilon^2\mu_i^2+|x-\xi_i|^2\,}
\right
\|_{L^q\big(\partial\Omega\setminus B_{(\varepsilon\mu_i)^{\beta/2}}(\xi_i)\big)}
\leq C\left[\varepsilon\mu_i
+
(\varepsilon\mu_i)^{1-\beta+\frac{\beta}{2q}}
\right],
\endaligned
$$
then
\begin{equation*}\label{}
\aligned
\left\|\frac{\partial Z_i}{\partial \nu}\right\|_{L^{q}(\partial\Omega)}
\leq C(\varepsilon\mu_i)^{\frac{\beta}{2q}}.
\endaligned
\end{equation*}
While if $\xi_i\in\Omega$,  by the definition of $\mathcal{O}_\varepsilon$ in (\ref{2.4}) we easily find
$$
\aligned
\left\|
\frac{\varepsilon\mu_i}{\varepsilon^2\mu_i^2+|x-\xi_i|^2}
\right\|_{L^{\infty}(\partial\Omega)}\leq
 \varepsilon\mu_i|\log\varepsilon|^{2\kappa},
\endaligned
$$
$$
\aligned
\left\|
\frac{(x-\xi_i)\cdot\nu(x)}{|x-\xi_i|^2}
\cdot
\frac{\varepsilon^2\mu_i^2}{\varepsilon^2\mu_i^2+|x-\xi_i|^2}
\right\|_{L^{\infty}(\partial\Omega)}
\leq
\varepsilon^2\mu_i^2|\log\varepsilon|^{3\kappa},
\endaligned
$$
then
\begin{equation*}\label{}
\aligned
\left\|\frac{\partial Z_i}{\partial \nu}\right\|_{L^{\infty}(\partial\Omega)}\leq C \varepsilon\mu_i|\log\varepsilon|^{2\kappa-1}.
\endaligned
\end{equation*}
As a consequence, from
elliptic regularity theory we have that for any $1<q<2$ and any $0<\theta<1/q$,
$$
\aligned
\left\|Z_i\right\|_{W^{1+\theta,q}(\Omega)}\leq
C\left(\big\|-\Delta_a
Z_i+Z_i\big\|_{L^q(\Omega)}+\left\|\frac{\partial Z_i}{\partial \nu}\right\|_{L^{q}(\partial\Omega)}\right)
\leq C(\varepsilon\mu_i)^{\beta(\frac{1}{q}-\frac{1}{2})}.
\endaligned
$$
By  Morrey's embedding theorem,
$$
\aligned
\left\|Z_i\right\|_{C^\tau(\overline{\Omega})}\leq C(\varepsilon\mu_i)^{\beta(\frac{1}{q}-\frac{1}{2})},
\endaligned
$$
where $0<\tau<1/2+1/q$,
which implies that expansion (\ref{2.19}) holds with $\alpha=2\beta(1/q-1/2)$.
\end{proof}
\vspace{1mm}
\vspace{1mm}
\vspace{1mm}
\vspace{1mm}

From Lemma 2.1 we can easily prove that away from each point $\xi_i$,
namely  $|x-\xi_i|\geq 1/|\log\varepsilon|^{2\kappa}$ for any $i=1,\ldots,m$,
\begin{equation}\label{2.20}
\aligned
U_{\xi}(x)=\frac{1}{p\gamma^{p-1}}
\sum\limits_{i=1}^{m}
\left\{
\left[
1
-\frac14
\sum_{j=1}^4\left(\frac{p-1}{p}\right)^j\frac{D^j_{\mu_i}}{\gamma^{jp}}
\right]
c_i G_a(x,\xi_i)
+O\left(
(\varepsilon\mu_i)^{\alpha/2}
\right)
\right\}.
\endaligned
\end{equation}
While if
$|x-\xi_i|< 1/|\log\varepsilon|^{2\kappa}$  with some $i$,  from the
fact that $H_a(\cdot,\xi_k)\in C^{\alpha}(\oo)$ for any $\xi_k\in\oo$ and
any
$\alpha\in(0,1)$
we find
$$
\aligned
PU_i(x)=&\,\frac{1}{p\gamma^{p-1}}
\left\{
p\gamma^{p}+
\omega_{\mu_i}\left(\frac{x-\xi_i}{\varepsilon}\right)
+\sum_{j=1}^4\left(\frac{p-1}{p}\right)^j\frac{1}{\gamma^{jp}}\omega^j_{\mu_{i}}\left(\frac{x-\xi_i}{\varepsilon}\right)
+\left[
1
-\frac14
\sum_{j=1}^4\left(\frac{p-1}{p}\right)^j\frac{D^j_{\mu_i}}{\gamma^{jp}}
\right]
c_i H_a(\xi_i,\xi_i)
\right.\\[1mm]
&\left.
-\log(8\mu_i^2)
+\left[
\sum_{j=1}^4\left(\frac{p-1}{p}\right)^j\frac{D^j_{\mu_i}}{\gamma^{jp}}
\right]\log(\varepsilon\mu_i)
+O\left(
|x-\xi_i|^\alpha
+(\varepsilon\mu_i)^{\alpha/2}\right)
\right\},
\endaligned
$$
and for any $k\neq i$,
$$
\aligned
PU_k(x)=\frac{1}{p\gamma^{p-1}}
\left\{
\left[
1
-\frac14
\sum_{j=1}^4\left(\frac{p-1}{p}\right)^j\frac{D^j_{\mu_k}}{\gamma^{jp}}
\right]
c_k G_a(\xi_i,\xi_k)
+O\left(|x-\xi_i|^\alpha+
(\varepsilon\mu_k)^{\alpha/2}
\right)
\right\}.
\endaligned
$$
Thus if
$|x-\xi_i|< 1/|\log\varepsilon|^{2\kappa}$,
\begin{eqnarray}\label{2.21}
U_{\xi}(x)=\frac{1}{p\gamma^{p-1}}
\left[
p\gamma^{p}+
\omega_{\mu_i}\left(\frac{x-\xi_i}{\varepsilon}\right)
+\sum_{j=1}^4\left(\frac{p-1}{p}\right)^j\frac{1}{\gamma^{jp}}\omega^j_{\mu_{i}}\left(\frac{x-\xi_i}{\varepsilon}\right)
+\,O\left(
|x-\xi_i|^\alpha
+\sum_{k=1}^m(\varepsilon\mu_k)^{\alpha/2}
\right)
\right]
\end{eqnarray}
will be a good  approximation for the  solution of problem (\ref{1.1})
provided that for each $i=1,\ldots,m$, the
concentration
parameter $\mu_i$   satisfies the nonlinear system
\begin{eqnarray}\label{2.22}
\log\big(8\mu_i^2\big)=\left[
1-
\frac{1}{4}
\sum_{j=1}^4\left(\frac{p-1}{p}\right)^j\frac{D^j_{\mu_i}}{\gamma^{jp}}
\right]
c_i H_a(\xi_i,\xi_i)
+\left[\sum_{j=1}^4\left(\frac{p-1}{p}\right)^j\frac{D^j_{\mu_i}}{\gamma^{jp}}
\right]\log(\varepsilon\mu_i)
\nonumber\\
+\sum_{k=1,\,k\neq i}^m
\left[
1-
\frac14\sum_{j=1}^4\left(\frac{p-1}{p}\right)^j\frac{D^j_{\mu_k}}{\gamma^{jp}}
\right]
c_k G_a(\xi_i,\xi_k).
\qquad\qquad\qquad\qquad\quad\,
\end{eqnarray}
It is necessary to  point out that from  {\upshape(\ref{2.7})},
{\upshape(\ref{2.14})},
{\upshape(\ref{2.15})}
 and the Implicit Function Theorem
we readily have that for
any  sufficiently small $\varepsilon$
and any points $\xi=(\xi_1,\ldots,\xi_m)\in\mathcal{O}_\varepsilon$,
there is a unique solution
$\mu=(\mu_1,\ldots,\mu_m)$ for  system {\upshape(\ref{2.22})}
under  assumption  {\upshape(\ref{2.6})}. Moreover,
for any $i=1,\ldots,m$,
\begin{equation*}\label{2.23}
\aligned
\big|D_{\xi}\log\mu_i\big|\leq C|\log\varepsilon|^{\kappa},
\endaligned
\end{equation*}
and
\begin{equation}\label{2.24}
\aligned
\log\big(8\mu_i^2\big)=\left\{
\frac{2(p-1)}{2-p}(1-\log8)+\frac{1}{2-p}\left[c_iH_a(\xi_i,\xi_i)
+\sum\large_{k=1,\,k\neq i}^m
c_k G_a(\xi_i,\xi_k)\right]
\right\}
\left[\,1+O\left(\frac{\log^2|\log\varepsilon|}{|\log\varepsilon|}\right)\right].
\endaligned
\end{equation}

\vspace{1mm}
\vspace{1mm}
\vspace{1mm}

\noindent{\bf Remark 2.2.}\,\,
Let us mention that $p\gamma^{p-1} U_\xi$ is a
positive  function over $\overline{\Omega}$. In fact,
for any $|x-\xi_i|< 1/|\log\varepsilon|^{2\kappa}$,
$$
\aligned
p\gamma^{p}+
\omega_{\mu_i}\left(\frac{x-\xi_i}{\varepsilon}\right)
+\sum_{j=1}^4\left(\frac{p-1}{p}\right)^j\frac{1}{\gamma^{jp}}\omega^j_{\mu_{i}}\left(\frac{x-\xi_i}{\varepsilon}\right)
&\geq  p\gamma^{p}+\log\frac{8}{\mu_i^2}
+\left(
-2+
\frac{p-1}{p}\frac{D^1_{\mu_i}}{2\gamma^{p}}
\right)
\log\left(1+\frac{1}{\varepsilon^2\mu_i^2|\log\varepsilon|^{4\kappa}}\right)+o\left(1\right)\\
&=8\kappa\log|\log\varepsilon|
+
\log\big(8\mu_i^2\big)
+\frac{1}4(p-1)D^1_{\mu_i}
+o\left(1\right)\\[1mm]
&=
8\kappa\log|\log\varepsilon|
+
(2-p)\log\mu_i^2+p\log8+2-2p
+o\left(1\right)\\[1mm]
&\geq
8\kappa\log|\log\varepsilon|
+
(2-p)\log\frac{1}{C^2}.
\endaligned
$$
%$$
%\aligned
%&p\gamma^{p}+
%\omega_{\mu_i}\left(\frac{x-\xi_i}{\varepsilon}\right)
%+\frac{p-1}p\frac{1}{\gamma^{p}}\omega^1_{\mu_{i}}\left(\frac{x-\xi_i}{\varepsilon}\right)
%+\left(\frac{p-1}{p}\right)^2\frac{1}{\gamma^{2p}}\omega^2_{\mu_{i}}\left(\frac{x-\xi_i}{\varepsilon}\right)\\
%&\geq p\gamma^{p}+\log8-2\log\mu_i+\left[
%-2+\frac{p-1}{p}\frac{1}{\gamma^p}\frac{D^1_{\mu_i}}{2}
%+\left(\frac{p-1}{p}\right)^2\frac{1}{\gamma^{2p}}\frac{D^2_{\mu_i}}{2}
%\right]\log\left(1+\frac{1}{\varepsilon^2\mu_i^2|\log\varepsilon|^{4\kappa}}\right)+o\big(1\big)\\
%&\geq\log8+2\log\mu_i+
%8\kappa\log|\log\varepsilon|
%+\frac{1}{4}
%(p-1)D^1_{\mu_i}
%+o\big(1\big)\\[1mm]
%&\geq
%8\kappa\log|\log\varepsilon|
%-C.
%\endaligned
%$$
Clearly, by (\ref{2.21}) we obtain
$p\gamma^{p-1}U_\xi(x)>4\kappa\log|\log\varepsilon|>0$ for any $|x-\xi_i|< 1/|\log\varepsilon|^{2\kappa}$.
Also, by the maximum principle we have that for each $\xi_i$,
 $G_a(x,\xi_i)>0$
over $\overline{\Omega}$ and hence by (\ref{2.20}),
$p\gamma^{p-1}U_\xi$ is a  positive  function over $\overline{\Omega}$.

\vspace{1mm}
\vspace{1mm}
\vspace{1mm}

Let us  perform  the change of variables
\begin{equation*}\label{2.25}
\aligned
\upsilon(y)=p\gamma^{p-1}u(\varepsilon y)-p\gamma^{p},\,\,\,
\,\quad\forall\,\,y\in\Omega_{\varepsilon}:=\varepsilon^{-1}\Omega.
\endaligned
\end{equation*}
By the definitions of
$\varepsilon$ and $\gamma$ in (\ref{1.5})-(\ref{2.7}), respectively,
we can rewrite equation (\ref{1.1}) in the following form
\begin{equation}\label{2.26}
\begin{array}{ll}
\left\{\aligned
&-\Delta_{a(\varepsilon y)}\upsilon+\varepsilon^2\upsilon=f(\upsilon)-p\gamma^{p}\varepsilon^2,
\,\,\,
\,\,
\textrm{in}\,\,\,\,\,\Omega_\varepsilon,\\[1mm]
&\frac{\partial \upsilon}{\partial\nu}=0
\qquad\qquad\qquad\qquad\qquad\qquad
\,
\textrm{on}\,\,\,\partial\Omega_\varepsilon,
\endaligned\right.
\end{array}
\end{equation}
where
\begin{equation}\label{2.27}
\aligned
f(\upsilon)=\left(1+\frac{\upsilon}{p\gamma^p}\right)^{p-1}
e^{\gamma^p\left[\left(1+\frac{\upsilon}{p\gamma^p}\right)^p-1\right]}.
\endaligned
\end{equation}
We
write $\xi_i'=\xi_i/\varepsilon$, $i=1,\ldots,m$
and define the initial approximate solution of (\ref{2.26}) as
\begin{equation}\label{2.28}
\aligned
V_{\xi'}(y)=p\gamma^{p-1}U_\xi(\varepsilon y)-p\gamma^{p},
\endaligned
\end{equation}
with
$\xi'=(\xi_1',\ldots,\xi_m')$
and
$U_\xi$  defined in (\ref{2.16}).  What remains of this paper is to
 look for solutions of problem (\ref{2.26}) in the form
$\upsilon=V_{\xi'}+\phi$, where $\phi$ will represent a higher-order correction.
In terms of $\phi$,
problem (\ref{2.26}) becomes
\begin{equation}\label{2.29}
\aligned
\left\{\aligned
&\mathcal{L}(\phi)=-\big[
E_{\xi'}+N(\phi)
\big]
\quad\textrm{in}\,\,\,\,\,\,\Omega_\varepsilon,\\
&\frac{\partial \phi}{\partial\nu}=0
\quad\qquad\qquad\qquad\,\,\,\,
\textrm{on}\,\,\,\,
\partial\Omega_\varepsilon,
\endaligned\right.\endaligned
\end{equation}
where
\begin{equation}\label{2.30}
\aligned
\mathcal{L}(\phi)=-\Delta_{a(\varepsilon y)}\phi+\varepsilon^2\phi-W_{\xi'}\phi\,\quad\,\textrm{with}\,\quad\,
W_{\xi'}=f'(V_{\xi'}),
\endaligned
\end{equation}
and
\begin{equation}\label{2.31}
\aligned
E_{\xi'}=-\Delta_{a(\varepsilon y)}V_{\xi'}+\varepsilon^2V_{\xi'}-f(V_{\xi'})+p\gamma^{p}\varepsilon^2,
\quad\quad\quad
N(\phi)=-\big[f(V_{\xi'}+\phi)-f(V_{\xi'})-f'(V_{\xi'})\phi\big].
\endaligned
\end{equation}

For any  $\xi=(\xi_1,\ldots,\xi_m)\in\mathcal{O}_{\varepsilon}$ and $h\in L^\infty(\Omega_\varepsilon)$, let us  introduce a weighted $L^\infty$-norm
defined as
\begin{equation}\label{2.32}
\aligned
\|h\|_{*}=
\sup_{y\in\Omega_\varepsilon}\left|\left(
\sum\limits_{i=1}^m\frac{\mu_i^\sigma}{(\mu_i+|y-\xi'_i|)^{2+\sigma}}
+\varepsilon^2
\right)^{-1}h(y)
\right|,
\endaligned
\end{equation}
where $\sigma>0$ is small   but  fixed,  independent of $\varepsilon$.
With respect to the $\|\cdot\|_{*}$-norm,
the error term $E_{\xi'}$ defined in (\ref{2.31}) can be estimated as follows.

\vspace{1mm}
\vspace{1mm}
\vspace{1mm}

\noindent{\bf Proposition 2.3.}\,\,{\it
There exists a constant $C>0$ such that
for  any $\xi=(\xi_1,\ldots,\xi_m)\in\mathcal{O}_{\varepsilon}$
and for  any $\varepsilon$ small enough,
\begin{equation}\label{2.33}
\aligned
\|E_{\xi'}\|_{*}\leq\frac{C}{\gamma^{4p}}=O\left(\frac{1}{|\log\varepsilon|^4}\right).
\endaligned
\end{equation}
}\noindent\begin{proof}
From (\ref{2.8}),  (\ref{2.16}), (\ref{2.17}) and (\ref{2.28}) we obtain
$$
\aligned
-\Delta_{a(\varepsilon y)}V_{\xi'}+\varepsilon^2V_{\xi'}+p\gamma^{p}\varepsilon^2
&=p\gamma^{p-1}\varepsilon^{2}\sum_{i=1}^m\left[-\Delta_{a}\big(U_i+H_i\big)+\big(U_i+H_i\big)\right]
=-p\gamma^{p-1}\varepsilon^{2}\sum_{i=1}^m\Delta U_i\\
&=-\varepsilon^2\sum_{i=1}^m\left[\Delta
\omega_{\varepsilon,\mu_i,\xi_i}(x)
+\frac{1}{\varepsilon^2}\sum_{j=1}^4\left(\frac{p-1}{p}\right)^j\frac{1}{\gamma^{jp}}\Delta\omega^j_{\mu_{i}}\left(\frac{x-\xi_i}{\varepsilon}\right)
\right].
\endaligned
$$
%$$
%\aligned
%-\Delta_{a(\varepsilon y)}V_{\xi'}+\varepsilon^2V_{\xi'}+p\gamma^{p}\varepsilon^2
%=p\gamma^{p-1}\varepsilon^{2}\sum_{i=1}^m\left[-\Delta_{a}\big(U_i+H_i\big)+\big(U_i+H_i\big)\right]
%=-p\gamma^{p-1}\varepsilon^{2}\sum_{i=1}^m\Delta U_i.
%\endaligned
%$$
Furthermore,  by  (\ref{2.1})-(\ref{2.3})  and (\ref{2.9}) we compute
\begin{eqnarray}\label{2.34}
-\Delta_{a(\varepsilon y)}V_{\xi'}+\varepsilon^2V_{\xi'}+p\gamma^{p}\varepsilon^2
=\sum_{i=1}^me^{\omega_{\mu_i}\left(y-\xi'_i\right)}\left[
1
+
\sum_{j=1}^4\left(\frac{p-1}{p}\right)^j\frac{1}{\gamma^{jp}}
\big(
\omega^j_{\mu_i}-f^j_{\mu_i}
\big)
\right]\left(y-\xi'_i\right).
\end{eqnarray}
%\begin{eqnarray}\label{2.34}
%-\Delta_{a(\varepsilon y)}V_{\xi'}+\varepsilon^2V_{\xi'}+p\gamma^{p}\varepsilon^2
%=\sum_{i=1}^m\left[
%e^{\omega_{\mu_i}\left(y-\xi'_i\right)}
%-\sum_{j=1}^4\left(\frac{p-1}{p}\right)^j\frac{1}{\gamma^{jp}}\Delta\omega^j_{\mu_{i}}
%(y-\xi'_i)\right]
%\qquad\,\,\,\,\,\quad
%\nonumber
%\\
%=\sum_{i=1}^me^{\omega_{\mu_i}\left(y-\xi'_i\right)}\left[
%1
%+
%\sum_{j=1}^4\left(\frac{p-1}{p}\right)^j\frac{1}{\gamma^{jp}}
%\big(
%\omega^j_{\mu_i}-f^j_{\mu_i}
%\big)
%\right]\left(y-\xi'_i\right).
%\end{eqnarray}
From (\ref{2.3}),  (\ref{2.6}) and (\ref{2.13})
we have that if $|y-\xi'_i|\geq 1/(\varepsilon|\log\varepsilon|^{2\kappa})$ for any $i=1,\ldots,m$,
$$
\aligned
\omega_{\mu_i}(y-\xi'_i)=4\log\varepsilon+O\big(\log|\log\varepsilon|\big),\,\qquad\qquad\,
\omega_{\mu_i}^j(y-\xi'_i)=-D^{j}_{\mu_i}\log\varepsilon+O\big(\log|\log\varepsilon|\big),
\quad
j=1,2,3,4,
\endaligned
$$
and then, by (\ref{2.10})-(\ref{2.121}),
\begin{equation}\label{2.35}
\aligned
-\Delta_{a(\varepsilon y)}V_{\xi'}+\varepsilon^2V_{\xi'}+p\gamma^{p}\varepsilon^2
=\left[\sum_{i=1}^me^{\omega_{\mu_i}(y-\xi'_i)}
\right]
O\left(|\log\varepsilon|^4
\right).
\endaligned
\end{equation}
On the other hand, in the same region, by (\ref{2.20}) and (\ref{2.28}) we obtain
\begin{equation}\label{2.36}
\aligned
1+\frac{V_{\xi'}(y)}{p\gamma^p}
=
\frac{p\gamma^{p-1}U_{\xi}(\varepsilon y)}{p\gamma^p}
=O\left(\frac{\log|\log\varepsilon|}{|\log\varepsilon|}
\right),
\endaligned
\end{equation}
and hence,
\begin{equation}\label{2.37}
\aligned
f(V_{\xi'})=\left(1+\frac{V_{\xi'}}{p\gamma^p}\right)^{p-1}e^{\gamma^p\left[\left(1+\frac{V_{\xi'}}{p\gamma^p}\right)^p-1\right]}
=\frac{O(\varepsilon^{\frac4p}\log^{p-1}|
\log\varepsilon|)}{|\log\varepsilon|^{p-1}}
\exp\left[O\left(\frac{\log^{p}|
\log\varepsilon|}{|\log\varepsilon|^{p-1}}\right)\right],
\endaligned
\end{equation}
which, together with  (\ref{2.6}) and (\ref{2.35}), easily yields
\begin{eqnarray}\label{2.38}
\left|\left(
\sum\limits_{i=1}^m\frac{\mu_i^\sigma}{(\mu_i+|y-\xi'_i|)^{2+\sigma}}
+\varepsilon^2
\right)^{-1}
E_{\xi'}(y)
\right|
\leq C
\varepsilon^{\frac{4-2p}{p}}\left(\frac{\log|\log\varepsilon|}{|\log\varepsilon|}
\right)^{p-1}
\exp\left[
O\left(
\frac{\log^{p}|
\log\varepsilon|}{|\log\varepsilon|^{p-1}}
\right)
\right]
=o\left(\frac{1}{\gamma^{4p}}
\right).
\end{eqnarray}
Let us  fix an index $i\in\{1,\ldots,m\}$
and the region $|y-\xi'_i|\leq1/(\varepsilon^\theta|\log\varepsilon|^{2\kappa})$ with
any $\theta<1$ but close enough to $1$.
From (\ref{2.21}), (\ref{2.28}) and  Taylor expansion we have that
in the ball  $|y-\xi'_i|<\mu_i|\log\varepsilon|^\tau$ with  $\tau\geq 10$ large but fixed,
$$
\aligned
\left(1+\frac{V_{\xi'}}{p\gamma^p}\right)^{p-1}
=&\,1+\frac{p-1}{p}\frac1{\gamma^p}\underbrace{\omega_{\mu_i}(y-\xi'_i)}\limits_{A_1}
+\left(\frac{p-1}{p}\right)^2\frac1{\gamma^{2p}}\underbrace{
\left[
\omega^{1}_{\mu_i}+\frac{p-2}{2(p-1)}(\omega_{\mu_i})^2\right]
(y-\xi'_i)}\limits_{A_2}\\
&+\left(\frac{p-1}{p}\right)^3\frac1{\gamma^{3p}}\underbrace{\left[
\omega^{2}_{\mu_i}+\frac{p-2}{p-1}\omega_{\mu_i}\omega^1_{\mu_i}
+\frac{(p-2)(p-3)}{6(p-1)^2}(\omega_{\mu_i})^3
\right]
(y-\xi'_i)}\limits_{A_3}
\\
&+\left(\frac{p-1}{p}\right)^4\frac1{\gamma^{4p}}\underbrace{\left[
\omega^{3}_{\mu_i}+\frac{p-2}{2(p-1)}(\omega^1_{\mu_i})^2
+\frac{p-2}{p-1}\omega_{\mu_i}\omega^2_{\mu_i}
+\frac{(p-2)(p-3)}{2(p-1)^2}(\omega_{\mu_i})^2\omega^1_{\mu_i}
\right.}\limits_{A_4}\\
&\underbrace{\left.
+\,\frac{(p-2)(p-3)(p-4)}{24(p-1)^3}(\omega_{\mu_i})^4
\right](y-\xi'_i)}\limits_{A_4}
+O\left(\frac{\log^\beta(\mu_i+|y-\xi'_i|)}{\gamma^{5p}}\right),
\endaligned
$$
%$$
%\aligned
%\left(1+\frac{V_{\xi'}}{p\gamma^p}\right)^{p-1}=&\,1+\frac{p-1}{p}\frac1{\gamma^p}\omega_{\mu_i}(y-\xi'_i)
%+\left(\frac{p-1}{p}\right)^2\frac1{\gamma^{2p}}\left[
%\omega^{1}_{\mu_i}+\frac{p-2}{2(p-1)}(\omega_{\mu_i})^2\right](y-\xi'_i)\\
%&+\left(\frac{p-1}{p}\right)^3\frac1{\gamma^{3p}}\left[
%\omega^{2}_{\mu_i}+\frac{p-2}{p-1}\omega_{\mu_i}\omega^1_{\mu_i}
%+\frac{(p-2)(p-3)}{6(p-1)^2}(\omega_{\mu_i})^3
%\right](y-\xi'_i)
%\\
%&+\left(\frac{p-1}{p}\right)^4\frac1{\gamma^{4p}}\left[
%\omega^{3}_{\mu_i}+\frac{p-2}{2(p-1)}(\omega^1_{\mu_i})^2
%+\frac{p-2}{p-1}\omega_{\mu_i}\omega^2_{\mu_i}
%+\frac{(p-2)(p-3)}{2(p-1)^2}(\omega_{\mu_i})^2\omega^1_{\mu_i}
%\right.\\
%&\left.
%+\,\frac{(p-2)(p-3)(p-4)}{24(p-1)^3}(\omega_{\mu_i})^4
%\right](y-\xi'_i)
%+O\left(\frac{\log^\beta(\mu_i+|y-\xi'_i|)}{\gamma^{5p}}\right),
%\endaligned
%$$
and
\begin{eqnarray}\label{2.380}
\gamma^p\left[\left(1+\frac{V_{\xi'}}{p\gamma^p}\right)^{p}-1\right]
=\omega_{\mu_i}(y-\xi'_i)+\frac{p-1}{p}\frac1{\gamma^p}
\underbrace{\left[\omega^{1}_{\mu_i}+\frac{1}{2}(\omega_{\mu_i})^2\right](y-\xi'_i)}\limits_{B_1}
\qquad\qquad\qquad\qquad\qquad\qquad\qquad\qquad\,\,\,\,
\nonumber\\
+\left(\frac{p-1}{p}\right)^2\frac1{\gamma^{2p}}\underbrace{
\left[
\omega^{2}_{\mu_i}+\omega_{\mu_i}\omega^{1}_{\mu_i}+\frac{p-2}{6(p-1)}(\omega_{\mu_i})^3\right]
(y-\xi'_i)}\limits_{B_2}
\qquad\qquad\qquad\qquad\qquad\qquad\qquad\,\,
\qquad\qquad
\nonumber\\
+\left(\frac{p-1}{p}\right)^3\frac1{\gamma^{3p}}\underbrace{
\left[
\omega^{3}_{\mu_i}
+\frac12(\omega^{1}_{\mu_i})^2
+
\omega_{\mu_i}\omega^{2}_{\mu_i}
+\frac{p-2}{2(p-1)}(\omega_{\mu_i})^2\omega^{1}_{\mu_i}
+\frac{(p-2)(p-3)}{24(p-1)^2}(\omega_{\mu_i})^4
\right]
(y-\xi'_i)}\limits_{B_3}
\qquad
\nonumber\\
+\left(\frac{p-1}{p}\right)^4\frac1{\gamma^{4p}}
\underbrace{\left[
\omega^{4}_{\mu_i}
+
\omega^1_{\mu_i}\omega^{2}_{\mu_i}
+\frac{p-2}{2(p-1)}(\omega_{\mu_i})^2\omega^{2}_{\mu_i}
+\frac{p-2}{2(p-1)}\omega_{\mu_i}(\omega^{1}_{\mu_i})^2
+\frac{(p-2)(p-3)}{6(p-1)^2}(\omega_{\mu_i})^3\omega^{1}_{\mu_i}
\right.}
\limits_{B_4}
\nonumber\\
\underbrace{+
\left.\frac{(p-2)(p-3)(p-4)}{120(p-1)^3}(\omega_{\mu_i})^5
\right](y-\xi'_i)}\limits_{B_4}
+
O\left(\frac{\log^\beta(\mu_i+|y-\xi'_i|)}{\gamma^{5p}}\right),
\qquad\qquad\qquad\qquad\qquad\qquad\qquad\,\,\,\,
\end{eqnarray}
%$$
%\aligned
%\gamma^p\left[\left(1+\frac{V_{\xi'}}{p\gamma^p}\right)^{p}-1\right]
%=&\,\omega_{\mu_i}(y-\xi'_i)+\frac{p-1}{p}\frac1{\gamma^p}\left[\omega^{1}_{\mu_i}+\frac{1}{2}(\omega_{\mu_i})^2\right](y-\xi'_i)
%\\
%+&\left(\frac{p-1}{p}\right)^2\frac1{\gamma^{2p}}\left[
%\omega^{2}_{\mu_i}+\omega_{\mu_i}\omega^{1}_{\mu_i}+\frac{p-2}{6(p-1)}(\omega_{\mu_i})^3\right](y-\xi'_i)
%\\
%+&\left(\frac{p-1}{p}\right)^3\frac1{\gamma^{3p}}\left[
%\omega^{3}_{\mu_i}
%+\frac12(\omega^{1}_{\mu_i})^2
%+
%\omega_{\mu_i}\omega^{2}_{\mu_i}
%+\frac{p-2}{2(p-1)}(\omega_{\mu_i})^2\omega^{1}_{\mu_i}
%+\frac{(p-2)(p-3)}{24(p-1)^2}(\omega_{\mu_i})^4
%\right](y-\xi'_i)\\
%+&\left(\frac{p-1}{p}\right)^4\frac1{\gamma^{4p}}
%\left[
%\omega^{4}_{\mu_i}
%+
%\omega^1_{\mu_i}\omega^{2}_{\mu_i}
%+\frac{p-2}{2(p-1)}(\omega_{\mu_i})^2\omega^{2}_{\mu_i}
%+\frac{p-2}{2(p-1)}\omega_{\mu_i}(\omega^{1}_{\mu_i})^2
%+\frac{(p-2)(p-3)}{6(p-1)^2}(\omega_{\mu_i})^3\omega^{1}_{\mu_i}
%\right.\\
%+&\,
%\left.\frac{(p-2)(p-3)(p-4)}{120(p-1)^3}(\omega_{\mu_i})^5
%\right]
%+
%O\left(\frac{\log^\beta(\mu_i+|y-\xi'_i|)}{\gamma^{5p}}\right),
%\endaligned
%$$
where $\beta>1$ is large but  fixed, independent of $\varepsilon$.
Then
\begin{eqnarray}\label{2.381}
e^{\gamma^p\left[\left(1+\frac{V_{\xi'}}{p\gamma^p}\right)^{p}-1\right]}
=e^{\omega_{\mu_i}\left(y-\xi'_i\right)}
\left\{1+\frac{p-1}{p}\frac1{\gamma^p}
B_1
+\left(\frac{p-1}{p}\right)^2\frac1{\gamma^{2p}}
\left[B_2+\frac12(B_1)^2\right]
+\left(\frac{p-1}{p}\right)^3\frac1{\gamma^{3p}}
\big[B_3+B_1B_2
\right.
\quad
\nonumber\\[1mm]
\left.\left.
+\frac16(B_1)^3\right]+\left(\frac{p-1}{p}\right)^4\frac1{\gamma^{4p}}\left[B_4
+\frac12(B_2)^2+B_1B_3+\frac12(B_1)^2B_2+\frac1{24}(B_1)^4\right]
+O\left(\frac{\log^\beta(\mu_i+|y-\xi'_i|)}{\gamma^{5p}}\right)
\right\}.
\end{eqnarray}
Furthermore, by the definition of $f(\cdot)$ in (\ref{2.27}) and
the definitions of $f^j_{\mu_i}$,
$j=1,2,3,4$ in (\ref{2.10})-(\ref{2.121}) we can compute
\begin{eqnarray}\label{2.39}
f(V_{\xi'})=e^{\omega_{\mu_i}\left(y-\xi'_i\right)}\,
\left\{1+\frac{p-1}{p}\frac1{\gamma^p}
\underbrace{\left(A_1+B_1\right)}\limits_{=\big(
\omega^1_{\mu_i}-f^1_{\mu_i}
\big)\left(y-\xi'_i\right)}
+\left(\frac{p-1}{p}\right)^2\frac1{\gamma^{2p}}
\underbrace{\left[A_2+A_1B_1+
B_2+\frac12(B_1)^2\right]}\limits_{=\big(
\omega^2_{\mu_i}-f^2_{\mu_i}
\big)\left(y-\xi'_i\right)}
\right.
\,\qquad\qquad\qquad\quad\,\,\,\,
\nonumber\\[1mm]
+\left(\frac{p-1}{p}\right)^3\frac1{\gamma^{3p}}
\underbrace{\left[
A_3
+A_2B_1
+A_1\left(B_2+\frac12(B_1)^2\right)+
B_3+B_1B_2
+\frac16(B_1)^3\right]}\limits_{=\big(
\omega^3_{\mu_i}-f^3_{\mu_i}
\big)\left(y-\xi'_i\right)}
\,\,\,+\,\left(\frac{p-1}{p}\right)^4\frac1{\gamma^{4p}}\times
\qquad\qquad\,
\nonumber\\[1mm]
\underbrace{\left[A_4
+A_3B_1+A_2\left(B_2+
\frac{(B_1)^2}2\right)
+A_1\left(B_3+B_1B_2
+\frac{(B_1)^3}6\right)+B_4
+\frac{(B_2)^2}2+B_1B_3+\frac{(B_1)^2B_2}2+\frac{(B_1)^4}{24}\right]}
\limits_{=\big(
\omega^4_{\mu_i}-f^4_{\mu_i}
\big)\left(y-\xi'_i\right)}
\nonumber\\[1mm]
\left.
+\,O\left(\frac{\log^\beta(\mu_i+|y-\xi'_i|)}{\gamma^{5p}}\right)
\right\}.
\qquad\qquad\qquad\qquad\qquad\qquad\qquad\qquad\qquad\qquad
\qquad\qquad\qquad\qquad\qquad\qquad\qquad
\end{eqnarray}
From (\ref{2.34}) and (\ref{2.39})   we find that
in the region
$|y-\xi'_i|<\mu_i|\log\varepsilon|^\tau$,
\begin{equation*}\label{2.40}
\aligned
E_{\xi'}=-\Delta_{a(\varepsilon y)}V_{\xi'}+\varepsilon^2V_{\xi'}-f(V_{\xi'})+p\gamma^{p}\varepsilon^2
=e^{\omega_{\mu_i}\left(y-\xi'_i\right)}
O\left(\frac{\log^\beta(\mu_i+|y-\xi'_i|)}{\gamma^{5p}}\right),
\endaligned
\end{equation*}
and by (\ref{2.7}) and (\ref{2.6}),
\begin{eqnarray}\label{2.41}
\left|\left(
\sum\limits_{i=1}^m\frac{\mu_i^\sigma}{(\mu_i+|y-\xi'_i|)^{2+\sigma}}
+\varepsilon^2
\right)^{-1}
E_{\xi'}(y)
\right|
\leq \frac{C}{\gamma^{5p}}\frac{\,\mu_i^{2-\sigma}\log^\beta\big(\mu_i+|y-\xi'_i|\big)\,}{\big(\mu_i+|y-\xi'_i|\big)^{2-\sigma}}
\leq \frac{C}{\gamma^{4p}}.
\end{eqnarray}
As in the remaining region
$\mu_i|\log\varepsilon|^\tau\leq
|y-\xi'_i|\leq1/(\varepsilon^\theta|\log\varepsilon|^{2\kappa})$ with
any $\theta<1$ but close enough to $1$,
by  (\ref{2.10})-(\ref{2.13}) and (\ref{2.34}) we  find
that there exists a constant $D>0$, independent of every $\theta<1$, such that
\begin{equation}\label{2.42}
\aligned
\big|-\Delta_{a(\varepsilon y)}V_{\xi'}+\varepsilon^2V_{\xi'}+p\gamma^{p}\varepsilon^2
\big|\leq D |\log\varepsilon|^4 e^{\omega_{\mu_i}(y-\xi'_i)}.
\endaligned
\end{equation}
On the other hand, in the same region,
by  (\ref{2.6}),  (\ref{2.13}), (\ref{2.21}) and (\ref{2.28}) we have that
$$
\aligned
1+\frac{V_{\xi'}}{p\gamma^p}=&1+\frac{1}{p\gamma^p}\left\{
\left[
1-
\frac{1}{4}
\sum_{j=1}^4\left(\frac{p-1}{p}\right)^j\frac{D^j_{\mu_i}}{\gamma^{jp}}
\right]\omega_{\mu_i}\left(y-\xi'_i\right)
+\left[\sum_{j=1}^4\left(\frac{p-1}{p}\right)^j\frac{D^j_{\mu_i}}{\gamma^{jp}}
\right]\frac14\log\left(\frac{8}{\mu_i^2}\right)
+O\left(\frac{\mu_i}{|y-\xi'_i|}\right)\right\}
\\
=&1+\frac{1}{p\gamma^p}\left\{
\left[
1+
O\left(\frac{\log|\log\varepsilon|}{|\log\varepsilon|}\right)
\right]\omega_{\mu_i}\left(y-\xi'_i\right)
+O\left(\frac{\log^2|\log\varepsilon|}{|\log\varepsilon|}\right)\right\},
\endaligned
$$
then
\begin{equation}\label{2.45}
\aligned
1-\theta+2\kappa\frac{\log|\log\varepsilon|}{|\log\varepsilon|}
+
O\left(\frac{1}{|\log\varepsilon|}\right)
\leq
1+\frac{V_{\xi'}}{p\gamma^p}
\leq
1-\tau\frac{\log|\log\varepsilon|}{|\log\varepsilon|}
+
O\left(\frac{1}{|\log\varepsilon|}\right).
\endaligned
\end{equation}
%Obviously,
%$$
%\aligned
%\left(1+\frac{V_{\xi'}}{p\gamma^p}\right)^{p-1}
%=O\left(1\right).
%\endaligned
%$$
%Moreover, by the Taylor expansion   we
%find that there exist a constant $D>0$, independent of $\theta$, such that
%\begin{equation}\label{2.46}
%\aligned
%e^{\gamma^p\left[\left(1+\frac{V_{\xi'}}{p\gamma^p}\right)^p-1\right]}\leq
%De^{\left[
%1-
%\frac{1}{4}
%\sum_{j=1}^4\left(\frac{p-1}{p}\right)^j\frac{D^j_{\mu_i}}{\gamma^{jp}}
%\right]
%\omega_{\mu_i}\left(y-\xi'_i\right)
%},
%\endaligned
%\end{equation}
%and then
%\begin{equation}\label{2.43}
%\aligned
%f(V_{\xi'})=\left(1+\frac{V_{\xi'}}{p\gamma^p}\right)^{p-1}e^{\gamma^p\left[\left(1+\frac{V_{\xi'}}{p\gamma^p}\right)^p-1\right]}
%=O\left(1\right)
%e^{\left[
%1+
%O\left(\frac{\log|\log\varepsilon|}{|\log\varepsilon|}\right)
%\right]
%\omega_{\mu_i}\left(y-\xi'_i\right)
%}.
%\endaligned
%\end{equation}
Furthermore, by the Taylor expansion   we
find that there exists a constant $D>0$, independent of every $\theta<1$, such that
\begin{equation}\label{2.46}
\aligned
e^{\gamma^p\left[\left(1+\frac{V_{\xi'}}{p\gamma^p}\right)^p-1\right]}\leq
De^{\left[
1-
\frac{1}{4}
\sum_{j=1}^4\left(\frac{p-1}{p}\right)^j\frac{D^j_{\mu_i}}{\gamma^{jp}}
\right]
\omega_{\mu_i}\left(y-\xi'_i\right)
},
\endaligned
\end{equation}
and
\begin{equation}\label{2.43}
\aligned
\left(1+\frac{V_{\xi'}}{p\gamma^p}\right)^{p-1}
\leq D\left(1+\frac{\log^{p-1}|\log\varepsilon|}{|\log\varepsilon|^{p-1}}\right).
\endaligned
\end{equation}
Hence in the region $\mu_i|\log\varepsilon|^\tau\leq
|y-\xi'_i|\leq1/(\varepsilon^\theta|\log\varepsilon|^{2\kappa})$ with
any $\theta<1$ but close enough to $1$, by (\ref{2.42}), (\ref{2.46}) and (\ref{2.43}),
$$
\aligned
\left|\left(
\sum\limits_{i=1}^m\frac{\mu_i^\sigma}{(\mu_i+|y-\xi'_i|)^{2+\sigma}}
+\varepsilon^2
\right)^{-1}
E_{\xi'}(y)
\right|
\leq C\left[|\log\varepsilon|^4
\left|\frac{y-\xi'_i}{\mu_i}\right|^{\sigma-2}
+\left(\frac{1}{\mu_i^2}\right)^{
O\left(\frac{\log|\log\varepsilon|}{|\log\varepsilon|}\right)}
\left|\frac{y-\xi'_i}{\mu_i}\right|^{\sigma-1+
O\left(\frac{\log|\log\varepsilon|}{|\log\varepsilon|}\right)}\right]
=o\left(\frac{1}{\gamma^{4p}}
\right),
\endaligned
$$
which, together with (\ref{2.38}) and (\ref{2.41}), implies the validity of
 estimate  (\ref{2.33}).
\end{proof}

\vspace{1mm}

\section{Analysis of the linearized operator}
In this section we perform the solvability theory for the linear operator $\mathcal{L}$
under the weighted $L^{\infty}$-norm introduced in (\ref{2.32}),
uniformly on $\xi\in\mathcal{O}_\varepsilon$.
Notice that
$\mathcal{L}(\phi)=-\Delta_{a(\varepsilon y)}\phi+\varepsilon^2\phi-W_{\xi'}\phi$,
where $W_{\xi'}=f'(V_{\xi'})$.
As in Proposition 2.3, we have  the following asymptotical expansions with respect to   $W_{\xi'}$
and $f''(V_{\xi'})$, respectively.

\vspace{1mm}
\vspace{1mm}
\vspace{1mm}
\vspace{1mm}

\noindent{\bf Proposition 3.1.}\,\,{\it
There exists a constant $D_0>0$ such that
for any $\xi=(\xi_1,\ldots,\xi_m)\in\mathcal{O}_{\varepsilon}$ and
for any $\varepsilon$ small enough,
\begin{equation}\label{3.1}
\aligned
\big|W_{\xi'}(y)\big|
\leq D_0\sum_{i=1}^me^{\left[
1-
\frac{1}{4}
\sum_{j=1}^4\left(\frac{p-1}{p}\right)^j\frac{D^j_{\mu_i}}{\gamma^{jp}}
\right]
\omega_{\mu_i}\left(y-\xi'_i\right)
}
\qquad
\textrm{and}
\qquad
\big|f''(V_{\xi'})\big|\leq D_0\sum_{i=1}^me^{\left[
1-
\frac{1}{4}
\sum_{j=1}^4\left(\frac{p-1}{p}\right)^j\frac{D^j_{\mu_i}}{\gamma^{jp}}
\right]
\omega_{\mu_i}\left(y-\xi'_i\right)
},
\endaligned
\end{equation}
uniformly in the region   $\mu_i|\log\varepsilon|^\tau\leq|y-\xi'_i|\leq1/(\varepsilon^\theta|\log\varepsilon|^{2\kappa})$ with
any $\theta<1$ but close enough to $1$.
While if $|y-\xi'_i|<\mu_i|\log\varepsilon|^\tau$ with  $\tau\geq10$  large but fixed,
then
\begin{equation}\label{3.2}
\aligned
W_{\xi'}(y)=\frac{8\mu_i^2}{(\mu_i^2+|y-\xi'_i|^2)^2}
\left\{1+\frac{p-1}{p}\frac1{\gamma^p}
\left[1+\omega^{1}_{\mu_i}+\frac{1}{2}(\omega_{\mu_i})^2+2\omega_{\mu_i}\right]\left(y-\xi'_i\right)
+O\left(\frac{\log^4(\mu_i+|y-\xi'_i|)}{\gamma^{2p}}\right)
\right\}.
\endaligned
\end{equation}
In addition,
\begin{equation}\label{3.3}
\aligned
\|W_{\xi'}\|_{*}\leq C
\,\,\qquad\ \,\,
\textrm{and}
\,\,\qquad\ \,\,
\|f''(V_{\xi'})\|_{*}\leq C.
\endaligned
\end{equation}}\noindent~
\begin{proof}
For the sake of simplicity, we consider the estimates for the potential $W_{\xi'}$ only.
By (\ref{2.27}) we can compute
$$
\aligned
W_{\xi'}=f'(V_{\xi'})=\frac{p-1}{p}\frac1{\gamma^p}\left(1+\frac{V_{\xi'}}{p\gamma^p}\right)^{p-2}
e^{\gamma^p\left[\left(1+\frac{V_{\xi'}}{p\gamma^p}\right)^p-1\right]}+
\left(1+\frac{V_{\xi'}}{p\gamma^p}\right)^{2(p-1)}
e^{\gamma^p\left[\left(1+\frac{V_{\xi'}}{p\gamma^p}\right)^p-1\right]}
:=I+J.
\endaligned
$$
If $|y-\xi'_i|<\mu_i|\log\varepsilon|^\tau$ with any $i\in\{1,\ldots,m\}$ and $\tau\geq10$ large but fixed,
by using (\ref{2.381})
and Taylor expansion we obtain
$$
\aligned
I
=& e^{\omega_{\mu_i}\left(y-\xi'_i\right)}\left\{1+\frac{p-1}{p}\frac1{\gamma^p}\left[\omega^{1}_{\mu_i}+\frac{1}{2}(\omega_{\mu_i})^2\right]\left(y-\xi'_i\right)
+O\left(\frac{\log^4(\mu_i+|y-\xi'_i|)}{\gamma^{2p}}\right)
\right\}\\[1mm]
&\times\frac{p-2}{p}\frac1{\gamma^p}\left[\frac{p-1}{p-2}+\frac{p-1}{p}\frac1{\gamma^p}\omega_{\mu_i}\left(y-\xi'_i\right)
+O\left(\frac{\log^2(\mu_i+|y-\xi'_i|)}{\gamma^{2p}}\right)
\right],
\endaligned
$$
and
$$
\aligned
J
=&e^{\omega_{\mu_i}\left(y-\xi'_i\right)}\left\{1+\frac{p-1}{p}\frac1{\gamma^p}
\left[\omega^{1}_{\mu_i}+\frac{1}{2}(\omega_{\mu_i})^2\right]\left(y-\xi'_i\right)
+O\left(\frac{\log^4(\mu_i+|y-\xi'_i|)}{\gamma^{2p}}\right)
\right\}
\\
&
\times
\left[1+\frac{p-1}{p}\frac2{\gamma^p}\omega_{\mu_i}\left(y-\xi'_i\right)
+O\left(\frac{\log^2(\mu_i+|y-\xi'_i|)}{\gamma^{2p}}\right)
\right],
\endaligned
$$
and hence
\begin{equation}\label{3.5}
\aligned
W_{\xi'}(y)=e^{\omega_{\mu_i}\left(y-\xi'_i\right)}
\left\{1+\frac{p-1}{p}\frac1{\gamma^p}
\left[1+\omega^{1}_{\mu_i}+\frac{1}{2}(\omega_{\mu_i})^2+2\omega_{\mu_i}\right]\left(y-\xi'_i\right)
+O\left(\frac{\log^4(\mu_i+|y-\xi'_i|)}{\gamma^{2p}}\right)
\right\}.
\endaligned
\end{equation}
While if  $\mu_i|\log\varepsilon|^\tau\leq
|y-\xi'_i|\leq1/(\varepsilon^\theta|\log\varepsilon|^{2\kappa})$
with  any $\theta<1$ but close enough to $1$, by  (\ref{2.45})
we find
$$
\aligned
\left(1+\frac{V_{\xi'}(y)}{p\gamma^p}\right)^{p-2}
=O\left(1\right)
\qquad\qquad
\textrm{and}
\qquad\qquad
\left(1+\frac{V_{\xi'}(y)}{p\gamma^p}\right)^{2(p-1)}
=O\left(1\right),
\endaligned
$$
and  by (\ref{2.46}),
\begin{equation}\label{3.4}
\aligned
\big|W_{\xi'}(y)\big|\leq
Ce^{\gamma^p\left[\left(1+\frac{V_{\xi'}}{p\gamma^p}\right)^p-1\right]}=
O\left(1\right)
e^{\left[
1-
\frac{1}{4}
\sum_{j=1}^4\left(\frac{p-1}{p}\right)^j\frac{D^j_{\mu_i}}{\gamma^{jp}}
\right]
\omega_{\mu_i}\left(y-\xi'_i\right)
}.
\endaligned
\end{equation}
Additionally, if
$|y-\xi'_i|\geq1/(\varepsilon|\log\varepsilon|^{2\kappa})$
for all $i=1,\ldots,m$, by  (\ref{2.36}) we deduce
$$
\aligned
I=\frac{O(\varepsilon^{\frac4p}\log^{p-2}|
\log\varepsilon|)}{|\log\varepsilon|^{p-1}}\exp\left[O\left(\frac{\log^{p}|
\log\varepsilon|}{|\log\varepsilon|^{p-1}}\right)
\right],
\,\quad
\quad\,J=\frac{O(\varepsilon^{\frac4p}\log^{2(p-1)}|
\log\varepsilon|)}{|\log\varepsilon|^{2(p-1)}}\exp\left[O\left(\frac{\log^{p}|
\log\varepsilon|}{|\log\varepsilon|^{p-1}}\right)\right],
\endaligned
$$
and so
\begin{equation}\label{3.5a}
\aligned
\big|W_{\xi'}(y)\big|=
\left(\frac{\log^{p-2}|
\log\varepsilon|}{|\log\varepsilon|^{p-1}}
+\frac{\log^{2(p-1)}|
\log\varepsilon|}{|\log\varepsilon|^{2(p-1)}}\right)O(\varepsilon^{\frac4p})
\exp\left[O\left(\frac{\log^{p}|
\log\varepsilon|}{|\log\varepsilon|^{p-1}}\right)\right].
\endaligned
\end{equation}
Jointing together  (\ref{3.5})-(\ref{3.5a}) and the definition of $\left\|\cdot\right\|_{*}$ in
(\ref{2.32}), we obtain the first estimate in  (\ref{3.3}).
\end{proof}

\vspace{1mm}
\vspace{1mm}

Given $h\in L^\infty(\Omega_\varepsilon)$ and points
$\xi=(\xi_1,\ldots,\xi_m)\in\mathcal{O}_{\varepsilon}$,
we consider
the following linear problem of finding a function $\phi$ and scalars
$c_{ij}\in\mathbb{R}$, $i=1,\ldots,m$, $j=1,J_i$,
such that
\begin{equation}\label{3.9}
\left\{\aligned
&\mathcal{L}(\phi)=-\Delta_{a(\varepsilon y)}\phi+\varepsilon^2\phi-W_{\xi'}\phi=h
+\frac1{a(\varepsilon y)}\sum\limits_{i=1}^m\sum\limits_{j=1}^{J_i}c_{ij}\chi_i\,Z_{ij}\,\,\ \,
\,\textrm{in}\,\,\,\,\,\,\Omega_\varepsilon,\\
&\frac{\partial\phi}{\partial\nu}=0\,\,\,\,\,\,\,\,
\ \ \ \ \ \ \ \ \ \ \ \ \ \,\,
\qquad\qquad\qquad\quad\qquad\qquad\qquad\qquad
\ \,\ \ \ \,\,\,\ \,
\ \,\textrm{on}\,\,\,\,\partial\Omega_{\varepsilon},\\[1mm]
&\int_{\Omega_\varepsilon}\chi_i\,Z_{ij}\phi=0
\,\qquad\qquad\qquad\quad
\qquad\qquad\qquad\forall\,\,i=1,\ldots,m,\,\,\,j=1, J_i,
\endaligned\right.
\end{equation}
where
$J_i=2$ if $i=1,\ldots,l$ while $J_i=1$ if $i=l+1,\ldots,m$,
and $Z_{ij}$, $\chi_i$, are  defined as follows:
let $R_0>0$ be a large but fixed number and
$\chi:\mathbb{R}\rightarrow[0,1]$ be
a smooth, non-increasing cut-off function  such
that  $\chi(r)=1$ if $r\leq
R_0$,  and $\chi(r)=0$ if $r\geq R_0+1$.
Set
\begin{equation}\label{3.6}
\aligned
Z_{0}(z)=\frac{|z|^2-1}{|z|^2+1},
\,\,\quad\qquad\quad\,\,
Z_{j}(z)=\frac{z_j}{|z|^2+1},\,\,\,\,j=1,\,2.
\endaligned
\end{equation}
For each $i=1,\ldots,l$, we
have $\xi_i\in\Omega$
and define
\begin{equation}\label{3.10}
\aligned
\chi_i(y)=\chi\left(
\frac{|y-\xi'_i|}{\mu_i}
\right),
\,\qquad\qquad\,
Z_{ij}(y)=\frac{1}{\mu_i}Z_j\left(\frac{y-\xi_i'}{\mu_i}\right),
\,\quad\,j=0,1,2.
\endaligned
\end{equation}
For each $i=l+1,\ldots,m$, we have  $\xi_i\in\po$ and
define a rotation map  $A_i: \mathbb{R}^2\mapsto\mathbb{R}^2$ such that
$A_i\nu_{\Omega}(\xi_i)=\nu_{\mathbb{R}_+^2}(0)$.
Let $\mathcal{G}(x_1)$ be the defining function
for the boundary $A_i(\po-\{\xi_i\})$ in a small
neighborhood  of the origin,
that is, there exist $R_1>0$, $\delta>0$ small and
a smooth function $\mathcal{G}:(-R_1,R_1)\mapsto\mathbb{R}$
satisfying $\mathcal{G}(0)=0$, $\mathcal{G}'(0)=0$ and such that
$A_i(\Omega-\{\xi_i\})\cap B_\delta(0,0)=\{(x_1,x_2):\,-R_1<x_1<R_1,\,x_2>\mathcal{G}(x_1)\}\cap B_\delta(0,0)$.
Then we consider the flattening change of variables
$F_i: B_\delta(0,0)\cap\overline{A_i(\Omega-\{\xi_i\})}\mapsto\mathbb{R}^2$
defined by
\begin{equation*}\label{3.11}
\aligned
F_i=(F_{i1}, F_{i2}),
\,\qquad\,\textrm{where}
\,\quad\,F_{i1}=x_1+\frac{x_2-\mathcal{G}(x_1)}{\,1+|\mathcal{G}'(x_1)|^2\,}\mathcal{G}'(x_1),
\ \qquad\,F_{i2}=x_2-\mathcal{G}(x_1).
\endaligned
\end{equation*}
Then for each $i=l+1,\ldots,m$, we set
\begin{equation*}\label{3.12}
\aligned
F_i^\varepsilon(y)=\frac{1}{\varepsilon}F_i\big(A_i(\varepsilon y-\xi_i)\big),
\endaligned
\end{equation*}
and define
\begin{equation}\label{3.13}
\aligned
\chi_i(y)=\chi\left(
\frac{1}{\mu_i}|F_i^\varepsilon(y)|
\right),
\,\qquad\qquad\,
Z_{ij}(y)=\frac{1}{\mu_i}Z_{j}\left(\frac{1}{\mu_i}F_i^\varepsilon(y)\right),
\,\quad\,j=0,1.
\endaligned
\end{equation}
Note that $F^\varepsilon_i$, $i=l+1,\ldots,m$,
preserves the  homogeneous Neumann boundary condition. Moreover,
\begin{equation}\label{3.14}
\aligned
\Delta_{a(\varepsilon y)}Z_{i0}+\frac{8\mu_i^2}{(\mu_i^2+|y-\xi'_i|^2)^2}Z_{i0}
=O\left(\frac{\varepsilon\mu_i}{
(\mu_i+|y-\xi'_i|)^3}\right),
\,\,\,
\,\ \,\forall\,\,i=1,\ldots,m.
\endaligned
\end{equation}

\vspace{1mm}
\vspace{1mm}
\vspace{1mm}
\vspace{1mm}

\noindent {\bf Proposition 3.2.}\,\,{\it
Let $m$ be a positive integer.
Then there exist constants $C>0$  and $\varepsilon_0>0$ such
that for any  $0<\varepsilon<\varepsilon_0$,    any points
$\xi=(\xi_1,\ldots,\xi_m)\in\mathcal{O}_{\varepsilon}$ and
any $h\in L^\infty(\Omega_\varepsilon)$,
there is a unique solution $\phi\in
L^\infty(\Omega_\varepsilon)$
and   $c_{ij}\in\mathbb{R}$,  $i=1,\ldots,m$, $j=1,J_i$ to problem {\upshape(\ref{3.9})}.
Moreover,
\begin{equation*}\label{3.15}
\aligned
\|\phi\|_{L^{\infty}(\Omega_\varepsilon)}\leq C|\log\varepsilon|\,\|h\|_{*}
\,\,\qquad\,\quad\,\,
\textrm{and}
\,\,\qquad\,\quad\,\,
\sum_{i=1}^m\sum_{j=1}^{J_i}\mu_i|c_{ij}|\leq C \|h\|_{*}.
\endaligned
\end{equation*}
}\indent The proof of this result  will be split into a series of lemmas which we state and prove as follows.

\vspace{1mm}
\vspace{1mm}
\vspace{1mm}

\noindent{\bf Lemma 3.3.}\,\,{\it There exist constants $R_1>0$ and $C>0$,
independent of $\varepsilon$, such that
for any sufficiently small
$\varepsilon$,  any points $\xi=(\xi_1,\ldots,\xi_m)\in\mathcal{O}_\varepsilon$ and any $\sigma\in(0,1)$,
there is a function
$$
\aligned
\psi:\,\,\Omega_\varepsilon\setminus\bigcup_{i=1}^mB_{R_1\mu_i}(\xi'_i)\,\,\mapsto\mathbb{R}
\endaligned
$$
smooth  and
positive so that
$$
\aligned
\mathcal{L}(\psi)=-\Delta_{a(\varepsilon y)}\psi+\varepsilon^2\psi-W_{\xi'}\psi&\geq\sum_{i=1}^m\frac{\mu_i^\sigma}{|y-\xi'_i|^{2+\sigma}}+\varepsilon^2
\,\,\,\quad\,\,
\textrm{in}\,\,\,\,\,\,\,\Omega_\varepsilon\setminus\bigcup_{i=1}^mB_{R_1\mu_i}(\xi'_i),\\
\frac{\partial\psi}{\partial\nu}&\geq0
\,\,\,\qquad\qquad\qquad\qquad\quad\,\,
\textrm{on}\,\,\,\,\partial\Omega_\varepsilon\setminus\bigcup_{i=1}^mB_{R_1\mu_i}(\xi'_i),\\
\psi&>0
\,\,\,\qquad\qquad\qquad\qquad\quad\,\,
\textrm{in}\,\,\,\,\,\,\,
\Omega_\varepsilon\setminus\bigcup_{i=1}^mB_{R_1\mu_i}(\xi'_i),\\
\psi&\geq1
\,\,\,\qquad\qquad\qquad\qquad\quad\,\,
\textrm{on}\,\,\,\,
\Omega_\varepsilon\cap\left(\bigcup_{i=1}^m\partial B_{R_1\mu_i}(\xi'_i)\right).
\endaligned
$$
Moreover, $\psi$ is   uniformly bounded, i.e.
$$
\aligned
1<\psi\leq C\,\,\,\quad\textrm{in}\,\,\,
\,\Omega_\varepsilon\setminus\bigcup_{i=1}^mB_{R_1\mu_i}(\xi'_i).
\endaligned
$$
}

\begin{proof}
Let us take
$$
\aligned
\psi=\sum_{i=1}^m
\left(1-
\frac{\mu_i^\sigma}{|y-\xi'_i|^\sigma}
\right)+
C_1
\Psi_0(y),
\endaligned
$$
where $\Psi_0$ is the uniformly bounded
 solution of
$$
\aligned
\left\{
\aligned
&-\Delta_{a(\varepsilon y)}\Psi_0+\varepsilon^2\Psi_0=\varepsilon^2
\,\quad\,\,\textrm{in}\,\,\,\,\Omega_\varepsilon,\\
&\frac{\partial\Psi_0}{\partial\nu}=\varepsilon
\,\qquad\qquad\qquad\qquad\,
\textrm{on}\,\,\,\,\partial\Omega_\varepsilon.
\endaligned
\right.
\endaligned
$$
Choosing the positive constant $C_1$ larger if necessary, it is directly checked that
 $\psi$ meets  all the conditions
of the lemma for $R_1$ large but $\varepsilon$ small enough.
\end{proof}

\vspace{1mm}
\vspace{1mm}

Given $h\in L^\infty(\Omega_\varepsilon)$
and $\xi=(\xi_1,\ldots,\xi_m)\in\mathcal{O}_\varepsilon$, let us consider the linear equation
\begin{equation}\label{3.16}
\left\{\aligned
&\mathcal{L}(\phi)=-\Delta_{a(\varepsilon y)}\phi+\varepsilon^2\phi-W_{\xi'}\phi=h
\,\,\ \,
\,\textrm{in}\,\,\,\,\,\,\Omega_\varepsilon,\\
&\frac{\partial\phi}{\partial\nu}=0\,\,\,\,\,
\qquad\qquad\qquad\qquad\qquad\qquad
\ \,\textrm{on}\,\,\,\,\partial\Omega_{\varepsilon}.
\endaligned\right.
\end{equation}

\vspace{1mm}
\vspace{1mm}

\noindent{\bf Lemma 3.4.}\,\,{\it There exist $R_0>0$ and  $\varepsilon_0>0$ such that
for any $0<\varepsilon<\varepsilon_0$
and any solution $\phi$ of {\upshape (\ref{3.16})} with the orthogonality conditions
\begin{equation}\label{3.17}
\aligned
\int_{\Omega_\varepsilon}\chi_iZ_{ij}\phi=0\,\,\,\,\,
\,\,\,\,\forall\,\,i=1,\ldots,m,\,\,j=0,1,J_i,
\endaligned
\end{equation}
one has
\begin{equation*}\label{3.18}
\aligned
\|\phi\|_{L^{\infty}(\Omega_\varepsilon)}\leq C
\|h\|_{*},
\endaligned
\end{equation*}
where $C>0$ is independent of $\varepsilon$.
}

\vspace{1mm}
\vspace{1mm}

\begin{proof}
Set $R_0=4R_1$,
$R_1$ being the constant in Lemma 3.3.
By  (\ref{2.6}) it
follows that $\varepsilon\mu_i=o(1/|\log\varepsilon|^\kappa)$
for $\varepsilon$ small enough,
and by (\ref{2.4}), all $B_{R_1\mu_i}(\xi'_i)$ are disjointed
for any points  $\xi=(\xi_1,\ldots,\xi_m)\in\mathcal{O}_\varepsilon$.
Let $h$ be bounded and $\phi$ a solution to (\ref{3.16}) satisfying (\ref{3.17}).
We define the  inner norm of $\phi$ by
$$
\aligned
\|\phi\|_i=\sup_{y\in\overline{\Omega}_\varepsilon\cap\left(\bigcup_{i=1}^mB_{R_1\mu_i}(\xi'_i)\right)}
|\phi(y)|,
\endaligned
$$
and claim that there is a constant $C>0$ independent of  $\varepsilon$ such that
\begin{equation}\label{3.171}
\aligned
\|\phi\|_{L^{\infty}(\Omega_\varepsilon)}\leq C\left(\|\phi\|_i+
\|h\|_{*}\right).
\endaligned
\end{equation}
Indeed, set
$$
\aligned
\widetilde{\phi}(y)=C_1\left(\|\phi\|_i+
\|h\|_{*}
\right)\psi(y)
\,\qquad\forall\,\,\,y\in\overline{\Omega}_\varepsilon\setminus\bigcup_{i=1}^mB_{R_1\mu_i}(\xi'_i),
\endaligned
$$
where $\psi$ is the positive, uniformly bounded barrier constructed by Lemma 3.3 and
the constant $C_1>0$ is chosen larger if necessary,  independent of $\varepsilon$.
Then for
$y\in\Omega_\varepsilon\setminus\bigcup_{i=1}^mB_{R_1\mu_i}(\xi'_i)$,
$$
\aligned
\mathcal{L}(\widetilde{\phi}\pm\phi)(y)\geq C_{1}\,\|h\|_{*}\left\{
\sum_{i=1}^m\frac{\mu_i^\sigma}{|y-\xi'_i|^{2+\sigma}}+\varepsilon^2
\right\}\pm
h(y)\geq|h(y)|\pm h(y)\geq0,
\endaligned
$$
for
$y\in\partial\Omega_\varepsilon\setminus\bigcup_{i=1}^mB_{R_1\mu_i}(\xi'_i)$,
$$
\aligned
\frac{\partial}{\partial\nu}(\widetilde{\phi}\pm\phi)(y)\geq 0,
\endaligned
$$
and for
$y\in\Omega_\varepsilon\cap\left(\bigcup_{i=1}^m\partial B_{R_1\mu_i}(\xi'_i)\right)$,
$$
\aligned
(\widetilde{\phi}\pm\phi)(y)>\|\phi\|_{i}\pm\phi(y)\geq
|\phi(y)|\pm\phi(y)\geq 0.
\endaligned
$$
From  the maximum principle (see
\cite{PW}), it follows that
$-\widetilde{\phi}\leq\phi\leq\widetilde{\phi}$ on
$\overline{\Omega}_\varepsilon\setminus\bigcup_{i=1}^mB_{R_1\mu_i}(\xi'_i)$,
which gives estimate (\ref{3.171}).

We prove the lemma by contradiction. Assume that there exist a  sequence
 $\varepsilon_n\rightarrow+\infty$,
points $\xi^n=(\xi_1^n,\ldots,\xi_m^n)\in\mathcal{O}_{\varepsilon_n}$,
functions $h_n$,  and associated solutions $\phi_n$ of
equation (\ref{3.16}) with orthogonality conditions (\ref{3.17})
such that
\begin{equation}\label{3.181}
\aligned
\|\phi_n\|_{L^{\infty}(\Omega_{\varepsilon_n})}=1
\,\,\quad\,\,
\textrm{and}
\,\,\quad\,\,\|h_n\|_{*}\rightarrow0,
\,\,\quad\,\,\textrm{as}\,\,\,\,n\rightarrow+\infty.
\endaligned
\end{equation}
For each $k\in\{1,\ldots,l\}$, we have $\xi_k^n\in\Omega$ and we consider
$\widehat{\phi}^n_k(z)=\phi_n\big(\mu_k^nz+(\xi^n_k)'\big)$,
where $\mu^n=(\mu^n_1,\ldots,\mu_m^n)$ and $(\xi^n_k)'=\xi^n_k/\varepsilon_n$.
Note that
$$
\aligned
h_n(y)=
\big(-\Delta_{a(\varepsilon_n y)}\phi_n+\varepsilon_n^2\phi-W_{(\xi^n)'}\phi_n\big)\big|_{y=\mu_{k}^n z+(\xi^n_k)'}
=(\mu_k^n)^{-2}\left[
-\Delta_{\widehat{a}_n}\widehat{\phi}_k^n
+\varepsilon_n^2(\mu_k^n)^{2}\widehat{\phi}_k^n
-(\mu_k^n)^{2}\widehat{W}^n\widehat{\phi}_k^n
\right](z),
\endaligned
$$
where
$$
\aligned
\widehat{a}_n(z)=a(\varepsilon_n\mu_{k}^nz+\xi^n_k),
\qquad\qquad
\widehat{W}^n(z)=W_{(\xi^n)'}(\mu_{k}^nz+(\xi^n_k)').
\endaligned
$$
By the expansion of $W_{(\xi^n)'}$ in (\ref{3.2})
and elliptic regularity,
$\widehat{\phi}^n_k$
converges uniformly over
compact sets to a bounded solution $\widehat{\phi}^{\infty}_k$ of equation
\begin{equation*}\label{3.7}
\aligned
\Delta
\phi+\frac{8}{(1+|z|^2)^2}\phi=0
\,\quad\textrm{in}\,\,\,\mathbb{R}^2,
\endaligned
\end{equation*}
which satisfies
\begin{equation}\label{3.201}
\aligned
\int_{\mathbb{R}^2}\chi Z_j\widehat{\phi}_k^{\infty}=0
\quad\,\,\,\textrm{for}\,\,\,\,j=0,\,1,\,2.
\endaligned
\end{equation}
However, by the result of \cite{BP,CL},
$\widehat{\phi}^{\infty}_k$ must be
a linear combination of $Z_j$, $j=0,1,2$.
Notice that $\int_{\mathbb{R}^2}\chi Z_jZ_{t}=0$ for $j\neq t$
and $\int_{\mathbb{R}^2}\chi Z_j^2>0$.
Hence (\ref{3.201}) implies  $\widehat{\phi}_k^{\infty}\equiv0$.

As for each $k\in\{l+1,\ldots,m\}$, we have $\xi_k^n\in\partial\Omega$ and we consider
$\widehat{\phi}^n_k(z)=\phi_n\big((A^n_k)^{-1}\mu_k^nz+(\xi^n_k)'\big)$,
where $A_k^n: \mathbb{R}^2\rightarrow\mathbb{R}^2$ is a rotation map
such that $A_k^n\nu_{\Omega_{\varepsilon_n}}\big((\xi_k^n)'\big)=\nu_{\mathbb{R}_+^2}\big(0\big)$.
Similarly to the above argument, we have that
$\widehat{\phi}^n_k$
converges uniformly over compact sets
to a bounded solution $\widehat{\phi}^{\infty}_k$ of equation
\begin{equation*}\label{3.8}
\aligned
\Delta
\phi+\frac{8}{(1+|z|^2)^2}\phi=0
\,\quad\textrm{in}\,\,\,\mathbb{R}^2_{+},
\,\qquad\qquad\,\frac{\partial\phi}{\partial\nu}=0
\,\quad\textrm{on}\,\,\,\partial\mathbb{R}^2_{+},
\endaligned
\end{equation*}
which satisfies
\begin{equation}\label{3.221}
\aligned
\int_{\mathbb{R}_{+}^2}\chi Z_j\widehat{\phi}_k^{\infty}=0
\quad\,\,\,\textrm{for}\,\,\,\,j=0,\,1.
\endaligned
\end{equation}
Then $\widehat{\phi}^{\infty}_k$ is
a linear combination of $Z_j$, $j=0,1$.
Notice that $\int_{\mathbb{R}_{+}^2}\chi Z_jZ_t=0$ for $j\neq t$
and $\int_{\mathbb{R}_{+}^2}\chi Z_j^2>0$.
Hence (\ref{3.221}) implies  $\widehat{\phi}_k^{\infty}=0$
and then
$\lim_{n\rightarrow+\infty}\|\phi_n\|_i=0$.
But by (\ref{3.171})-(\ref{3.181}),
$\liminf_{n\rightarrow+\infty}\|\phi_n\|_i>0$,
which is a contradiction.
\end{proof}

\vspace{1mm}
\vspace{1mm}
\vspace{1mm}

\noindent{\bf Lemma 3.5.}\,\,{\it For $\varepsilon>0$  small enough, if
$\phi$ solves {\upshape (\ref{3.16})} and satisfies
\begin{equation}\label{3.18}
\aligned
\int_{\Omega_\varepsilon}\chi_iZ_{ij}\phi=0\,\,\,\,\,
\,\,\,\,\forall\,\,i=1,\ldots,m,\,\,j=1,J_i,
\endaligned
\end{equation}
then
\begin{equation}\label{3.19}
\aligned
\|\phi\|_{L^{\infty}(\Omega_\varepsilon)}\leq C|\log\varepsilon|\, \|h\|_{*},
\endaligned
\end{equation}
where $C>0$ is independent of $\varepsilon$.}

\vspace{1mm}

\begin{proof}
According to the results in Lemma 3.4 of \cite{DW} and Lemma 4.5 of \cite{D},
for simplicity  we  consider the validity of estimate (\ref{3.19}) only
when the $m$ concentration points $\xi=(\xi_1,\ldots,\xi_m)\in\mathcal{O}_\varepsilon$ satisfy
the relation $|\xi_i-\xi_k|\leq2d$ for any $i,k=1,\ldots,m$, $i\neq k$ and for
any $d>0$ sufficiently small, fixed and independent of $\varepsilon$.
Let $R>R_0+1$ be a large but fixed number. For any $i=1,\ldots,m$, we define
\begin{equation}\label{3.20}
\aligned
\widehat{Z}_{i0}(y)=Z_{i0}(y)-\frac1{\mu_i}
+a_{i0}G_a(\varepsilon y,\xi_i),
\endaligned
\end{equation}
where
\begin{equation}\label{3.21}
\aligned
a_{i0}=\frac1{\mu_i\big[H_a(\xi_i,\xi_i)-\frac{4}{c_i}\log(\varepsilon \mu_i R)\big]}.
\endaligned
\end{equation}
From estimate (\ref{2.6})  and definitions  (\ref{3.10}) and (\ref{3.13}) we have
\begin{equation}\label{3.22}
\aligned
C_1|\log\varepsilon|\leq-\log(\varepsilon \mu_i R)
\leq C_2|\log\varepsilon|,
\endaligned
\end{equation}
and
\begin{equation}\label{3.23}
\aligned
\widehat{Z}_{i0}(y)=O\left(
\frac{\,G_a(\varepsilon y,\xi_i)\,}{\mu_i|\log\varepsilon|}
\right).
\endaligned
\end{equation}
Let   $\eta_1$ and $\eta_2$ be  radial smooth cut-off functions in $\mathbb{R}^2$ such that
$$
\aligned
&0\leq\eta_1\leq1;\,\,\,\ \,\,\,|\nabla\eta_1|\leq C\,\,\ \textrm{in}\,\,\,\mathbb{R}^2;
\,\,\,\ \,\,\,\eta_1\equiv1\,\,\ \textrm{in}\,\,\,B_R(0);\,\,\,\,\ \,\,\,
\,\eta_1\equiv0\,\,\ \textrm{in}\,\,\,\mathbb{R}^2\setminus B_{R+1}(0);\\[1mm]
&0\leq\eta_2\leq1;\,\,\,\ \,\,\,|\nabla\eta_2|\leq C\,\,\ \textrm{in}\,\,\,\mathbb{R}^2;
\,\,\,\ \,\,\,\eta_2\equiv1\,\,\ \textrm{in}\,\,\,B_{3d}(0);\,\,\,
\,\,\,\,\,\,\eta_2\equiv0\,\,\ \textrm{in}\,\,\,\mathbb{R}^2\setminus B_{6d}(0).
\endaligned
$$
Denote that for any $i=1,\ldots,l$,
\begin{equation}\label{3.24}
\aligned
\eta_{i1}(y)=
\eta_1\left(\frac{1}{\mu_i}\big|y-\xi_i'\big|\right),
\,\,\qquad\quad\,\,
\eta_{i2}(y)=
\eta_2\left(\varepsilon\big|y-\xi'_i\big|\right),
\endaligned
\end{equation}
and for any $i=l+1,\ldots,m$,
\begin{equation}\label{3.25}
\aligned
\eta_{i1}(y)=
\eta_1\left(\frac{1}{\mu_i}\big|F_i^\varepsilon(y)\big|\right),
\,\,\qquad\quad\,\,
\eta_{i2}(y)=
\eta_2\left(\varepsilon\big|F_i^\varepsilon(y)\big|\right).
\endaligned
\end{equation}

Now define
\begin{equation}\label{3.26}
\aligned
\widetilde{Z}_{i0}(y)=\eta_{i1}Z_{i0}+(1-\eta_{i1})\eta_{i2}\widehat{Z}_{i0}.
\endaligned
\end{equation}
Given $\phi$ satisfying (\ref{3.16}) and (\ref{3.18}), let
\begin{equation}\label{3.27}
\aligned
\widetilde{\phi}=\phi+\sum\limits_{i=1}^{m}d_i\widetilde{Z}_{i0}+\sum_{i=1}^m\sum\limits_{j=1}^{J_i}e_{ij}\chi_iZ_{ij}.
\endaligned
\end{equation}
We can adjust
$d_i$ and $e_{ij}$  such that
$\widetilde{\phi}$ satisfies the orthogonality conditions
\begin{equation}\label{3.28}
\aligned
\int_{\Omega_\varepsilon}\chi_iZ_{ij}\widetilde{\phi}=0,
\,\,\quad\,\,\,\,i=1,\ldots,m,\,\,j=0,1,J_i.
\endaligned
\end{equation}
Indeed,
testing (\ref{3.27}) by $\chi_iZ_{ij}$, $i=1,\ldots,m$,
$j=0,1,J_i$  and using (\ref{3.18}), (\ref{3.28})  and the fact that
$\chi_i\chi_k\equiv0$ if $i\neq k$, we find
\begin{equation}\label{3.29}
\aligned
d_i\int_{\Omega_\varepsilon}\chi_iZ_{i0}\widetilde{Z}_{i0}
+\sum_{k\neq i}^md_k\int_{\Omega_\varepsilon}\chi_iZ_{i0}\widetilde{Z}_{k0}
+\sum\limits_{t=1}^{J_i}e_{it}\int_{\Omega_\varepsilon}\chi_i^2Z_{i0}Z_{it}
=-\int_{\Omega_\varepsilon}\chi_iZ_{i0}\phi,
\endaligned
\end{equation}
\begin{equation}\label{3.30}
\aligned
d_i\int_{\Omega_\varepsilon}\chi_iZ_{ij}\widetilde{Z}_{i0}
+\sum_{k\neq i}^md_k\int_{\Omega_\varepsilon}\chi_iZ_{ij}\widetilde{Z}_{k0}
+\sum_{t=1}^{J_i}e_{it}\int_{\Omega_\varepsilon}\chi^2_iZ_{ij}Z_{it}=0,
\,\quad\ \ \,\ \,j=1,\,J_i.
\endaligned
\end{equation}
Note that for any $i=1,\ldots,l,\,$   $j=1,2\,$ and $\,t=1,2$,
$$
\aligned
\int_{\Omega_\varepsilon}\chi_iZ_{i0}\widetilde{Z}_{i0}=
\int_{\mathbb{R}^2}\chi Z^2_{0}=C_0>0,
\qquad\qquad\qquad\qquad
\int_{\Omega_\varepsilon}\chi_i^2Z_{i0}Z_{it}=
\int_{\mathbb{R}^2}\chi^2 Z_{0}Z_t=0,
\endaligned
$$
$$
\aligned
\int_{\Omega_\varepsilon}\chi_iZ_{ij}\widetilde{Z}_{i0}=
\int_{\mathbb{R}^2}\chi Z_{j}Z_{0}=0,
\qquad\qquad\qquad\qquad
\int_{\Omega_\varepsilon}\chi^2_iZ_{ij}Z_{it}=
\int_{\mathbb{R}^2}\chi^2 Z_{j}Z_{t}
=C_j\delta_{jt},
\endaligned
$$
where $\delta_{jt}$ denotes the Kronecker's symbol, but for any $i=l+1,\ldots,m$ and $j=t=J_i=1$,
$$
\aligned
\int_{\Omega_\varepsilon}\chi_iZ_{i0}\widetilde{Z}_{i0}=
\int_{\mathbb{R}_{+}^2}\chi Z^2_{0}[1+O\big(\varepsilon\mu_i|z|\big)]=\frac{C_0}2+O\left(\varepsilon\mu_i\right),
\qquad
\int_{\Omega_\varepsilon}\chi_i^2Z_{i0}Z_{i1}=
\int_{\mathbb{R}_{+}^2}\chi^2 Z_{0}Z_1[1+O\big(\varepsilon\mu_i|z|\big)]=O\left(\varepsilon\mu_i\right),
\endaligned
$$
$$
\aligned
\int_{\Omega_\varepsilon}\chi_iZ_{i1}\widetilde{Z}_{i0}=
\int_{\mathbb{R}_{+}^2}\chi  Z_{1} Z_{0}  [1+O\big(\varepsilon\mu_i|z|\big) ] =O\left(\varepsilon\mu_i\right),
\,\,\qquad\,\,
\int_{\Omega_\varepsilon}\chi^2_iZ_{i1}^2=
\int_{\mathbb{R}_{+}^2}\chi^2  Z_{1}^2  [1+O\big(\varepsilon\mu_i|z|\big) ]
=\frac{C_1}{2} +O\left(\varepsilon\mu_i\right).
\endaligned
$$
From (\ref{3.23}) and (\ref{3.26})
it follows that  for any $i=1,\ldots,m$ and
$j=0,1,J_i$,
$$
\aligned
\int_{\Omega_\varepsilon}\chi_iZ_{ij}\widetilde{Z}_{k0}=
O\left(\frac{\mu_i\log|\log\varepsilon|}{\mu_k|\log\varepsilon|}
\right),\,
\quad\,\forall\,\,k\neq i.
\endaligned
$$
Hence by (\ref{3.30})  we can get that for any $i=1,\ldots,m$
and $j=1,J_i$,
\begin{equation*}\label{3.31}
\aligned
e_{ij}=\left(-d_i\int_{\Omega_\varepsilon}\chi_iZ_{ij}\widetilde{Z}_{i0}
-\sum_{k\neq i}^md_k\int_{\Omega_\varepsilon}\chi_iZ_{ij}\widetilde{Z}_{k0}
\right)
\left/
\int_{\Omega_\varepsilon}\chi^2_iZ^2_{ij},
\right.
\endaligned
\end{equation*}
and then
\begin{equation}\label{3.32}
\aligned
|e_{ij}|\leq
C\left(\varepsilon\mu_i|d_i|+
\sum_{k\neq i}^m
\frac{\mu_i\log|\log\varepsilon|}{\mu_k|\log\varepsilon|}|d_k|
\right).
\endaligned
\end{equation}
We need just to consider $d_i$. From (\ref{3.29}) it follows
that for any $i=1,\ldots,l$,
\begin{equation}\label{3.33}
\aligned
d_iC_0
+\sum_{k\neq i}^md_kO\left(\frac{\mu_i\log|\log\varepsilon|}{\mu_k|\log\varepsilon|}
\right)
=-\int_{\Omega_\varepsilon}\chi_iZ_{i0}\phi,
\endaligned
\end{equation}
and for any  $i=l+1,\ldots,m$,
\begin{equation}\label{3.34}
\aligned
\frac12d_iC_0\big[1+O\big(\varepsilon\mu_i\big)\big]
+\sum_{k\neq i}^md_kO\left(\frac{\mu_i\log|\log\varepsilon|}{\mu_k|\log\varepsilon|}
\right)+e_{i1}O\big(\varepsilon\mu_i\big)
=-\int_{\Omega_\varepsilon}\chi_iZ_{i0}\phi,
\endaligned
\end{equation}
where $e_{i1}$ satisfies (\ref{3.32}).
We denote $\mathcal{A}$ the coefficient matrix of equations (\ref{3.33})-(\ref{3.34})
with respect to $(d_1,\ldots,d_m)$. By the
above estimates,
$\mathcal{M}^{-1}\mathcal{A}\mathcal{M}$ is diagonally dominant, so invertible, where
$\mathcal{M}=\diag(\mu_1,\ldots,\mu_m)$. Hence $\mathcal{A}$ is invertible and
$(d_1,\ldots,d_m)$ is well defined.

Estimate (\ref{3.19}) is a  direct consequence of the following two claims.

\vspace{1mm}
\vspace{1mm}
\vspace{1mm}

\noindent{\bf Claim 1.}\,\,{\it
Let $\mathcal{L}=-\Delta_{a(\varepsilon y)}+\varepsilon^2-W_{\xi'}$, then for any $i=1,\ldots,m$
and $j=1, J_i$,
\begin{equation}\label{3.35}
\aligned
\big\|\mathcal{L}(\chi_iZ_{ij})\big\|_{*}\leq\frac{C}{\mu_i},
\,\quad\quad\quad\ \,\,\,\quad\,
\big\|\mathcal{L}(\widetilde{Z}_{i0})\big\|_{*}\leq
C\frac{\log^2 |\log\varepsilon|}{\mu_i|\log\varepsilon|}.
\endaligned
\end{equation}
}

\noindent{\bf Claim 2.}\,\,{\it
For any $i=1,\ldots,m$ and $j=1,J_i$,
\begin{equation*}\label{3.36}
\aligned
|d_i|\leq C\mu_i |\log\varepsilon|\,\|h\|_{*},
\,\quad\qquad\quad\,\,\quad
|e_{ij}|\leq C\mu_i\log(|\log\varepsilon|)\,\|h\|_{*}.
\endaligned
\end{equation*}
}

In fact, by  the definition of
$\widetilde{\phi}$ in (\ref{3.27}) we obtain
\begin{equation}\label{3.37}
\aligned\left\{\aligned
&
\mathcal{L}(\widetilde{\phi})=h+\sum\limits_{i=1}^{m}d_i\mathcal{L}(\widetilde{Z}_{i0})
+\sum_{i=1}^m\sum_{j=1}^{J_i}e_{ij}\mathcal{L}(\chi_iZ_{ij})
\,\quad\,\textrm{in}\,\,\,\ \,\,\Omega_\varepsilon,\\
&\frac{\partial\widetilde{\phi}}{\partial\nu}=0
\,\qquad\qquad\qquad\qquad\quad\qquad\,
\,\qquad\qquad\qquad\,\,
\,\textrm{on}\,\,\,\,\,\partial\Omega_\varepsilon.
\endaligned\right.
\endaligned
\end{equation}
Since (\ref{3.28}) holds, the previous lemma allows us to conclude
\begin{equation}\label{3.38}
\aligned
\|\widetilde{\phi}\|_{L^{\infty}(\Omega_\varepsilon)}\leq
C\left\{\|h\|_{*}
+\sum\limits_{i=1}^{m}|d_i|\big\|\mathcal{L}(\widetilde{Z}_{i0})\big\|_{*}
+\sum_{i=1}^m\sum_{j=1}^{J_i}|e_{ij}|\big\|\mathcal{L}(\chi_i Z_{ij})\big\|_{*}
\right\}\leq C\log^2(|\log\varepsilon|)\,\|h\|_{*}.
\endaligned
\end{equation}
Using the definition of $\widetilde{\phi}$
again and the fact that
\begin{equation}\label{3.39}
\aligned
\big\|\widetilde{Z}_{i0}\big\|_{L^{\infty}(\Omega_\varepsilon)}\leq\frac{C}{\mu_i}
\quad\quad\,\textrm{and}\,\quad\quad
\big\|\chi_iZ_{ij}\big\|_{L^{\infty}(\Omega_\varepsilon)}\leq\frac{C}{\mu_i},
\,\,\quad\,\forall\,\,\,i=1,\ldots,m,\,\,j=1,J_i,
\endaligned
\end{equation}
estimate (\ref{3.19}) then follows from estimate (\ref{3.38}) and Claim 2.

\vspace{1mm}
\vspace{1mm}
\vspace{1mm}
\vspace{1mm}

\noindent{\bf Proof of Claim 1.}
Let us first denote that  $z_i:=y-\xi'_i$  for any $i=1,\ldots,l$, but
$z_i:=F_i^\varepsilon(y)$  for any $i=l+1,\ldots,m$.
For any $i=l+1,\ldots,m$,
due to $F_i^\varepsilon(\xi'_i)=(0,0)$ and $\nabla F_i^\varepsilon(\xi'_i)=A_i$, we find
\begin{equation}\label{3.43}
\aligned
z_i=F_i^\varepsilon(y)=\frac{1}{\varepsilon}F_i\big(A_i(\varepsilon y-\xi_i)\big)
=A_i(y-\xi_i')\big\{1+O\big(\varepsilon A_i(y-\xi_i')\big)\big\},
\endaligned
\end{equation}
and
\begin{equation}\label{3.42}
\aligned
\nabla_y=A_i\nabla_{z_i}+O(\varepsilon|z_i|)
\nabla_{z_i},
\,\quad\,\,\,\quad\,
\quad
\,\quad\,\,\,\quad\,
-\Delta_y=-\Delta_{z_i}+O(\varepsilon|z_i|)
\nabla_{z_i}^2+O(\varepsilon)\nabla_{z_i}.
\endaligned
\end{equation}
Then for any $i=1,\ldots,m$
and $j=1, J_i$, by
 (\ref{3.2}), (\ref{3.10}) and (\ref{3.13})
we have that in  the region
$|z_i|\leq\mu_i(R_0+1)$,
$$
\aligned
\mathcal{L}(Z_{ij})=&
\big(-\Delta_y-W_{\xi'}\big)\left[\frac1{\mu_i}Z_{j}\left(\frac{z_i}{\mu_i}\right)
\right]
-\varepsilon\nabla_y\log a(\varepsilon y)\nabla_y\left[\frac1{\mu_i}Z_{j}\left(\frac{z_i}{\mu_i}\right)
\right]+\frac{\varepsilon^2}{\mu_i}Z_{j}\left(\frac{y-\xi'_i}{\mu_i}\right)\\[1mm]
=&
O\left(\frac1{\mu_i|\log\varepsilon|}\cdot\frac{8\mu_i^2}{(\mu_i^2+|y-\xi'_i|^2)^2}\right)
+O\left(
\frac{\varepsilon}{\mu_{i}^2+|y-\xi'_i|^2}
\right)+O\left(
\frac{\varepsilon^2}{(\mu_i^2+|y-\xi'_i|^2)^{1/2}}
\right).
\endaligned
$$
Hence
\begin{eqnarray*}\label{3.40}
\mathcal{L}(\chi_i Z_{ij})=\,
\chi_i\mathcal{L}(Z_{ij})
-2\nabla\chi_i\nabla Z_{ij}
-Z_{ij}\big[\Delta\chi_i
+\varepsilon\nabla\log a(\varepsilon y)\nabla\chi_i\big]
\qquad\qquad\qquad\qquad\qquad\qquad\qquad\qquad
\,\,\,\,\,
\nonumber\\[1mm]
=O\left(\frac1{\mu_i|\log\varepsilon|}\cdot\frac{8\mu_i^2}{(\mu_i^2+|y-\xi'_i|^2)^2}\right)
+O\left(\frac1{\mu_i}\cdot
\frac{1}{\mu_i^2+|y-\xi'_i|^2}
\right)
+O\left(\frac1{\mu_i^2}\cdot
\frac{1}{(\mu_{i}^2+|y-\xi'_i|^2)^{1/2}}
\right),
\end{eqnarray*}
which, together with
the definition of  $\|\cdot\|_*$ in (\ref{2.32}),  implies
$\big\|\mathcal{L}(\chi_iZ_{ij})\big\|_{*}=O\left(1/\mu_i\right)$
for all $i=1,\ldots,m$ and $j=1,J_i$.

Let us prove the second inequality in (\ref{3.35}).
Consider  four regions
$$
\aligned
\Omega_{1}=\left\{y\in\Omega_\varepsilon\big|\,|z_i|\leq\mu_iR\right\},
\,\quad\,\quad\quad\quad\quad\quad\quad\quad\,
\Omega_{2}=\left\{y\in\Omega_\varepsilon\big|\,\mu_iR<|z_i|\leq\mu_i(R+1)\right\},\\[1mm]
\Omega_{3}=\left\{y\in\Omega_\varepsilon\left|\,\mu_i(R+1)<|z_i|\leq3d/\varepsilon\right.\right\},
\,\quad\quad\quad\quad\,\,\,\,\,\,\,
\Omega_{4}=\left\{y\in\Omega_\varepsilon\left|\,3d/\varepsilon<|z_i|\leq
6d/\varepsilon\right.\right\}.
\quad\quad
\endaligned
$$
Notice  first that
\begin{equation}\label{3.46}
\aligned
\left|
Z_{i0}-\frac{1}{\mu_i}
\right|=\frac{2\mu_i}{\mu_i^2+|z_i|^2}
=O\left(\frac{\mu_i}{(\mu_i+|y-\xi_i'|)^2}\right),
\endaligned
\end{equation}
and for any $\mu_iR<|z_i|\leq6d/\varepsilon$,
\begin{equation}\label{3.47}
\aligned
Z_{i0}-\widehat{Z}_{i0}=\frac1{\mu_i}
-a_{i0}G_a(\varepsilon y,\xi_i)=
\frac1{\mu_i\big[H_a(\xi_i,\xi_i)-\frac{4}{c_i}\log(\varepsilon \mu_i R)\big]}\left[
\frac{4}{c_i}\log\frac{|y-\xi_i'|}{\mu_iR}+O\big(
\varepsilon^\alpha|y-\xi_i'|^\alpha\big)
\right],
\endaligned
\end{equation}
and for any $|z_i|\leq\mu_i(R+1)$, by (\ref{3.2}),  (\ref{3.10}) and  (\ref{3.13}),
\begin{equation}\label{3.48}
\aligned
\left[\frac{8\mu_i^2}{(\mu_i^2+|y-\xi'_i|^2)^2}-W_{\xi'}\right]Z_{i0}
=O\left(\frac{\log^2\mu_i}{\mu_i|\log\varepsilon|}\cdot\frac{8\mu_i^2}{(\mu_i^2+|y-\xi'_i|^2)^2}\right).
\endaligned
\end{equation}
In $\Omega_1$,
$$
\aligned
\mathcal{L}(\widetilde{Z}_{i0})=\mathcal{L}(Z_{i0})
=\left[-\Delta_{a(\varepsilon y)}Z_{i0}-\frac{8\mu_i^2}{(\mu_i^2+|y-\xi'_i|^2)^2}Z_{i0}
\right]+\left[\frac{8\mu_i^2}{(\mu_i^2+|y-\xi'_i|^2)^2}-W_{\xi'}\right]Z_{i0}
+\varepsilon^2Z_{i0}.
\endaligned
$$
By (\ref{3.14}) and  (\ref{3.48}),
\begin{equation}\label{3.49}
\aligned
\big|\mathcal{L}(\widetilde{Z}_{i0})(y)\big|=
O\left(\frac{\log^2\mu_i}{\mu_i^3|\log\varepsilon|}\right),
\qquad\forall\,\,y\in\Omega_1.
\endaligned
\end{equation}
In $\Omega_2$,
\begin{eqnarray*}
\mathcal{L}(\widetilde{Z}_{i0})
=\left[-\Delta_{a(\varepsilon y)}Z_{i0}-\frac{8\mu_i^2}{(\mu_i^2+|y-\xi'_i|^2)^2}Z_{i0}
\right]+\left[\frac{8\mu_i^2}{(\mu_i^2+|y-\xi'_i|^2)^2}-W_{\xi'}\right]Z_{i0}
+\varepsilon^2\left(
Z_{i0}-\frac1{\mu_i}
\right)
&&\nonumber\\
+W_{\xi'}(1-\eta_{i1})(Z_{i0}-\widehat{Z}_{i0})
+\frac{\varepsilon^2}{\mu_i}\eta_{i1}
-2\nabla\eta_{i1}\nabla(Z_{i0}-\widehat{Z}_{i0})-(Z_{i0}-\widehat{Z}_{i0})\Delta_{a(\varepsilon y)}\eta_{i1}.
\qquad\quad
&&
\end{eqnarray*}
Using (\ref{3.22}) and (\ref{3.47}) we conclude  that for any $\mu_iR<|z_i|\leq\mu_i(R+1)$,
\begin{equation}\label{3.50}
\aligned
|Z_{i0}-\widehat{Z}_{i0}|=O\left(
\frac{1}{\mu_iR|\log\varepsilon|}
\right)
\,\,\quad\quad\,\textrm{and}\,\,\quad\quad\,
|\nabla\big(Z_{i0}-\widehat{Z}_{i0}\big)|=O\left(
\frac{1}{\mu_i^2R|\log\varepsilon|}
\right).
\endaligned
\end{equation}
Moreover, $|\nabla\eta_{i1}|=O\left(1/\mu_i\right)$ and $|\Delta_{a(\varepsilon y)}\eta_{i1}|=O\left(1/\mu_i^{2}\right)$.
From  (\ref{3.14}), (\ref{3.46}) and (\ref{3.48}) we can derive that
\begin{equation}\label{3.51}
\aligned
\big|\mathcal{L}(\widetilde{Z}_{i0})(y)\big|
=O\left(
\frac{1}{\mu_i^3R|\log\varepsilon|}
\right),
\qquad\forall\,\,y\in\Omega_2.
\endaligned
\end{equation}
In $\Omega_3$, by (\ref{3.14}),  (\ref{3.46}) and (\ref{3.47}),
$$
\aligned
\mathcal{L}(\widetilde{Z}_{i0})=&\mathcal{L}(\widehat{Z}_{i0})
=\mathcal{L}(Z_{i0})-\mathcal{L}(Z_{i0}-\widehat{Z}_{i0})
\\[1.6mm]
=&\left[\frac{8\mu_i^2}{(\mu_i^2+|y-\xi'_i|^2)^2}-W_{\xi'}\right]Z_{i0}
+W_{\xi'}\left[\frac1{\mu_i}
-a_{i0}G_a(\varepsilon y, \xi_i)\right]
+O\left(\frac{\varepsilon\mu_i}{
(\mu_i+|y-\xi'_i|)^3}\right)
+O\left(\frac{\varepsilon^2\mu_i}{(\mu_i+|y-\xi_i'|)^2}\right)
\\[1.3mm]
\equiv&\mathcal{A}_1+\mathcal{A}_2+O\left(\frac{\varepsilon\mu_i}{
(\mu_i+|y-\xi'_i|)^3}\right)
+O\left(\frac{\varepsilon^2\mu_i}{(\mu_i+|y-\xi_i'|)^2}\right).
\endaligned
$$
For the estimation of the first two terms, we split  $\Omega_3$ into some subregions:
$$
\aligned
&\Omega_{3,i}=\left\{y\in\Omega_3\left|\,\,\mu_i(R+1)<|z_i|\leq
\frac14\mu_i|\log\varepsilon|^\tau\,\right.\right\},\\
\Omega_{3,k}=&\left\{y\in\Omega_3\left|\,\,|z_k|\leq\frac14\mu_k|\log\varepsilon|^\tau\,\right.\right\},
\,\,\,\,k\neq i,\quad\,\,\textrm{and}\quad\,\,
\widetilde{\Omega}_3=\Omega_3\setminus\bigcup_{t=1}^m\Omega_{3,t}.
\endaligned
$$
In $\Omega_{3,i}$,
by (\ref{2.13})  and   (\ref{3.2}) we find
$$
\aligned
\mathcal{A}_1=&
\frac{8\mu^2_i}{(\mu_{i}^2+|y-\xi'_i|^2)^2}
\left\{\frac{p-1}{p}\frac1{\gamma^p}
\left[1+\omega^{1}_{\mu_i}+\frac{1}{2}(\omega_{\mu_i})^2+2\omega_{\mu_i}\right]\left(y-\xi'_i\right)
+
O\left(\frac{\log^4(\mu_i+|y-\xi'_i|)}{\gamma^{2p}}\right)
\right\}O\left(\frac1{\mu_i}\right)\\
=&\frac{8\mu^2_i}{(\mu_{i}^2+|y-\xi'_i|^2)^2}
O\left(\frac{\log^2(\mu_i+|y-\xi'_i|)}{\mu_i|\log\varepsilon|}\right),
\endaligned
$$
and by (\ref{3.22}) and (\ref{3.47}),
$$
\aligned
\mathcal{A}_2=
\frac{8\mu^2_i}{(\mu_i^2+|y-\xi'_i|^2)^2} O\left(\frac{\log|y-\xi_i'|-\log\mu_iR+\varepsilon^\alpha|y-\xi_i'|^\alpha}
{\mu_i|\log\varepsilon|}\right),
\endaligned
$$
which implies
\begin{equation}\label{3.52}
\aligned
\big|\mathcal{L}(\widetilde{Z}_{i0})(y)\big|=\frac{8\mu^2_i}{(\mu_{i}^2+|y-\xi'_i|^2)^2}
O\left(\frac{\log^2(\mu_i+|y-\xi'_i|)}{\mu_i|\log\varepsilon|}
+\frac{\log|y-\xi_i'|-\log\mu_iR}
{\mu_i|\log\varepsilon|}
\right),
\quad\forall\,\,y\in\Omega_{3,i}.
\endaligned
\end{equation}
In $\widetilde{\Omega}_{3}$, by (\ref{3.1}),
\begin{equation}\label{3.52a}
\aligned
\big|\mathcal{L}(\widetilde{Z}_{i0})(y)\big|=
\sum_{k=1}^m\left[\frac{8\mu^2_k}{(\mu_{k}^2+|y-\xi'_k|^2)^2}\right]^{1+
O\left(\frac{\log|\log\varepsilon|}{|\log\varepsilon|}\right)}
O\left(\frac1{\mu_i}+
\frac{\,\log|y-\xi_i'|-\log\mu_iR\,}
{\mu_i|\log\varepsilon|}
\right),
\quad\forall\,\,y\in\widetilde{\Omega}_{3}.
\endaligned
\end{equation}
As in $\Omega_{3,k}$ with $k\neq i$, by (\ref{3.2}), (\ref{3.14}), (\ref{3.23})
and (\ref{3.46}),
\begin{eqnarray}\label{3.53}
\mathcal{L}(\widetilde{Z}_{i0})=
\frac{8\mu_i^2}{(\mu_{i}^2+|y-\xi'_i|^2)^2}Z_{i0}
-\left[
\Delta_{a(\varepsilon y)} Z_{i0}+\frac{8\mu_i^2}{(\mu_{i}^2+|y-\xi'_i|^2)^2}Z_{i0}
\right]
+\varepsilon^2\left(
Z_{i0}-\frac1{\mu_i}
\right)-W_{\xi'}\widehat{Z}_{i0}
&&\nonumber\\[1.5mm]
=O\left(\frac{8\mu^2_k}{(\mu^2_k+|y-\xi'_k|^2)^2}
\cdot\frac{\log|\log\varepsilon|}{\mu_i|\log\varepsilon|}
\right).
\qquad\qquad\qquad\qquad\qquad\qquad
\qquad\quad\qquad\qquad\qquad\,\,\,
&&
\end{eqnarray}
Finally, in $\Omega_4$,
$$
\aligned
\mathcal{L}(\widetilde{Z}_{i0})=&\frac{8\mu_i^2}{(\mu_i^2+|y-\xi'_i|^2)^2}\eta_{i2}Z_{i0}
-\eta_{i2}\left[\Delta_{a(\varepsilon y)} Z_{i0}
+\frac{8\mu_i^2}{(\mu_i^2+|y-\xi'_i|^2)^2}Z_{i0}
\right]+\varepsilon^2\eta_{i2}\left(
Z_{i0}-\frac1{\mu_i}
\right)
\\[1.5mm]
&-\eta_{i2}W_{\xi'}\widehat{Z}_{i0}
-2\nabla\eta_{i2}\nabla \widehat{Z}_{i0}-\widehat{Z}_{i0}\Delta_{a(\varepsilon y)}\eta_{i2}.
\endaligned
$$
Note that $W_{\xi'}=O(\varepsilon^{\frac{4}{p}-\sigma})$ in $\Omega_4$.
Moreover,
$|\nabla\eta_{i2}|=O\left(\varepsilon/d\right)$,
$|\Delta_{a(\varepsilon y)}\eta_{i2}|=O\left(\varepsilon^2/d^2\right)$,
\begin{equation}\label{3.54}
\aligned
|\widehat{Z}_{i0}|=O\left(
\frac{|\log d|}{\mu_i|\log\varepsilon|}
\right)\,\,\qquad\quad\,\textrm{and}\,\,\quad\qquad\,
|\nabla\widehat{Z}_{i0}|=O\left(
\frac{\varepsilon}{d\mu_i|\log\varepsilon|}
\right).
\endaligned
\end{equation}
By (\ref{3.14}) and (\ref{3.46}),
\begin{equation}\label{3.55}
\aligned
\big|\mathcal{L}(\widetilde{Z}_{i0})(y)\big|
=O\left(
\frac{\varepsilon^2|\log d|}{\mu_id^2|\log\varepsilon|}
\right),
\qquad\forall\,\,y\in\Omega_4.
\endaligned
\end{equation}
Combining   (\ref{3.49}), (\ref{3.51}), (\ref{3.52}), (\ref{3.52a}), (\ref{3.53}) and (\ref{3.55}), we arrive at
$$
\aligned
\big\|\mathcal{L}(\widetilde{Z}_{i0})\big\|_{*}\leq
C\left(
\frac{\log^2\mu_i}{\mu_i|\log\varepsilon|}
+\frac{\log |\log\varepsilon|}{\mu_i|\log\varepsilon|}
\right)=O\left(
\frac{\log^2 |\log\varepsilon|}{\mu_i|\log\varepsilon|}\right),
\,\,\quad\,\forall\,\,i=1,\ldots,m.
\endaligned
$$

\vspace{1mm}
\vspace{1mm}

\noindent{\bf Proof of Claim 2.}
Testing equation (\ref{3.37}) against
$a(\varepsilon y)\widetilde{Z}_{i0}$  and using estimates (\ref{3.38})-(\ref{3.39}), we
find
$$
\aligned
\sum_{k=1}^md_k&
\int_{\Omega_\varepsilon}a(\varepsilon y)\widetilde{Z}_{k0}\mathcal{L}(\widetilde{Z}_{i0})
\\
=&-\int_{\Omega_\varepsilon} a(\varepsilon y)h\widetilde{Z}_{i0}
+\int_{\Omega_\varepsilon}a(\varepsilon y)\widetilde{\phi}\mathcal{L}(\widetilde{Z}_{i0})-\sum_{k=1}^m\sum_{t=1}^{J_k}e_{kt}
\int_{\Omega_\varepsilon}a(\varepsilon y)\chi_kZ_{kt}\mathcal{L}(\widetilde{Z}_{i0})
\\[1mm]
\leq&C\frac{\|h\|_{*}}{\mu_i}
+C\big\|\mathcal{L}(\widetilde{Z}_{i0})\big\|_{*}
\left(\|\widetilde{\phi}\|_{L^{\infty}(\Omega_\varepsilon)}
+\sum_{k=1}^m\sum_{t=1}^{J_k}\frac{1}{\mu_k}|e_{kt}|\right)
\\[1mm]
\leq&C\frac{\|h\|_{*}}{\mu_i}
+
C\big\|\mathcal{L}(\widetilde{Z}_{i0})\big\|_{*}
\left[\|h\|_{*}+
\sum\limits_{k=1}^{m}|d_k|\big\|\mathcal{L}(\widetilde{Z}_{k0})\big\|_{*}
+\sum_{k=1}^m\sum_{t=1}^{J_k}|e_{kt}|\left(\frac{1}{\mu_k}
+\big\|\mathcal{L}(\chi_kZ_{kt})\big\|_{*}
\right)
\right],
\endaligned
$$
where we have used that
$$
\aligned
\int_{\Omega_\varepsilon}\frac{\mu_i^\sigma}{(|y-\xi'_i|+\mu_i)^{2+\sigma}}\leq C,
\,\,\,\quad\,\,\,\forall\,\,\,i=1,\ldots,m.
\endaligned
$$
From  estimates (\ref{3.32}) and  (\ref{3.35})  it follows   that for any $i=1,\ldots,m$,
\begin{equation}\label{3.56}
\aligned
|d_i|\left|
\int_{\Omega_\varepsilon}a(\varepsilon y)\widetilde{Z}_{i0}\mathcal{L}(\widetilde{Z}_{i0})
\right|
\leq
C\frac{\|h\|_{*}}{\mu_i}
+C\sum_{k=1}^m\frac{|d_k|\log^4|\log\varepsilon|}{\mu_i\mu_k|\log\varepsilon|^2}
+\sum_{k\neq i}^m\left|d_k
\int_{\Omega_\varepsilon}a(\varepsilon y)\widetilde{Z}_{k0}\mathcal{L}(\widetilde{Z}_{i0})
\right|.
\endaligned
\end{equation}
Observe that
\begin{eqnarray*}\label{3.45}
\mathcal{L}(\widetilde{Z}_{i0})=\eta_{i1}\mathcal{L}(Z_{i0}-\widehat{Z}_{i0})
+\eta_{i2}\mathcal{L}(\widehat{Z}_{i0})
-(Z_{i0}-\widehat{Z}_{i0})\Delta_{a(\varepsilon y)}\eta_{i1}
-2\nabla\eta_{i1}\nabla(Z_{i0}-\widehat{Z}_{i0})
-2\nabla\eta_{i2}\nabla\widehat{Z}_{i0}-\widehat{Z}_{i0}\Delta_{a(\varepsilon y)}\eta_{i2}.
\end{eqnarray*}
Then by (\ref{3.20}) and (\ref{3.26}),
\begin{equation*}\label{3.57}
\aligned
\int_{\Omega_\varepsilon}a(\varepsilon y)\widetilde{Z}_{i0}\mathcal{L}(\widetilde{Z}_{i0})
:=K+I,
\endaligned
\end{equation*}
where
$$
\aligned
K=\int_{\Omega_\varepsilon}a(\varepsilon y)\widetilde{Z}_{i0}\left[
-(Z_{i0}-\widehat{Z}_{i0})\Delta_{a(\varepsilon y)}\eta_{i1}-2\nabla\eta_{i1}\nabla(Z_{i0}-\widehat{Z}_{i0})
-2\nabla\eta_{i2}\nabla\widehat{Z}_{i0}-\widehat{Z}_{i0}\Delta_{a(\varepsilon y)}\eta_{i2}
\right],
\endaligned
$$
and
$$
\aligned
I=&\int_{\Omega_\varepsilon}a(\varepsilon y)\widetilde{Z}_{i0}\left[
\eta_{i1}\mathcal{L}(Z_{i0}-\widehat{Z}_{i0})
+\eta_{i2}\mathcal{L}(\widehat{Z}_{i0})
\right]\\[0.5mm]
=&\int_{\Omega_\varepsilon}a(\varepsilon y)\eta_{i2}^2\left\{Z_{i0}-\big(1-\eta_{i1}\big)
\left[\frac1{\mu_i}-a_{i0}G_a(\varepsilon y,\xi_i)\right]\right\}\times\left\{
W_{\xi'}\big(1-\eta_{i1}\big)\left[\frac1{\mu_i}-a_{i0}G_a(\varepsilon y,\xi_i)\right]
\right.\\[1mm]
&\left.
-\left[\Delta_{a(\varepsilon y)}Z_{i0}+\frac{8\mu_i^2}{(\mu_i^2+|y-\xi'_i|^2)^2}Z_{i0}
\right]+\left[\frac{8\mu_i^2}{(\mu_i^2+|y-\xi'_i|^2)^2}-W_{\xi'}\right]Z_{i0}
+\varepsilon^2\left(
Z_{i0}-\frac{1}{\mu_i}\right)
+\frac{\varepsilon^2}{\mu_i}\eta_{i1}
\right\}.
\endaligned
$$
Let us first estimate the expression $K$.
Integrating by parts the first term and the last term of $K$ respectively, we find
$$
\aligned
K=&-\int_{\Omega_2}a(\varepsilon y)\widehat{Z}_{i0}\nabla\eta_{i1}\nabla(Z_{i0}-\widehat{Z}_{i0})
+\int_{\Omega_2}a(\varepsilon y)(Z_{i0}-\widehat{Z}_{i0})^2|\nabla\eta_{i1}|^2
\\
&+\int_{\Omega_2}a(\varepsilon y)(Z_{i0}-\widehat{Z}_{i0})\nabla\eta_{i1}\nabla\widehat{Z}_{i0}
+\int_{\Omega_4}a(\varepsilon y)
|\widehat{Z}_{i0}|^2|\nabla\eta_{i2}|^2
\\[1.2mm]
=&K_{21}+K_{22}+K_{23}+K_4.
\endaligned
$$
From (\ref{3.6}), (\ref{3.10}), (\ref{3.13}),  (\ref{3.42}) and (\ref{3.50})  we have that
$|\nabla\widehat{Z}_{i0}|=O\big(\frac1{\mu_i^2 R^3}\big)$
and
$|\nabla\eta_{i1}|=O\big(\frac1{\mu_i}\big)$ in $\Omega_2$. Then
$$
\aligned
K_{22}=O\left(\frac{1}{\mu_i^2R|\log\varepsilon|^2}\right)
\,\,\quad\qquad\quad\,\,
\textrm{and}
\,\,\quad\qquad\quad\,\,
K_{23}=O\left(\frac{1}{\mu^2_iR^3|\log\varepsilon|}\right).
\endaligned
$$
By (\ref{3.54}),
$$
\aligned
K_4=O\left(\frac{|\log d|^2}{\mu_i^2|\log\varepsilon|^2}\right).
\endaligned
$$
Since $\widehat{Z}_{i0}=Z_{i0}\big[1+O\big(\frac{1}{R|\log\varepsilon|}\big)\big]$
in $\Omega_2$, by (\ref{2.6}),  (\ref{3.6}), (\ref{3.10}), (\ref{3.13}), (\ref{3.20}), (\ref{3.21}), (\ref{3.43}) and (\ref{3.42}) we can derive that
\begin{eqnarray}\label{3.581}
\nonumber
K=-\frac{a_{i0}}{\mu_i^2}
\int_{\{\mu_iR<|z_i|\leq\mu_i(R+1)\}}
\frac{1}{|y-\xi'_i|}a(\varepsilon y)Z_0\left(\frac{z_i}{\mu_i}\right)
\eta_1'\left(\frac{|z_i|}{\mu_i}\right)
\left[\frac{4}{c_i}+o(1)\right]dy
+O\left(\frac{1}{\mu^2_iR^3|\log\varepsilon|}\right)
\\[1.2mm]
\nonumber
=-\frac{c_i a_{i0}}{4\mu_i}\int_{R}^{R+1}
a(\xi_i)\eta_1'(r)\left[
\frac{4}{c_i}+O\left(\frac1{r^2}\right)
\right]dr+O\left(\frac{1}{\mu^2_iR^3|\log\varepsilon|}\right)
\qquad\qquad\qquad\qquad\qquad\qquad\,\,\,\,
\\[1mm]
=\frac14\frac{c_i a(\xi_i)}{\mu_i^2|\log\varepsilon|}\left[1+O\left(\frac1{R^2}\right)\right].
\qquad\qquad\qquad\qquad\qquad\qquad\qquad\,
\qquad\qquad\qquad\qquad\qquad\qquad\qquad
\end{eqnarray}
Next, we analyze the expression $I$.
From (\ref{2.6}), (\ref{3.1}),
(\ref{3.2}), (\ref{3.6}), (\ref{3.10}),   (\ref{3.13}), (\ref{3.14}),  (\ref{3.46}) and  (\ref{3.47})
we can estimate
$$
\aligned
\int_{|z_i|\leq\frac14\mu_i|\log\varepsilon|^\tau}
a(\varepsilon y)\eta_{i2}^2\left\{Z_{i0}-(1-\eta_{i1})
\left[\frac1{\mu_i}-a_{i0}G_a(\varepsilon y,\xi_i)\right]\right\}\times
W_{\xi'}\big(1-\eta_{i1}\big)\left[\frac1{\mu_i}-a_{i0}G_a(\varepsilon y,\xi_i)\right]
dy=O\left(\frac{1}{\mu_i^2R|\log\varepsilon|}\right),
\endaligned
$$
and
$$
\aligned
\int_{|z_i|\leq\frac14\mu_i|\log\varepsilon|^\tau}
\,a(\varepsilon y)\eta_{i2}^2\left\{Z_{i0}-(1-\eta_{i1})
\left[\frac1{\mu_i}-a_{i0}G_a(\varepsilon y,\xi_i)\right]\right\}\times
\left[\Delta_{a(\varepsilon y)}Z_{i0}+\frac{8\mu_i^2}{(\mu_i^2+|y-\xi'_i|^2)^2}Z_{i0}
\right]dy=O\left(\frac{\varepsilon}{\mu_i}\right),
\endaligned
$$
and
$$
\aligned
\int_{|z_i|\leq\frac14\mu_i|\log\varepsilon|^\tau}
\,a(\varepsilon y)\eta_{i2}^2\left\{Z_{i0}-(1-\eta_{i1})
\left[\frac1{\mu_i}-a_{i0}G_a(\varepsilon y,\xi_i)\right]\right\}\times
\left[\varepsilon^2\left(
Z_{i0}-\frac{1}{\mu_i}\right)
+\frac{\varepsilon^2}{\mu_i}\eta_{i1}\right]dy
=O\left(\varepsilon^2|\log\varepsilon|\right),
\endaligned
$$
and
$$
\aligned
&\int_{|z_i|\leq\frac14\mu_i|\log\varepsilon|^\tau}
\,a(\varepsilon y)\eta_{i2}^2(1-\eta_{i1})
\left[\frac1{\mu_i}-a_{i0}G_a(\varepsilon y,\xi_i)\right]
\left[\frac{8\mu_i^2}{(\mu_i^2+|y-\xi'_i|^2)^2}-W_{\xi'}\right]Z_{i0}dy\\[1mm]
=&\int_{\mu_iR<|z_i|\leq\frac14\mu_i|\log\varepsilon|^\tau}
\,\frac{8\mu^2_ia(\varepsilon y)}{(\mu_{i}^2+|y-\xi'_i|^2)^2}
\left\{\frac{p-1}{p}\frac1{\gamma^p}
\left[1+\omega^{1}_{\mu_i}+\frac{1}{2}(\omega_{\mu_i})^2+2\omega_{\mu_i}\right]\left(y-\xi'_i\right)
+O\left(\frac{\log^4(\mu_i+|y-\xi'_i|)}{\gamma^{2p}}\right)\right\}\\[1.5mm]
&\,\,\,\,
\times
O\left(\frac{\log|y-\xi_i'|-\log(\mu_iR)+
\varepsilon^\alpha|y-\xi_i'|^\alpha}{\mu_i^2|\log\varepsilon|}\right)dy
\\[0.5mm]
=&\,O\left(
\frac{\log^2|\log\varepsilon|}{\mu_i^2R|\log\varepsilon|^2}
\right).
\endaligned
$$
But  by   (\ref{3.23}),  (\ref{3.52a})   and (\ref{3.53}),
$$
\aligned
&\int_{|z_i|>\frac14\mu_i|\log\varepsilon|^\tau}
\,a(\varepsilon y)\widetilde{Z}_{i0}\left[
\eta_{i1}\mathcal{L}(Z_{i0}-\widehat{Z}_{i0})
+\eta_{i2}\mathcal{L}(\widehat{Z}_{i0})
\right]dy\\[0.5mm]
=&\sum_{k\neq i}^m\int_{\Omega_{3,k}}
a(\varepsilon y)\eta_{i2}^2\widehat{Z}_{i0}\mathcal{L}(\widehat{Z}_{i0})dy
+
\int_{\widetilde{\Omega}_{3}\cup\Omega_4}
a(\varepsilon y)\eta_{i2}^2\widehat{Z}_{i0}\mathcal{L}(\widehat{Z}_{i0})dy
\\
=&\sum_{k\neq i}^mO\left(\int_{0}^{\frac14\mu_k|\log\varepsilon|^\tau}
\frac{8\mu^2_k}{(\mu_{k}^2+r^2)^2}
\frac{\log^2 (|\log\varepsilon|)}{\mu_i^2|\log\varepsilon|^2}rdr
\right)+O
\left(\int_{\frac14\mu_i|\log\varepsilon|^\tau}^{+\infty}
\left[\frac{8\mu_i^2}{(\mu_i^2+r^2)^2}\right]^{1+
O\left(\frac{\log|\log\varepsilon|}{|\log\varepsilon|}\right)}
\frac{\log\varepsilon r}{\mu_i^2|\log\varepsilon|}rdr
\right)\\[0.5mm]
=&\,O\left(\frac{\log^2(|\log\varepsilon|)}{\mu^2_i|\log\varepsilon|^2}\right).
\endaligned
$$
So
$$
\aligned
I=&\int_{|z_i|\leq\frac14\mu_i|\log\varepsilon|^\tau}
\,a(\varepsilon y)\eta_{i2}^2Z_{i0}^2
\left[\frac{8\mu_i^2}{(\mu_i^2+|y-\xi'_i|^2)^2}-W_{\xi'}\right]dy
+O\left(\frac{1}{\mu_i^2R|\log\varepsilon|}\right)\\[1mm]
=&-\frac{p-1}{p}\frac1{\gamma^p}
\int_{|z_i|\leq\frac14\mu_i|\log\varepsilon|^\tau}
\frac{8a(\varepsilon y)}{(\mu_i^2+|y-\xi'_i|^2)^2}
\left[Z_0\left(\frac{z_i}{\mu_i}\right)\right
]^2
\left[1+\omega^{1}_{\mu_i}+\frac{1}{2}(\omega_{\mu_i})^2+2\omega_{\mu_i}\right]\left(y-\xi'_i\right)
dy
+O\left(\frac{1}{\mu_i^2R|\log\varepsilon|}\right).
\endaligned
$$
Owing to the  relation regarding $\omega^1_{\mu_i}$ (see Lemma  A.2 of \cite{Z1})
\begin{equation}\label{3.57a}
\aligned
\int_{\mathbb{R}^2}
\frac{8}{(\mu_i^2+|z|^2)^2}
\left[Z_0\left(\frac{z}{\mu_i}\right)\right
]^2
\left[
1+\omega^{1}_{\mu_i}+\frac{1}{2}(\omega_{\mu_i})^2+2\omega_{\mu_i}
\right](z)dz=\frac{8\pi}{\mu_i^2},
\endaligned
\end{equation}
 by (\ref{2.6}) and (\ref{3.43}) we find
\begin{equation}\label{3.58}
\aligned
I=
\frac{1-p}{4}\frac{c_i a(\xi_i)}{\mu_i^2|\log\varepsilon|}
+O\left(\frac{1}{\mu_i^2R|\log\varepsilon|}\right),
\quad\forall\,\,i=1,\ldots,m.
\endaligned
\end{equation}
Combining  estimates (\ref{3.581}) and (\ref{3.58}), we arrive at
\begin{equation}\label{3.59}
\aligned
\int_{\Omega_\varepsilon}a(\varepsilon y)\widetilde{Z}_{i0}\mathcal{L}(\widetilde{Z}_{i0})
=K+I=\frac{2-p}{4}\frac{c_i a(\xi_i)}{\mu_i^2|\log\varepsilon|}\left[1+O\left(\frac1{R}\right)\right],
\quad\forall\,\,i=1,\ldots,m.
\endaligned
\end{equation}

According to (\ref{3.56}), we still need  to
calculate
$\int_{\Omega_\varepsilon}a(\varepsilon y)\widetilde{Z}_{k0}\mathcal{L}(\widetilde{Z}_{i0})$
with $k\neq i$.
From the previous estimates of
$\mathcal{L}(\widetilde{Z}_{i0})$  and $\widetilde{Z}_{k0}$, we can easily compute
$$
\aligned
\int_{\Omega_1}a(\varepsilon y)\widetilde{Z}_{k0}\mathcal{L}(\widetilde{Z}_{i0})
=O\left(\frac{R^2(\log^2\mu_i)
\log|\log\varepsilon|}{\mu_i\mu_k|\log\varepsilon|^2}\right),
\,\qquad\quad\quad\,\,
\int_{\Omega_2}a(\varepsilon y)\widetilde{Z}_{k0}\mathcal{L}(\widetilde{Z}_{i0})
=O\left(\frac{\log|\log\varepsilon|}{\mu_i\mu_k|\log\varepsilon|^2}\right),
\endaligned
$$
$$
\aligned
\int_{\Omega_4}a(\varepsilon y)\widetilde{Z}_{k0}\mathcal{L}(\widetilde{Z}_{i0})=
O\left(\frac{|\log d|^2}{\mu_i\mu_k|\log\varepsilon|^2}\right),
\,\,\,\,\quad\qquad\,\,\,\,
\int_{\Omega_{3,i}\cup\widetilde{\Omega}_{3}}a(\varepsilon y)\widetilde{Z}_{k0}\mathcal{L}(\widetilde{Z}_{i0})
=O\left(\frac{\log^3|\log\varepsilon|}{\mu_i\mu_k|\log\varepsilon|^2}\right),
\endaligned
$$
and
$$
\aligned
\int_{\Omega_{3,l}}a(\varepsilon y)\widetilde{Z}_{k0}\mathcal{L}(\widetilde{Z}_{i0})
=O\left(\frac{\log^2 |\log\varepsilon|}{\mu_i\mu_k|\log\varepsilon|^2}\right)
\quad\quad\textrm{for all}\,\,\,l\neq i\,\,\,\textrm{and}\,\,\,l\neq k.
\endaligned
$$
It remains to consider the integral over $\Omega_{3,k}$. Using (\ref{3.26}) and  an integration by parts,  we have
$$
\aligned
\int_{\Omega_{3,k}}a(\varepsilon y)\widetilde{Z}_{k0}\mathcal{L}(\widetilde{Z}_{i0})
=\int_{\Omega_{3,k}}a(\varepsilon y)\widetilde{Z}_{i0}\mathcal{L}(\widetilde{Z}_{k0})
-\int_{\partial\Omega_{3,k}}a(\varepsilon y)\widehat{Z}_{k0}\frac{\partial\widehat{Z}_{i0}}{\partial\nu}
+\int_{\partial\Omega_{3,k}}a(\varepsilon y)\widehat{Z}_{i0}\frac{\partial\widehat{Z}_{k0}}{\partial\nu}.
\endaligned
$$
As above, we get
$$
\aligned
\int_{\Omega_{3,k}}a(\varepsilon y)\widetilde{Z}_{i0}\mathcal{L}(\widetilde{Z}_{k0})
=O\left(\frac{\log^3|\log\varepsilon|}{\mu_i\mu_k|\log\varepsilon|^2}\right).
\endaligned
$$
On $\partial\Omega_{3,k}$, by (\ref{2.4}) and (\ref{3.23}),
$$
\aligned
\,\,\,\,\,\,\,\widehat{Z}_{k0}
=O\left(\frac{1}{\mu_k}\right),
\,\,\qquad\,\qquad\ \qquad\,\qquad\,\,
\widehat{Z}_{i0}
=O\left(\frac{\log |\log\varepsilon|}{\mu_i|\log\varepsilon|}\right),
\quad
\endaligned
$$
and
$$
\aligned
|\nabla\widehat{Z}_{k0}|
=O\left(\frac{1}{\mu_k^2|\log\varepsilon|^{1+\tau}}\right),
\,\qquad\,\,\,\qquad\,
|\nabla\widehat{Z}_{i0}|
=O\left(\frac{\,\varepsilon |\log\varepsilon|^{\kappa-1}\,}{\mu_i}\right).
\endaligned
$$
Then
$$
\aligned
\int_{\Omega_{3,k}}a(\varepsilon y)\widetilde{Z}_{k0}\mathcal{L}(\widetilde{Z}_{i0})
=O\left(\frac{\log^3|\log\varepsilon|}{\mu_i\mu_k|\log\varepsilon|^2}\right).
\endaligned
$$
From the above estimates  we find
\begin{equation}\label{3.60}
\aligned
\int_{\Omega_\varepsilon}a(\varepsilon y)\widetilde{Z}_{k0}\mathcal{L}(\widetilde{Z}_{i0})
=O\left(\frac{\log^3 |\log\varepsilon|}{\mu_i\mu_k|\log\varepsilon|^2}\right),\,
\,\quad\,\,\textrm{if}\,\,\,\,i\neq k.
\endaligned
\end{equation}
Furthermore, substituting (\ref{3.59})-(\ref{3.60})
into (\ref{3.56}), we obtain
\begin{equation*}\label{}
\aligned
\frac{|d_i|}{\mu_i}
\leq
C|\log\varepsilon|\,\|h\|_{*}
+C\sum_{k=1}^m\frac{|d_k|}{\mu_k}\frac{\log^4|\log\varepsilon|}{|\log\varepsilon|}.
\endaligned
\end{equation*}
Using linear algebra arguments, we  can conclude Claim 2 for $d_i$ and complete the proof by (\ref{3.32}).
\end{proof}

\vspace{1mm}
\vspace{1mm}
\vspace{1mm}

\noindent {\bf
Proof of Proposition 3.2.}
Let us first prove that
for any $\phi$, $c_{ij}$ solutions of
problem (\ref{3.9}), the a priori estimate
\begin{equation}\label{3.62}
\aligned
\|\phi\|_{L^{\infty}(\Omega_\varepsilon)}\leq C|\log\varepsilon|\,\|h\|_{*}
\endaligned
\end{equation}
holds. In fact,  Lemma  3.5 gives
$$
\aligned
\|\phi\|_{L^{\infty}(\Omega_\varepsilon)}\leq C|\log\varepsilon|\,\left(
\|h\|_{*}+\sum_{i=1}^m\sum_{j=1}^{J_i}|c_{ij}|\cdot\|\chi_iZ_{ij}\|_*
\right)\leq C|\log\varepsilon|\,\left(
\|h\|_{*}+\sum_{i=1}^m\sum_{j=1}^{J_i}\mu_i|c_{ij}|
\right).
\endaligned
$$
As in Lemma 3.4, arguing by contradiction to (\ref{3.62}), we
assume further that
\begin{equation}\label{3.63}
\aligned
\|\phi_n\|_{L^{\infty}(\Omega_{\varepsilon_n})}=1,
\,\,\ \,\quad\,\,|\log\varepsilon_n|\,\|h_n\|_{*}\rightarrow0,
\,\,\ \,\quad\,\,
|\log\varepsilon_n|\,\sum_{i=1}^m\sum_{j=1}^{J_i}\mu_i^n|c_{ij}^n|\geq\delta>0
\,\,\,\quad\,\,\,
\textrm{as}\,\,\,n\rightarrow+\infty.
\endaligned
\end{equation}
We omit the dependence on $n$.
It suffices to prove that $|c_{ij}|\leq C\mu_i^{-1}\|h\|_{*}$.
To this end, we multiply  (\ref{3.9}) by $a(\varepsilon y)\eta_{i2}Z_{ij}$,
with $\eta_{i2}$ the cut-off function
defined in (\ref{3.24})-(\ref{3.25}),
and integrate by parts to find that for
any
$i=1,\ldots,m$ and $j=1,J_i$,
\begin{equation}\label{3.64}
\aligned
\int_{\Omega_\varepsilon}
a(\varepsilon y)\phi \mathcal{L}(\eta_{i2}Z_{ij})=
\int_{\Omega_\varepsilon}
a(\varepsilon y)h\eta_{i2}Z_{ij}+
\sum_{k=1}^m\sum_{t=1}^{J_k}c_{kt}
\int_{\Omega_\varepsilon}\chi_kZ_{kt}\eta_{i2}Z_{ij}.
\endaligned
\end{equation}
From (\ref{3.2}), (\ref{3.6}), (\ref{3.10}), (\ref{3.13}), (\ref{3.43}) and (\ref{3.42}) we can compute
$$
\aligned
\mathcal{L}(\eta_{i2}Z_{ij})=&
\eta_{i2}\mathcal{L}(Z_{ij})
-Z_{ij}\Delta_{a(\varepsilon y)}\eta_{i2}
-2\nabla\eta_{i2}\nabla Z_{ij}
\\[1mm]
=&\left[\frac{8\mu_i^2}{(\mu_{i}^2+|y-\xi'_i|^2)^2}-W_{\xi'}\right]\eta_{i2}Z_{ij}
+\varepsilon^2\eta_{i2}Z_{ij}
+\eta_{i2}\left[-\Delta_{a(\varepsilon y)}Z_{ij}
-\frac{8\mu_i^2}{(\mu_{i}^2+|y-\xi'_i|^2)^2}Z_{ij}\right]
+O\left(\frac{\,\varepsilon^3\,}{d^3}\right)\\[1mm]
\equiv&\mathcal{B}
+\varepsilon^2\eta_{i2}Z_{ij}
+O\left(
\frac{\varepsilon}{(|y-\xi'_i|+\mu_i)^2}
\right)
+O\left(\frac{\,\varepsilon^3\,}{d^3}\right).
\endaligned
$$
For the estimation of the first term, we decompose  $\supp(\eta_{i2})$ into several pieces:
$$
\aligned
\widehat{\Omega}_{1k}=\left\{y\in\Omega_\varepsilon\big|\,|z_k|\leq
\frac14\mu_k|\log\varepsilon|^\tau\right\},
\,\,\ \,\,\forall\,\,k=1,\ldots,m,
\endaligned
$$
$$
\aligned
\widehat{\Omega}_{2}=\left\{y\in\Omega_\varepsilon\big|\,
|z_i|\leq\frac{6d}{\varepsilon},
\,\quad\,|z_k|>
\frac14\mu_k|\log\varepsilon|^\tau,
\,\,\,\,\,\,\,\,\,k=1,\ldots,m\right\},
\endaligned
$$
where $z_k=y-\xi'_k$  for $k=1,\ldots,l$, but
$z_k=F_k^\varepsilon(y)$  for $k=l+1,\ldots,m$.
From (\ref{2.4}),  (\ref{2.6})  and (\ref{3.43}) we obtain
\begin{equation}\label{3.65}
\aligned
|y-\xi'_i|\geq |\xi_i'-\xi_k'|-|y-\xi_k'|\geq
|\xi_i'-\xi_k'|-\frac12\mu_k|\log\varepsilon|^\tau
\geq
\frac{1}{2\varepsilon |\log\varepsilon|^\kappa},
\endaligned
\end{equation}
uniformly in $\widehat{\Omega}_{1k}$, $k\neq i$.
In $\widehat{\Omega}_{1i}$, by (\ref{3.2}),
(\ref{3.6})  (\ref{3.10}) and (\ref{3.13})
we have that for any $i=1,\ldots,l$ and $j=1,2$,
$$
\aligned
\mathcal{B}=-
\frac{8\mu^2_i(y-\xi'_i)_j}{(\mu_{i}^2+|y-\xi'_i|^2)^{3}}
\left\{
\frac{p-1}{p}\frac1{\gamma^p}
\left[1+\omega^{1}_{\mu_i}+\frac{1}{2}(\omega_{\mu_i})^2+2\omega_{\mu_i}\right]\left(y-\xi'_i\right)
+O\left(\frac{\log^4(\mu_i+|y-\xi'_i|)}{\gamma^{2p}}\right)
\right\},
\endaligned
$$
and for any $i=l+1,\ldots,m$ and $j=J_i=1$,
$$
\aligned
\mathcal{B}=-
\frac{8\mu^2_i}{(\mu_{i}^2+|y-\xi'_i|^2)^{2}}
\frac{(F_i^\varepsilon(y))_1}{\mu_i^2+|F_i^\varepsilon(y)|^2}
\left\{
\frac{p-1}{p}\frac1{\gamma^p}
\left[1+\omega^{1}_{\mu_i}+\frac{1}{2}(\omega_{\mu_i})^2+2\omega_{\mu_i}\right]\left(y-\xi'_i\right)
+O\left(\frac{\log^4(\mu_i+|y-\xi'_i|)}{\gamma^{2p}}\right)
\right\}.
\endaligned
$$
In $\widehat{\Omega}_{1k}$, $k\neq i$, by   (\ref{3.65}),
$$
\aligned
\mathcal{B}=
O\left(\frac{8\varepsilon\mu^2_k |\log\varepsilon|^\kappa}
{(\mu_{k}^2+|y-\xi'_k|^2)^{2}}\right).
\endaligned
$$
In $\widehat{\Omega}_{2}$, by (\ref{3.1}),
$$
\aligned
\mathcal{B}=
\sum_{k=1}^mO\left(\left(\frac{\mu_k^2}{|y-\xi'_k|^4}\right)^{1+
o\left(1\right)}\cdot\frac{1}{\mu_i|\log\varepsilon|^\tau}
\right).
\endaligned
$$
Hence
\begin{equation}\label{3.66}
\aligned
\int_{\Omega_\varepsilon}
a(\varepsilon y)\phi \mathcal{L}(\eta_{i2}Z_{ij})=-\frac{p-1}{p}\frac1{\gamma^p}\frac{1}{\mu_i}a(\xi_i)E_{j}(\widehat{\phi}_i)
+O\left(\frac{1}{\mu_i|\log\varepsilon|^2}\|\phi\|_{L^{\infty}(\Omega_\varepsilon)}\right),
\endaligned
\end{equation}
where for any $i=1,\ldots,l$ and $j=1,2$, $\widehat{\phi}_i(z)=\phi\big(\xi_i'+\mu_iz\big)$ and
$$
\aligned
E_{j}(\widehat{\phi}_i)=\int_{B_{\frac18|\log\varepsilon|^\tau}\left(0\right)}
\frac{8z_j}
{(|z|^{2}+1)^3}\widehat{\phi}_i(z)\left[1+\omega^{1}_{\mu_i}+\frac{1}{2}(\omega_{\mu_i})^2+2\omega_{\mu_i}\right](\mu_i|z|)
dz,
\endaligned
$$
but for any $i=l+1,\ldots,m$ and $j=1$,
$\widehat{\phi}_i(z)=\phi\big((F_i^\varepsilon)^{-1}(\mu_iz)\big)$ and
$$
\aligned
E_{j}(\widehat{\phi}_i)=\int_{\mathbb{R}_+^2\bigcap
B_{\frac18|\log\varepsilon|^\tau}\left(0\right)}
\frac{8z_j}
{(|z|^{2}+1)^3}\widehat{\phi}_i(z)\left[1+\omega^{1}_{\mu_i}+\frac{1}{2}(\omega_{\mu_i})^2+2\omega_{\mu_i}\right](\mu_i|z|)
dz.
\endaligned
$$
On the other hand, if $1\leq k=i\leq l$,
\begin{equation}\label{3.68}
\aligned
\int_{\Omega_\varepsilon}\chi_kZ_{kt}\eta_{i2}Z_{ij}=
\int_{\mathbb{R}^2}\chi Z_{t} Z_{j}dz=D_t\delta_{tj},
\endaligned
\end{equation}
and if $l+1\leq k=i\leq m$,
\begin{equation}\label{3.69}
\aligned
\int_{\Omega_\varepsilon}\chi_kZ_{k1}\eta_{i2}Z_{i1}=
\int_{\mathbb{R}_{+}^2}\chi Z^2_{1} \big[1+O\left(\varepsilon\mu_i|z|\right)\big] dz=\frac12D_1\big[1+O\left(\varepsilon\mu_i\right)\big],
\endaligned
\end{equation}
and  if $k\neq i$,  by   (\ref{3.65}),
\begin{equation}\label{3.70}
\aligned
\int_{\Omega_\varepsilon}\chi_kZ_{kt}\eta_{i2}Z_{ij}=O\left(\mu_k\varepsilon |\log\varepsilon|^\kappa\right).
\endaligned
\end{equation}
In addition, due to
$\|\eta_{i2}Z_{ij}\|_{L^{\infty}(\Omega_\varepsilon)}\leq C\mu_i^{-1}$, we obtain
\begin{equation}\label{3.67}
\aligned
\int_{\Omega_\varepsilon}a(\varepsilon y)h\eta_{i2}Z_{ij}=O\left(\frac{1}{\mu_i}\|h\|_{*}\right).
\endaligned
\end{equation}
As a consequence,  replacing estimates (\ref{3.66})-(\ref{3.67}) to  (\ref{3.64}),  we have that
for any $i=1,\ldots,m$ and $j=1,J_i$,
$$
\aligned
D_jc_{ij}+O\left(\sum\limits_{k\neq i}^m\sum\limits_{t=1}^{J_k}\varepsilon\mu_k |\log\varepsilon|^\kappa|c_{kt}|
\right)
= O\left(
\frac{1}{\mu_i}\|h\|_{*}
+\frac1{\mu_i|\log\varepsilon|}\|\phi\|_{L^{\infty}(\Omega_\varepsilon)}
\right).
\endaligned
$$
Then
\begin{equation}\label{3.71}
\aligned
\sum_{i=1}^m\sum_{j=1}^{J_i}\mu_i|c_{ij}|
= O\left(
\frac1{|\log\varepsilon|}\|\phi\|_{L^{\infty}(\Omega_\varepsilon)}
\right)+
 O\left(
\|h\|_{*}
\right).
\endaligned
\end{equation}
From the first two assumptions in (\ref{3.63}) we get
$\sum\limits_{i=1}^m\sum\limits_{j=1}^{J_i}\mu_i|c_{ij}|=o\left(1\right)$.
As in contradiction arguments of
Lemma 3.4, we can derive that for any $i=1,\ldots,l$,
$$
\aligned
\widehat{\phi}_i\rightarrow
C_i\frac{|z|^2-1}{|z|^2+1}
\,\,\ \,\,\,
\,\,\textrm{uniformly in}\,\,\,C_{loc}^0(\mathbb{R}^2),
\endaligned
$$
but for any $i=l+1,\ldots,m$,
$$
\aligned
\widehat{\phi}_i\rightarrow
C_i\frac{|z|^2-1}{|z|^2+1}
\,\,\ \,\,\,
\,\,\textrm{uniformly in}\,\,\,C_{loc}^0(\mathbb{R}_{+}^2),
\endaligned
$$
with some constant $C_i$.
In view of the odd function  $\frac{8z_j}
{(|z|^{2}+1)^3}$ with $j=1,2$, by (\ref{8.2}) and
 Lebesgue's theorem we have that
$$
\aligned
E_{j}(\widehat{\phi}_i)
\longrightarrow
0,
\quad\,\,\,\,\,\forall\,\,\,
i=1,\ldots,m,\,\,\,j=1,J_i.
\endaligned
$$
Hence by replacing estimates (\ref{3.66})-(\ref{3.67}) to  (\ref{3.64}) again  we have a better estimate
$$
\aligned
\sum_{i=1}^m\sum_{j=1}^{J_i}\mu_i|c_{ij}|=o\left(\frac1{|\log\varepsilon|}\right)+O(\|h\|_{*}),
\endaligned
$$
which is impossible because of  the last assumption in (\ref{3.63}). So estimate (\ref{3.62})
is established and then by (\ref{3.71}),  we find
\begin{equation*}\label{3.72}
\aligned
|c_{ij}|\leq C\frac{1}{\mu_i}\|h\|_{*}.
\endaligned
\end{equation*}
\indent Let us consider the Hilbert space
$$
\aligned
H_{\xi}=\left\{\phi\in H^1(\Omega_\varepsilon)\left|\,\,
\int_{\Omega_\varepsilon}\chi_iZ_{ij}\phi=0
\quad
\textrm{for any}\,\,\,i=1,\ldots,m,\,\,j=1,J_i;
\quad
\frac{\partial\phi}{\partial\nu}=0\quad\textrm{on}\,\,\,\po_\varepsilon
\right.
\right\}
\endaligned
$$
with the norm
$\|\phi\|_{H^1(\Omega_\varepsilon)}^2=\int_{\Omega_\varepsilon}a(\varepsilon y)\big(|\nabla\phi|^2+\varepsilon^2\phi^2\big)$.
Equation (\ref{3.9}) is equivalent to find
$\phi\in H_\xi$ such that
$$
\aligned
\int_{\Omega_\varepsilon}a(\varepsilon y)\big(\nabla\phi\nabla\psi+\varepsilon^2\phi\psi\big)
-\int_{\Omega_\varepsilon}a(\varepsilon y)W_{\xi'}\phi\psi
=\int_{\Omega_\varepsilon}a(\varepsilon  y)h\psi\,\,\quad\,\,\forall
\,\,\psi\in H_\xi.
\endaligned
$$
By Fredholm's alternative this is equivalent to the uniqueness of solutions to this
problem, which in turn follows from estimate (\ref{3.62}).
The proof is complete.
\qquad\qquad\qquad\qquad\qquad\qquad\qquad\qquad\qquad\qquad\quad
\qquad\qquad\qquad\qquad\qquad\qquad\qquad\qquad\qquad\quad
$\square$

\vspace{1mm}
\vspace{1mm}
\vspace{1mm}
\vspace{1mm}

\noindent{\bf Remark 3.6.}\,\,Given $h\in L^\infty(\Omega_\varepsilon)$ with
$\|h\|_*<\infty$, let $\phi$ be the solution of equation (\ref{3.9}) given by Proposition 3.2.
Multiplying  (\ref{3.9}) by $a(\varepsilon y)\phi$ and integrating by parts, we get
$$
\aligned
\|\phi\|_{H^1(\Omega_\varepsilon)}^2
=\int_{\Omega_\varepsilon}a(\varepsilon y)W_{\xi'}\phi^2
+\int_{\Omega_\varepsilon}a(\varepsilon  y)h\phi.
\endaligned
$$
By Proposition 3.1 we find
\begin{equation*}\label{3.73}
\aligned
\|\phi\|_{H^1(\Omega_\varepsilon)}\leq C\big(\|\phi\|_{L^{\infty}(\Omega_\varepsilon)}+\|h\|_*\big).
\endaligned
\end{equation*}

\vspace{1mm}
\vspace{1mm}
\vspace{1mm}

\noindent{\bf Remark 3.7.}\,\,The result of Proposition 3.2 implies that the unique solution $\phi=\mathcal{T}(h)$
of equation (\ref{3.9}) defines a continuous linear map from the Banach space
$\mathcal{C}_*$ of all functions $h$ in $L^\infty$ for which $\|h\|_{*}<\infty$,
into $L^\infty$. It is necessary to point out that
the operator $\mathcal{T}$ is differentiable  with
respect to the variables $\xi=(\xi_1,\ldots,\xi_m)$ in $\mathcal{O}_\varepsilon$.
More precisely, if  we fix $h\in\mathcal{C}_*$ and set $\phi=T(h)$, then by formally
computing the derivative of $\phi$ with respect to $\xi'=(\xi'_1,\ldots,\xi'_m)$
and using the delicate estimate $\|\partial_{(\xi'_{i})_j}W_{\xi'}\|_{*}=O\left(1\right)$
we can obtain the a priori  estimate
\begin{equation*}\label{3.74}
\aligned
\|\partial_{(\xi'_{i})_j}\mathcal{T}(h)\|_{L^{\infty}(\Omega_\varepsilon)}\leq C|\log\varepsilon|^2\,\|h\|_{*},
\quad\,\,\,\,\,\forall\,\,\,
i=1,\ldots,m,\,\,\,j=1,J_i.
\endaligned
\end{equation*}
%In fact, the proof of (\ref{3.74}) is similar to that in Lemma 4.6 of \cite{D},
%but the  delicate ingredient is to estimate $\|\partial_{(\xi'_{i})_j}W_{\xi'}\|_{*}$.
%Notice that $\partial_{(\xi'_{i})_j}W_{\xi'}=f''(V_{\xi'})\partial_{(\xi'_{i})_j}V_{\xi'}$.
%Similar to the proof in Lemma 2.1, by (\ref{2.7}) and  (\ref{2.23})
%we can derive that $\partial_{(\xi'_{i})_j}\big[p\gamma^{p-1}H_k(\varepsilon y)\big]=o\big(1\big)$
%uniformly over $\overline{\Omega}$. Furthermore, by (\ref{2.16}) and (\ref{2.28}) we find
%\begin{equation}\label{3.91}
%\aligned
%\partial_{(\xi'_{i})_j}V_{\xi'}(y)=\frac{4}{\mu_i}Z_j\left(\frac{y-\xi_i'}{\mu_i}\right)
%+o\big(1\big).
%\endaligned
%\end{equation}
%This, combined with the facts that $\frac{1}{\mu_i}\leq C$
%and $\|f''(V_{\xi'})\|_{*}=O\left(1\right)$,  implies
%\begin{equation}\label{3.92}
%\aligned
%\|\partial_{(\xi'_{i})_j}V_{\xi'}\|_{L^{\infty}(\Omega_\varepsilon)}=O\left(1\right)
%\,\,\qquad\ \ \ \,\,
%\textrm{and}
%\,\,\qquad\ \ \ \,\,
%\|\partial_{(\xi'_{i})_j}W_{\xi'}\|_{*}=O\left(1\right).
%\endaligned
%\end{equation}

\vspace{1mm}
\vspace{1mm}

\section{The nonlinear problem}
Consider the nonlinear problem: for any points
$\xi=(\xi_1,\ldots,\xi_m)\in\mathcal{O}_\varepsilon$,
we find a function $\phi$ and scalars $c_{ij}$,
$i=1,\ldots,m$, $j=1,J_i$  such that
\begin{equation}\label{4.1}
\left\{\aligned
&\mathcal{L}(\phi)=-\big[
E_{\xi'}+N(\phi)
\big]
+\frac1{a(\varepsilon y)}\sum\limits_{i=1}^m\sum\limits_{j=1}^{J_i}c_{ij}\chi_i\,Z_{ij}
\,\,\ \,\textrm{in}\,\,\,\,\,\,\Omega_\varepsilon,\\
&\frac{\partial\phi}{\partial\nu}=0\,\,\,\,\,\,\,
\qquad\qquad\qquad\quad\qquad\qquad\qquad\qquad
\qquad
\textrm{on}\,\,\,\,\partial\Omega_{\varepsilon},\\[1mm]
&\int_{\Omega_\varepsilon}\chi_i\,Z_{ij}\phi=0
\qquad\qquad\qquad\quad\,
\forall\,\,i=1,\ldots,m,\,\,\,j=1, J_i,
\endaligned\right.
\end{equation}
where $W_{\xi'}=f'(V_{\xi'})$ satisfies (\ref{3.1})-(\ref{3.3}),
and $E_{\xi'}$, $N(\phi)$ are defined in
(\ref{2.31}). We have the following result.

\vspace{1mm}
\vspace{1mm}
\vspace{1mm}
\vspace{1mm}

\noindent{\bf Proposition 4.1.}\,\,{\it
Let $m$ be a positive integer.
Then there exist constants $C>0$ and  $\varepsilon_0>0$ such
that for any $0<\varepsilon<\varepsilon_0$ and any points
$\xi=(\xi_1,\ldots,\xi_m)\in\mathcal{O}_\varepsilon$,
problem {\upshape(\ref{4.1})} admits
a unique solution
$\phi_{\xi'}$
for some coefficients $c_{ij}(\xi')$,
$i=1,\ldots,m$, $j=1,J_i$, such that
\begin{equation}\label{4.2}
\aligned
\|\phi_{\xi'}\|_{L^{\infty}(\Omega_\varepsilon)}\leq\frac{C}{|\log\varepsilon|^3},\,\,\quad\,\quad\,\,
\sum_{i=1}^m\sum_{j=1}^{J_i}\mu_i|c_{ij}(\xi')|\leq\frac{C}{|\log\varepsilon|^4}\quad\,\quad\textrm{and}\quad\,\quad
\|\phi_{\xi'}\|_{H^1(\Omega_\varepsilon)}\leq\frac{C}{|\log\varepsilon|^3}.
\endaligned
\end{equation}
Furthermore, the map $\xi'\mapsto\phi_{\xi'}$ is a $C^1$-function in $C(\overline{\Omega}_\varepsilon)$ and $H^1(\Omega_\varepsilon)$,
precisely for any $i=1,\ldots,m$ and $j=1,J_i$,
\begin{equation}\label{4.3}
\aligned
\|\partial_{(\xi'_{i})_j}\phi_{\xi'}\|_{L^{\infty}(\Omega_\varepsilon)}\leq\frac{C}{|\log\varepsilon|^2},
\endaligned
\end{equation}
where $\xi':=(\xi'_1,\ldots,\xi'_m)=(\frac1{\varepsilon}\xi_1,\ldots,\frac1{\varepsilon}\xi_m)$.
}

\vspace{1mm}

\begin{proof}
Proposition $3.2$ and Remarks 3.6-3.7 allow us to apply the Contraction Mapping
Theorem
and the Implicit Function Theorem
to find a unique solution for problem (\ref{4.1})
satisfying (\ref{4.2})-(\ref{4.3}).
Since it is a standard procedure, we omit the details,
see Lemmas 4.1-4.2 in \cite{DKM} for a similar proof.
We just mention that
$\|N(\phi)\|_{*}\leq
C\|\phi\|_{L^{\infty}(\Omega_\varepsilon)}^2$
and
$\|\partial_{(\xi'_{i})_j}E_{\xi'}\|_{*}
\leq
C|\log\varepsilon|^{-3}$.
\end{proof}

\vspace{1mm}

\section{Variational reduction}
Since problem (\ref{4.1}) has been solved, we  find a solution of problem (\ref{2.29})
and hence to the original equation  (\ref{1.1}) if we match  $\xi'$ with
the coefficient $c_{ij}(\xi')$ in (\ref{4.1}) so  that
\begin{equation}\label{5.1}
\aligned
c_{ij}(\xi')=0\,\quad\,\,\textrm{for all}\,\,\,i=1,\ldots,m,\,\,\,j=1,J_i.
\endaligned
\end{equation}
We consider the  functional $J_\lambda$ defined in (\ref{1.2})
and
take its finite-dimensional restriction
\begin{equation*}\label{5.2}
\aligned
F_\lambda(\xi)=J_\lambda\left(\big(U_\xi+\widetilde{\phi}_{\xi}\big)(x)
\right)\,\ \,\,\ \ \,\forall\,\,\xi=(\xi_1,\ldots,\xi_m)\in\mathcal{O}_\varepsilon,
\endaligned
\end{equation*}
where
\begin{equation*}\label{5.3}
\aligned
\big(U_\xi+\widetilde{\phi}_{\xi}\big)(x)=\gamma+\frac1{p\gamma^{p-1}}\big(V_{\xi'}+\phi_{\xi'}\big)
\left(\frac{x}{\varepsilon}\right),
\,\quad\,x\in\Omega,
\endaligned
\end{equation*}
with  $V_{\xi'}$   defined in (\ref{2.28}) and   $\phi_{\xi'}$ the unique solution to problem
(\ref{4.1}) given by Proposition 4.1.
Define
\begin{equation*}\label{5.4}
\aligned
I_\varepsilon(\upsilon)=\frac12\int_{\Omega_\varepsilon}a(\varepsilon y)\left(
|\nabla \upsilon|^2+\varepsilon^2\upsilon^2
\right)dy-\int_{\Omega_\varepsilon}a(\varepsilon y)e^{\gamma^p\left[\left(1+\frac{\upsilon}{p\gamma^p}\right)^p-1\right]}dy,
\,\quad\,\upsilon\in H^1(\Omega_\varepsilon).
\endaligned
\end{equation*}
Then by (\ref{1.5})-(\ref{2.7}),
\begin{equation}\label{5.5}
\aligned
I_\varepsilon\big(V_{\xi'}+\phi_{\xi'}\big)=
p^2\gamma^{2(p-1)}F_\lambda(\xi)
\qquad\,\,
\textrm{and}
\qquad\,\,
I_\varepsilon\big(V_{\xi'}+\phi_{\xi'}\big)-I_\varepsilon(V_{\xi'}
)=
p^2\gamma^{2(p-1)}\big[F_\lambda(\xi)
-
J_\lambda\left(U_\xi\right)\big].
\endaligned
\end{equation}

\vspace{1mm}
\vspace{1mm}

\noindent{\bf Proposition 5.1.}\,\,{\it The function $F_\lambda:\mathcal{O}_\varepsilon\mapsto\mathbb{R}$
is of class $C^1$.
Moreover, for all $\lambda$ sufficiently small,
if  $D_{\xi}F_\lambda(\xi)=0$,  then $\xi'=\xi/\varepsilon$ satisfies {\upshape (\ref{5.1})},
that is,
$U_\xi+\widetilde{\phi}_{\xi}$ is a solution of equation {\upshape(\ref{1.1})}.
}

\vspace{1mm}

\begin{proof}
Since the map
$\xi'\mapsto\phi_{\xi'}$ is a $C^1$-function in $C(\overline{\Omega}_\varepsilon)$ and $H^1(\Omega_\varepsilon)$,
we can check that $F_\lambda(\xi)$ is a $C^1$-function of $\xi$ in $\mathcal{O}_\varepsilon$.
Assume that   $\phi_{\xi'}$ solves problem (\ref{4.1})
and $D_{\xi}F_\lambda(\xi)=0$. Then  by (\ref{5.5}), we have that for any $k=1,\ldots,m$ and $t=1,J_k$,
\begin{eqnarray}\label{5.6}
&&0
=I'_\varepsilon\big(V_{\xi'}+\phi_{\xi'}\big)\partial_{(\xi'_k)_t}\big(V_{\xi'}+\phi_{\xi'}\big)\nonumber\\[1mm]
&&\,\,\,\,\,=\sum\limits_{i=1}^m\sum\limits_{j=1}^{J_i}c_{ij}(\xi')\int_{\Omega_\varepsilon}\chi_i Z_{ij}
\partial_{(\xi'_k)_t}V_{\xi'}
-\sum\limits_{i=1}^m\sum\limits_{j=1}^{J_i}c_{ij}(\xi')\int_{\Omega_\varepsilon}\phi_{\xi'}
\partial_{(\xi'_k)_t}\big(\chi_i Z_{ij}\big).
\end{eqnarray}
Recall  that  $D_{\xi'}V_{\xi'}(y)=p\gamma^{p-1}D_{\xi'}U_\xi(\varepsilon y)$.
From (\ref{2.1}), (\ref{2.3}), (\ref{2.8})  and (\ref{2.16})
we know that
$$
\aligned
\partial_{(\xi'_k)_t}V_{\xi'}(y)
=\sum_{i=1}^{m}
\partial_{(\xi'_k)_t}\left[
\omega_{\mu_i}\left(y-\xi'_i\right)
+\sum_{j=1}^4\left(\frac{p-1}{p}\right)^j\frac{1}{\gamma^{jp}}
\omega^j_{\mu_{i}}
\left(y-\xi'_i\right)
+p\gamma^{p-1}H_i(\varepsilon y)\right].
\endaligned
$$
From the fact that
$|\partial_{(\xi'_k)_t}\log\mu_i|=O\left(\varepsilon|\log\varepsilon|^{\kappa}\right)$
for any $i=1,\ldots,m$,
we have that
$$
\aligned
\partial_{(\xi'_k)_t}\omega_{\mu_i}\left(y-\xi'_i\right)=
\frac4{\mu_i}\delta_{ki}
Z_{t}\left(
\frac{y-\xi'_i}{\mu_i}
\right)
+O\left(\varepsilon|\log\varepsilon|^{\kappa}\right),
\endaligned
$$
and for each $j=1,2,3,4$,
$$
\aligned
\partial_{(\xi'_k)_t}\omega^j_{\mu_{i}}
\left(y-\xi'_i\right)=
-\frac1{\mu_i}\delta_{ki}\left[D^j_{\mu_i}Z_{t}\left(
\frac{y-\xi'_i}{\mu_i}
\right)+
O\left(\frac{\mu_i^2}{|y-\xi'_i|^2+\mu_i^2}\right)
\right]
+
O\left(\varepsilon|\log\varepsilon|^{\kappa}\right).
\endaligned
$$
As in the proof of Lemma 2.1, by the elliptic regularity of the equation we can  prove that
$$
\aligned
\partial_{(\xi'_k)_t}\big[p\gamma^{p-1}H_i(\varepsilon y)\big]
=O\left(\varepsilon|\log\varepsilon|^{\kappa}\right).
\endaligned
$$
Then
$$
\aligned
\partial_{(\xi'_k)_t}V_{\xi'}(y)=
\frac{1}{\mu_k}\left\{
\left[
4-\sum_{j=1}^4\left(\frac{p-1}{p}\right)^j\frac{D^j_{\mu_k}}{\gamma^{jp}}
\right]
Z_{t}\left(
\frac{y-\xi'_k}{\mu_k}
\right)
+
O\left(\frac{1}{|\log\varepsilon|}\right)
\right\}
+O\left(\varepsilon|\log\varepsilon|^{\kappa}\right).
\endaligned
$$
So
\begin{equation}\label{5.7}
\aligned
\int_{\Omega_\varepsilon}\chi_i Z_{ij}
\partial_{(\xi'_k)_t}V_{\xi'}=&
\delta_{ki}\left[
\frac{c_k}{2\pi}\delta_{tj}\int_{\mathbb{R}^2}\chi(|z|)Z^2_t(z)dz+O\left(\frac{\log|\log\varepsilon|}{|\log\varepsilon|}\right)
\right]
+(1-\delta_{ki})O\left(
\frac{\mu_i}{|\xi'_i-\xi'_k|}
\right)
+O\left(\varepsilon\mu_i|\log\varepsilon|^{\kappa}\right).
\endaligned
\end{equation}
On the other hand,  by (\ref{3.6}), (\ref{3.10}),  (\ref{3.13})  and (\ref{3.43})
we can compute
$$
\aligned
\big|\partial_{(\xi'_k)_t}\big(\chi_i Z_{ij}\big)\big|
=O\left(
\frac1{\mu_i}\varepsilon|\log\varepsilon|^{\kappa}
+\frac1{\mu_i^2}\delta_{ki}
\right).
\endaligned
$$
Then
\begin{equation}\label{5.8}
\aligned
\int_{\Omega_\varepsilon}\phi_{\xi'}
\partial_{(\xi'_k)_t}\big(\chi_i Z_{ij}\big)=
\|\phi_{\xi'}\|_{L^{\infty}(\Omega_\varepsilon)}O\big(\varepsilon\mu_i|\log\varepsilon|^{\kappa}
+
\delta_{ki}
\big)=O\left(\frac{1}{|\log\varepsilon|^3}\right).
\endaligned
\end{equation}
Hence by (\ref{5.7})-(\ref{5.8}),
equations (\ref{5.6}) can be written as, for
each $k=1,\ldots,m$  and  $t=1,J_k$,
$$
\aligned
c_{kt}(\xi')
\left[\frac{c_k}{2\pi}\int_{\mathbb{R}^2}\chi(|z|)Z_t^2(z)dz
\right]+
\sum\limits_{i=1}^m\sum\limits_{j=1}^{J_i}|c_{ij}(\xi')|
O\left(\frac{\log|\log\varepsilon|}{|\log\varepsilon|}
\right)=0,
\endaligned
$$
which is a strictly diagonal dominant system. This implies that $c_{kt}(\xi')=0$
for each $k=1,\ldots,m$  and  $t=1,J_k$.
\end{proof}

\vspace{1mm}

In order to solve for critical points of the function $F_\lambda$, a delicate ingredient  is the
expected uniformly $C^1$-closeness between  the functions $I_\varepsilon\big(V_{\xi'}+\phi_{\xi'}\big)$
and $I_\varepsilon\big(V_{\xi'}\big)$, which will
be applied in the proof of our main theorems.

\vspace{1mm}
\vspace{1mm}
\vspace{1mm}

\noindent{\bf Proposition 5.2.}\,\,{\it
For any points $\xi=(\xi_1,\ldots,\xi_m)\in\mathcal{O}_\varepsilon$
and for any $\varepsilon$ small enough,
the
following expansion  uniformly holds
$$
\aligned
I_\varepsilon\big(V_{\xi'}+\phi_{\xi'}\big)=
I_\varepsilon(V_{\xi'})
+
\theta_\varepsilon(\xi'),
\endaligned
$$
where
$$
\aligned
|\theta_\varepsilon(\xi')|+\|\nabla\theta_\varepsilon(\xi')\|
=
O\left(\frac{1}{|\log\varepsilon|^6}\right).
\endaligned
$$
}

\begin{proof}
Using $DI_\varepsilon(V_{\xi'}+\phi_{\xi'})[\phi_{\xi'}]=0$, a Taylor expansion and an integration by parts give
$$
\aligned
I_\varepsilon\big(V_{\xi'}+\phi_{\xi'}\big)-I_\varepsilon(V_{\xi'}
)
&=\int_0^1D^2I_\varepsilon(V_{\xi'}+t\phi_{\xi'})\phi_{\xi'}^2(1-t)dt\\
&=\int_0^1\left\{
\int_{\Omega_\varepsilon}
a(\varepsilon y)\big[f'(V_{\xi'})
-f'(V_{\xi'}+t\phi_{\xi'})\big]\phi_{\xi'}^2-
a(\varepsilon y)\big[E_{\xi'}+N(\phi_{\xi'})\big]
\phi_{\xi'}
\right\}(1-t)dt,
\endaligned
$$
so we get
$$
\aligned
\theta_\varepsilon(\xi')=I_\varepsilon\big(V_{\xi'}+\phi_{\xi'}\big)-I_\varepsilon(V_{\xi'})
=O\left(\frac{1}{|\log\varepsilon|^7}\right),
\endaligned
$$
taking into account  $\|\phi_{\xi'}\|_{L^{\infty}(\Omega_\varepsilon)}\leq C|\log\varepsilon|^{-3}$,
$\|E_{\xi'}\|_{*}\leq C|\log\varepsilon|^{-4}$ and
$\|N(\phi_{\xi'})\|_{*}\leq C|\log\varepsilon|^{-6}$
and
(\ref{3.3}).
Let us differentiate with respect to $\xi'$,
$$
\aligned
\partial_{\xi'}\left[I_\varepsilon\big(V_{\xi'}+\phi_{\xi'}\big)-I_\varepsilon(V_{\xi'}
)\right]
=\int_0^1\left\{
\int_{\Omega_\varepsilon}
a(\varepsilon y)\partial_{\xi'}\big\{\big[f'(V_{\xi'})
-f'(V_{\xi'}+t\phi_{\xi'})\big]\phi_{\xi'}^2
-
\big[E_{\xi'}+N(\phi_{\xi'})\big]
\phi_{\xi'}\big\}
\right\}(1-t)dt.
\endaligned
$$
From estimates
$\|\partial_{\xi'}\phi_{\xi'}\|_{L^{\infty}(\Omega_\varepsilon)}\leq C|\log\varepsilon|^{-2}$,
$\|\partial_{\xi'}E_{\xi'}\|_{*}\leq C|\log\varepsilon|^{-3}$,
$\|\partial_{\xi'}N(\phi_{\xi'})\|_{*}\leq C|\log\varepsilon|^{-5}$
and
$\|\partial_{\xi'}W_{\xi'}\|_{*}\leq C$ we find
$$
\aligned
\partial_{\xi'}\theta_\varepsilon(\xi')=\partial_{\xi'}\left[I_\varepsilon\big(V_{\xi'}+\phi_{\xi'}\big)-I_\varepsilon(V_{\xi'}
)\right]
=O\left(\frac{1}{|\log\varepsilon|^6}\right).
\endaligned
$$
The continuity in $\xi'$ of all these expressions is inherited from that of $\phi_{\xi'}$ and its derivatives in
$\xi'$ in the $L^\infty$ norm.
\end{proof}

\section{Expansion of the energy}
In this section we will give an asymptotic  estimate of  $J_\lambda(U_\xi)$
where $U_\xi$ is the approximate solution defined in (\ref{2.16}) and
$J_\lambda$ is the energy functional (\ref{1.2}) associated to problem (\ref{1.1}).

We have

\vspace{1mm}
\vspace{1mm}

\noindent{\bf Proposition 6.1.}\,\,{\it
Let $m$ be a positive integer. With the choice
{\upshape(\ref{2.22})} for the parameters $\mu_i$, there exists $\varepsilon_0>0$
such that for
any  $0<\varepsilon<\varepsilon_0$
and any points $\xi=(\xi_1,\ldots,\xi_m)\in\mathcal{O}_\varepsilon$,
the following expansion uniformly holds
\begin{eqnarray}\label{6.1}
J_\lambda\left(U_\xi\right)=
\frac{1}{2p^2\gamma^{2(p-1)}}
\sum_{i=1}^mc_ia(\xi_i)
\left\{
4\big|\log\varepsilon\big|
+2\log8-4
-
c_i H_a(\xi_i,\xi_i)
-\sum_{k=1,\,k\neq i}^mc_k
G_a(\xi_i,\xi_k)
+O\left(\frac{\log^\beta|\log\varepsilon|}{|\log\varepsilon|}
\right)\right\},
\end{eqnarray}
where $\beta>1$ is a large but  fixed number, independent of $\varepsilon$.
}

\noindent
\begin{proof}
Observe that
\begin{equation}\label{6.2}
\aligned
J_\lambda\left(U_\xi\right)=
\frac12\int_{\Omega}a(x)\big(|\nabla U_\xi|^2+ U_\xi^2\big)
-\frac{\lambda}{p}\int_{\Omega}a(x)e^{(U_\xi)^p}:=I_A-I_B.
\endaligned
\end{equation}
Let us analyze the behavior of $I_A$.
From the definition of $U_\xi$ in (\ref{2.16}) we get
$$
\aligned
I_A=&\frac12\int_{\Omega}a(x)\big(-\Delta_a U_\xi+U_\xi\big)U_\xi dx
=\frac12\sum\limits_{i=1}^{m}\int_{\Omega}a(x)\big(-\Delta U_i\big)U_\xi dx\nonumber\\
=&\frac1{2p\gamma^{p-1}\varepsilon^2}
\sum\limits_{i=1}^{m}\int_{\Omega}a(x)\left[
e^{\omega_{\mu_i}\left(\frac{x-\xi_i}{\varepsilon}\right)}
-\sum_{j=1}^4\left(\frac{p-1}{p}\right)^j\frac{1}{\gamma^{jp}}\Delta\omega^j_{\mu_{i}}
\left(\frac{x-\xi_i}{\varepsilon}\right)\right]
U_{\xi} dx\nonumber\\
=&\frac{1}{2p\gamma^{p-1}\varepsilon^2}
\sum\limits_{i=1}^{m}\int_{\Omega}a(x)
e^{\omega_{\mu_i}\left(\frac{x-\xi_i}{\varepsilon}\right)}\left[
1
+\sum_{j=1}^4\left(\frac{p-1}{p}\right)^j\frac{1}{\gamma^{jp}}
\big(
\omega^j_{\mu_i}-f^j_{\mu_i}
\big)\left(\frac{x-\xi_i}{\varepsilon}\right)
\right]U_\xi dx.
\endaligned
$$
%$$
%\aligned
%I_A=&\frac12\int_{\Omega}a(x)\big(-\Delta_a U_\xi+U_\xi\big)U_\xi dx
%=\frac12\sum\limits_{i=1}^{m}\int_{\Omega}a(x)\big(-\Delta U_i\big)U_\xi dx\nonumber\\
%=&\frac1{2p\gamma^{p-1}\varepsilon^2}
%\sum\limits_{i=1}^{m}\int_{\Omega}a(x)\left[
%e^{\omega_{\mu_i}\left(\frac{x-\xi_i}{\varepsilon}\right)}
%-\sum_{j=1}^4\left(\frac{p-1}{p}\right)^j\frac{1}{\gamma^{jp}}\Delta\omega^j_{\mu_{i}}
%\left(\frac{x-\xi_i}{\varepsilon}\right)\right]
%U_{\xi} dx\nonumber\\
%=&\frac{1}{2p\gamma^{p-1}\varepsilon^2}
%\sum\limits_{i=1}^{m}\int_{\Omega}a(x)
%e^{\omega_{\mu_i}\left(\frac{x-\xi_i}{\varepsilon}\right)}\left[
%1
%+\sum_{j=1}^4\left(\frac{p-1}{p}\right)^j\frac{1}{\gamma^{jp}}
%\big(
%\omega^j_{\mu_i}-f^j_{\mu_i}
%\big)\left(\frac{x-\xi_i}{\varepsilon}\right)
%\right]U_\xi dx.
%\endaligned
%$$
From  (\ref{2.20})-(\ref{2.21})  we can compute
$$
\aligned
I_A=&\frac{1}{2p^2\gamma^{2(p-1)}\varepsilon^2}
\sum\limits_{i,\,k=1}^{m}\int_{\Omega\cap B_{1/|\log\varepsilon|^{2\kappa}}(\xi_k)}a(x)
e^{\omega_{\mu_i}\left(\frac{x-\xi_i}{\varepsilon}\right)}\left[
1
+\sum_{j=1}^4\left(\frac{p-1}{p}\right)^j\frac{1}{\gamma^{jp}}
\big(
\omega^j_{\mu_i}-f^j_{\mu_i}
\big)\left(\frac{x-\xi_i}{\varepsilon}\right)
\right]\quad\ \ \,\nonumber\\
&\times
\left[
p\gamma^{p}+
\omega_{\mu_k}\left(\frac{x-\xi_k}{\varepsilon}\right)
+\sum_{j=1}^4\left(\frac{p-1}{p}\right)^j\frac{1}{\gamma^{jp}}\omega^j_{\mu_{k}}\left(\frac{x-\xi_k}{\varepsilon}\right)
+\,O\left(
|x-\xi_k|^\alpha
+\sum_{i=1}^m(\varepsilon\mu_i)^{\alpha/2}
\right)
\right]dx
+O\left(\frac{\varepsilon}{p^2\gamma^{2(p-1)}}\right)
\nonumber\\
=&\frac{1}{2p^2\gamma^{2(p-1)}}
\left\{
\frac{1}{\varepsilon^2}
\sum\limits_{k=1}^{m}\int_{\Omega\cap B_{1/|\log\varepsilon|^{2\kappa}}(\xi_k)}a(x)
e^{\omega_{\mu_k}\left(\frac{x-\xi_k}{\varepsilon}\right)}\left[
1
+\frac{p-1}{p}\frac{1}{\gamma^{p}}
\big(
\omega^1_{\mu_k}-f^1_{\mu_k}
\big)\left(\frac{x-\xi_k}{\varepsilon}\right)
\right]
\times
\left[
p\gamma^{p}+
\omega_{\mu_k}\left(\frac{x-\xi_k}{\varepsilon}\right)
\right]
dx
\right.
\nonumber\\
&\left.
+
O\left(\frac{\log^\beta|\log\varepsilon|}{|\log\varepsilon|}\right)
\right\}.
\endaligned
$$
On one hand, by using the relation $p\gamma^{p}=-4\log\varepsilon$
and the change of   variables $\varepsilon\mu_k z=x-\xi_k$ we obtain
$$
\aligned
&\frac{1}{\varepsilon^2}
\int_{\Omega\cap B_{1/|\log\varepsilon|^{2\kappa}}(\xi_k)}a(x)
e^{\omega_{\mu_k}\left(\frac{x-\xi_k}{\varepsilon}\right)}
\left[
p\gamma^{p}+
\omega_{\mu_k}\left(\frac{x-\xi_k}{\varepsilon}\right)
\right]
dx\\
&\quad=
\int_{\Omega\cap B_{1/|\log\varepsilon|^{2\kappa}}(\xi_k)}
a(x)\frac{8\varepsilon^2\mu^2_k}{(\varepsilon^2\mu_k^2+|x-\xi_k|^2)^2}
\log\frac{8\mu^2_k}{(\varepsilon^2\mu_k^2+|x-\xi_k|^2)^2}dx\\
&\quad=
\int_{\Omega_{\varepsilon\mu_k}\cap B_{1/(\varepsilon\mu_k|\log\varepsilon|^{2\kappa})}(0)}
a(\varepsilon\mu_k z+\xi_k)\frac{8}{(1+|z|^2)^2}
\left[\log\frac{8}{(1+|z|^2)^2}-
\log\left(\varepsilon^4\mu_k^2\right)
\right]dz,
\endaligned
$$
where $\Omega_{\varepsilon\mu_k}=\frac{1}{\varepsilon\mu_k}(\Omega-\{\xi_k\})$. But
$$
\aligned
\int_{\Omega_{\varepsilon\mu_k}\cap B_{1/(\varepsilon\mu_k|\log\varepsilon|^{2\kappa})}(0)}
a(\varepsilon\mu_k z+\xi_k)\frac{8}{(1+|z|^2)^2}
=c_k a(\xi_k)+O(\varepsilon\mu_k),
\endaligned
$$
and
$$
\aligned
\int_{\Omega_{\varepsilon\mu_k}\cap B_{1/(\varepsilon\mu_k|\log\varepsilon|^{2\kappa})}(0)}
a(\varepsilon\mu_k z+\xi_k)\frac{8}{(1+|z|^2)^2}
\log\frac{1}{(1+|z|^2)^2}=-2c_k a(\xi_k)+O(\varepsilon\mu_k).
\endaligned
$$
Then
\begin{equation}\label{6.3}
\aligned
\frac{1}{\varepsilon^2}
\int_{\Omega\cap B_{1/|\log\varepsilon|^{2\kappa}}(\xi_k)}a(x)
e^{\omega_{\mu_k}\left(\frac{x-\xi_k}{\varepsilon}\right)}
\left[
p\gamma^{p}+
\omega_{\mu_k}\left(\frac{x-\xi_k}{\varepsilon}\right)
\right]
dx
=c_ka(\xi_k)\big[\log8-\log\left(\varepsilon^4\mu_k^2\right)-2\big]+O(\varepsilon\mu_k).
\endaligned
\end{equation}
%Similarly,
%$$
%\aligned
%&\frac{1}{\varepsilon^2}
%\int_{\Omega\cap B_{1/|\log\varepsilon|^{2\kappa}}(\xi_k)}a(x)
%e^{\omega_{\mu_k}\left(\frac{x-\xi_k}{\varepsilon}\right)}
%\left[
%\frac{p-1}{p}\frac{1}{\gamma^{p}}
%\big(
%\omega^1_{\mu_k}-f^1_{\mu_k}
%\big)\left(\frac{x-\xi_k}{\varepsilon}\right)
%\right]
%\left[
%p\gamma^{p}+
%\omega_{\mu_k}\left(\frac{x-\xi_k}{\varepsilon}\right)
%\right]
%dx\\
%=&\,(p-1)a(\xi_k) \int_{\Omega_{\varepsilon\mu_k}\cap B_{1/(\varepsilon\mu_k|\log\varepsilon|^{2\kappa})}(0)}
%e^{\omega_{\mu_k}(\mu_k z)}
%\big(
%\omega^1_{\mu_k}-f^1_{\mu_k}
%\big)(\mu_k z)
%\,
%d (\mu_k z)
%+O\left(\frac{\log^\beta|\log\varepsilon|}{|\log\varepsilon|}\right).
%\endaligned
%$$
%$$
%\aligned
%&\frac{1}{\varepsilon^2}
%\int_{\Omega\cap B_{1/|\log\varepsilon|^{2\kappa}}(\xi_k)}a(x)
%e^{\omega_{\mu_k}\left(\frac{x-\xi_k}{\varepsilon}\right)}
%\left[
%\frac{p-1}{p}\frac{1}{\gamma^{p}}
%\big(
%\omega^1_{\mu_k}-f^1_{\mu_k}
%\big)\left(\frac{x-\xi_k}{\varepsilon}\right)
%\right]
%\left[
%p\gamma^{p}+
%\omega_{\mu_k}\left(\frac{x-\xi_k}{\varepsilon}\right)
%\right]
%dx\\
%=&\int_{\Omega_{\varepsilon\mu_k}\cap B_{1/(\varepsilon\mu_k|\log\varepsilon|^{2\kappa})}(0)}
%a(\varepsilon\mu_k z+\xi_k)\frac{8(p-1)}{(1+|z|^2)^2}
%\left\{
%\left[
%\frac12\big(\upsilon_{\infty}\big)^2
%-\omega^0_{\infty}
%\right](z)
%+(1-2\log\mu_k)\left(
%\frac{1-|z|^2}{|z|^2+1}\log8-\frac{2|z|^2}{|z|^2+1}
%\right)\right.\\
%&\left.
%-\frac{\,4(\log^2\mu_k-\log\mu_k)\,}{\,|z|^2+1\,}
%+2\log^2\mu_k-2\log\mu_k\right\}dz+O\left(\frac{\log^\beta|\log\varepsilon|}{|\log\varepsilon|}\right).
%\endaligned
%$$
On the other hand,
from (\ref{2.9}),  (\ref{2.15}) and the first definition of $D^1_{\mu_k}$ in (\ref{2.14})   we can compute
$$
\aligned
\int_{\Omega_{\varepsilon\mu_k}\cap B_{1/(\varepsilon\mu_k|\log\varepsilon|^{2\kappa})}(0)}
e^{\omega_{\mu_k}(\mu_k z)}
\big(
\omega^1_{\mu_k}-f^1_{\mu_k}
\big)(\mu_k z)
\,
d (\mu_k z)
=
&\,
\frac{c_k}{8\pi}
\int_{\mathbb{R}^2}
e^{\omega_{\mu_k}(\mu_k z)}
\big(
\omega^1_{\mu_k}-f^1_{\mu_k}
\big)(\mu_k z)
\,
d (\mu_k z)+O\left(
\varepsilon \mu_k^2
\right)\\[1mm]
=&\,
-\frac{c_k}{8\pi}
\int_{\mathbb{R}^2}
\Delta\big[\omega^1_{\mu_k}(\mu_kz)\big]dz
+O\left(
\varepsilon \mu_k^2
\right)
\\[1mm]
=&\,
-\frac{1}{4}c_k D^1_{\mu_k}+O\left(
\varepsilon \mu_k^2
\right).
\endaligned
$$
%$$
%\aligned
%\int_{\Omega_{\varepsilon\mu_k}\cap B_{1/(\varepsilon\mu_k|\log\varepsilon|^{2\kappa})}(0)}
%a(\varepsilon\mu_k z+\xi_k)\frac{8|z|^2}{(1+|z|^2)^3}
%=\frac12c_k a(\xi_k)+O(\varepsilon\mu_k).
%\endaligned
%$$
Furthermore,  by the explicit expression of $D^1_{\mu_k}$ in  (\ref{2.15}) we find
\begin{eqnarray}\label{6.4}
\frac{1}{\varepsilon^2}
\int_{\Omega\cap B_{1/|\log\varepsilon|^{2\kappa}}(\xi_k)}a(x)
e^{\omega_{\mu_k}\left(\frac{x-\xi_k}{\varepsilon}\right)}
\left[
\frac{p-1}{p}\frac{1}{\gamma^{p}}
\big(
\omega^1_{\mu_k}-f^1_{\mu_k}
\big)\left(\frac{x-\xi_k}{\varepsilon}\right)
\right]
\left[
p\gamma^{p}+
\omega_{\mu_k}\left(\frac{x-\xi_k}{\varepsilon}\right)
\right]
dx
\nonumber
\\[1mm]
=\,(p-1)a(\xi_k) \int_{\Omega_{\varepsilon\mu_k}\cap B_{1/(\varepsilon\mu_k|\log\varepsilon|^{2\kappa})}(0)}
e^{\omega_{\mu_k}(\mu_k z)}
\big(
\omega^1_{\mu_k}-f^1_{\mu_k}
\big)(\mu_k z)
\,
d (\mu_k z)
+O\left(\frac{\log^\beta|\log\varepsilon|}{|\log\varepsilon|}\right)
\qquad\quad
\nonumber
\\[1mm]
=\,(p-1)c_k a(\xi_k)(2\log\mu_k+2-\log8)
+O\left(\frac{\log^\beta|\log\varepsilon|}{|\log\varepsilon|}\right).
\qquad\qquad\qquad\qquad\quad
\qquad\qquad\qquad\qquad\ \ \,
\end{eqnarray}
%\begin{eqnarray}\label{6.4}
%\frac{1}{\varepsilon^2}
%\int_{\Omega\cap B_{1/|\log\varepsilon|^{2\kappa}}(\xi_k)}a(x)
%e^{\omega_{\mu_k}\left(\frac{x-\xi_k}{\varepsilon}\right)}
%\left[
%\frac{p-1}{p}\frac{1}{\gamma^{p}}
%\big(
%\omega^1_{\mu_k}-f^1_{\mu_k}
%\big)\left(\frac{x-\xi_k}{\varepsilon}\right)
%\right]
%\left[
%p\gamma^{p}+
%\omega_{\mu_k}\left(\frac{x-\xi_k}{\varepsilon}\right)
%\right]
%dx
%\nonumber
%\\
%=(
%2\log\mu_k+2
%-\log8
%)(p-1)c_k a(\xi_k)
%+O\left(\frac{\log^\beta|\log\varepsilon|}{|\log\varepsilon|}\right).
%\end{eqnarray}
%\begin{eqnarray}\label{6.4}
%\frac{1}{\varepsilon^2}
%\int_{\Omega\cap B_{1/|\log\varepsilon|^{2\kappa}}(\xi_k)}a(x)
%e^{\omega_{\mu_k}\left(\frac{x-\xi_k}{\varepsilon}\right)}
%\left[
%\frac{p-1}{p}\frac{1}{\gamma^{p}}
%\big(
%\omega^1_{\mu_k}-f^1_{\mu_k}
%\big)\left(\frac{x-\xi_k}{\varepsilon}\right)
%\right]
%\left[
%p\gamma^{p}+
%\omega_{\mu_k}\left(\frac{x-\xi_k}{\varepsilon}\right)
%\right]
%dx
%\nonumber
%\\
%=(p-1)c_ka(\xi_k)\left\{
%\frac{1}{8\pi}\int_{\mathbb{R}^2}
%\frac{8}{(1+|z|^2)^2}
%\left[
%\frac12\big(\upsilon_{\infty}\big)^2
%-\omega^0_{\infty}
%\right](z)dz
%-1+2\log\mu_k\right\}
%+O\left(\frac{\log^\beta|\log\varepsilon|}{|\log\varepsilon|}\right).
%\end{eqnarray}
Hence  by (\ref{6.3}) and (\ref{6.4}),
\begin{equation}\label{6.5}
\aligned
I_A=\frac1{2p^2\gamma^{2(p-1)}}
\sum\limits_{k=1}^{m}c_k a(\xi_k)\left
\{4|\log\varepsilon|+
(p-2)(2\log\mu_k+2-\log8)
+
O\left(\frac{\log^\beta|\log\varepsilon|}{|\log\varepsilon|}\right)
\right
\}.
\endaligned
\end{equation}
Regarding the expression $I_B$, by  (\ref{2.28}) we have
$$
\aligned
I_B
=\frac{\lambda\varepsilon^2}{p}\left[\sum\limits_{i=1}^m
\left(\int_{\Omega_\varepsilon\cap B_{\mu_i|\log\varepsilon|^\tau}(\xi'_i)}
+\int_{\Omega_\varepsilon\cap \big( B_{1/(\varepsilon|\log\varepsilon|^{2\kappa})}(\xi'_i)\setminus
B_{\mu_i|\log\varepsilon|^\tau}(\xi'_i)\big)}
\right)
+\int_{\Omega_\varepsilon\setminus\cup_{i=1}^mB_{1/(\varepsilon|\log\varepsilon|^{2\kappa})}(\xi'_i)}
\right]a(\varepsilon y)
e^{\gamma^p\left(
1+\frac{V_{\xi'}(y)}{p\gamma^p}
\right)^p}dy.
\endaligned
$$
By (\ref{1.5}), (\ref{2.7}) and (\ref{2.36}),
$$
\aligned
\frac{\lambda\varepsilon^2}{p}
\int_{\Omega_\varepsilon\setminus\cup_{i=1}^mB_{1/(\varepsilon|\log\varepsilon|^{2\kappa})}(\xi'_i)}
a(\varepsilon y)
e^{\gamma^p\left(
1+\frac{V_{\xi'}(y)}{p\gamma^p}
\right)^p}dy=
O(\lambda)\exp\left[
O\left(\frac{\log^p|\log\varepsilon|}{|\log\varepsilon|^{p-1}}\right)
\right].
\endaligned
$$
By (\ref{1.5}), (\ref{2.7}) and  (\ref{2.46}),
$$
\aligned
\frac{\lambda\varepsilon^2}{p}&
\int_{\Omega_\varepsilon\cap \big( B_{1/(\varepsilon|\log\varepsilon|^{2\kappa})}(\xi'_i)\setminus
B_{\mu_i|\log\varepsilon|^\tau}(\xi'_i)\big)}
a(\varepsilon y)
e^{\gamma^p\left(
1+\frac{V_{\xi'}(y)}{p\gamma^p}
\right)^p}dy\\
&
=\frac{1}{p^2\gamma^{2(p-1)}}
\int_{\Omega_\varepsilon\cap \big( B_{1/(\varepsilon|\log\varepsilon|^{2\kappa})}(\xi'_i)\setminus
B_{\mu_i|\log\varepsilon|^\tau}(\xi'_i)\big)}
a(\varepsilon y)
e^{\gamma^p\left[\left(
1+\frac{V_{\xi'}(y)}{p\gamma^p}
\right)^p-1\right]}dy\\
&
\leq\frac{D}{p^2\gamma^{2(p-1)}}
\int_{\Omega_\varepsilon\cap \big( B_{1/(\varepsilon|\log\varepsilon|^{2\kappa})}(\xi'_i)\setminus
B_{\mu_i|\log\varepsilon|^\tau}(\xi'_i)\big)}
e^{\left[
1+
O\left(\frac{\log|\log\varepsilon|}{|\log\varepsilon|}\right)
\right]
\omega_{\mu_i}\left(y-\xi'_i\right)}dy\\
&=\frac{1}{p^2\gamma^{2(p-1)}}O\left(\frac{1}{\,|\log\varepsilon|^{\tau}}\right).
\endaligned
$$
By  (\ref{2.6}), (\ref{2.13}) and  (\ref{2.381}),
$$
\aligned
\frac{\lambda\varepsilon^2}{p}&\int_{\Omega_\varepsilon\cap
B_{\mu_i|\log\varepsilon|^\tau}(\xi'_i)}
a(\varepsilon y)
e^{\gamma^p\left(
1+\frac{V_{\xi'}(y)}{p\gamma^p}
\right)^p}dy\\
&
=\frac{1}{p^2\gamma^{2(p-1)}}
\int_{\Omega_\varepsilon\cap B_{\mu_i|\log\varepsilon|^\tau}(\xi'_i)}
a(\varepsilon y)
e^{\gamma^p\left[\left(
1+\frac{V_{\xi'}(y)}{p\gamma^p}
\right)^p-1\right]}dy\\
&=
\frac{1}{p^2\gamma^{2(p-1)}}
\int_{\Omega_\varepsilon\cap B_{\mu_i|\log\varepsilon|^\tau}(\xi'_i)}
a(\varepsilon y)
e^{\omega_{\mu_i}\left(y-\xi'_i\right)}
\left[1
+
O\left(\frac{\log^\beta|\log\varepsilon|}{|\log\varepsilon|}\right)
\right]
dy\\
&=\frac{1}{p^2\gamma^{2(p-1)}}c_i a(\xi_i)\left[
1
+
O\left(\frac{\log^\beta|\log\varepsilon|}{|\log\varepsilon|}\right)
\right].
\endaligned
$$
Hence
\begin{equation}\label{6.10}
\aligned
I_B=\frac{1}{p^2\gamma^{2(p-1)}}
\sum_{i=1}^m
c_i a(\xi_i)\left[
1
+
O\left(\frac{\log^\beta|\log\varepsilon|}{|\log\varepsilon|}\right)
\right].
\endaligned
\end{equation}
Submitting (\ref{6.5})-(\ref{6.10}) into (\ref{6.2}), we obtain
$$
\aligned
J_\lambda\left(U_\xi\right)=
\frac1{2p^2\gamma^{2(p-1)}}
\sum\limits_{i=1}^{m}c_i a(\xi_i)\left
\{4|\log\varepsilon|+
(p-2)(2\log\mu_i+2-\log8)-2
+
O\left(\frac{\log^\beta|\log\varepsilon|}{|\log\varepsilon|}\right)
\right
\},
\endaligned
$$
which, together with the expansion of $\log(8\mu_i^2)$
in (\ref{2.24}), implies (\ref{6.1}).
\end{proof}

\section{Proofs of theorems}
\noindent {\bf Proof of Theorem 1.1.}
We will look for  a solution of problem (\ref{1.1})
in the form $u_\lambda=U_{\xi}+\widetilde{\phi}_{\xi}$, where
the concentration points $\xi=(\xi_1,\ldots,\xi_m)$
are determined by the parametrization
\begin{equation*}\label{7.1}
\aligned
\xi_i\equiv\xi_i(\mathbf{s},\mathbf{t})=s_i-\frac{t_i}{|\log\varepsilon|}\nu(s_i),
\quad i=1,\ldots,l,
\qquad\quad\textrm{but}\quad\qquad
\xi_i\equiv\xi_i(\mathbf{s},\mathbf{t})=s_i,
\quad i=l+1,\ldots,m,
\endaligned
\end{equation*}
where $\mathbf{s}=(s_1,\ldots,s_m)\in(\po)^m$ and $\mathbf{t}=(t_1,\ldots,t_l)\in(\mathbb{R}_{+})^l$
belong to the configuration space
$$
\aligned
\Lambda_d=\left\{
(\mathbf{s},\,\mathbf{t})\in(\po)^m\times(\mathbb{R}_{+})^l\,
\big|\,|s_i-s_k|>2d,
\quad\,d<t_{j}<1/d,
\quad\forall\,\,\,i, k=1,\ldots,m,\,\,\,j=1,\ldots,l,\,\,\,i\neq k
\right\},
\endaligned
$$
for any $d>0$ small and  independent of $\varepsilon>0$.
Notice that  if $(\mathbf{s},\mathbf{t})$
is a critical point of the reduced energy
$\widehat{F}_\lambda\big(\mathbf{s},\mathbf{t}\big):=F_\lambda\big(\xi(\mathbf{s},\mathbf{t})\big)$
in $\Lambda_d$, then the function
$u_\lambda=U_{\xi(\mathbf{s},\mathbf{t})}+\widetilde{\phi}_{\xi(\mathbf{s},\mathbf{t})}$
is a solution of problem (\ref{1.1}) with
the qualitative properties (\ref{1.19}) and (\ref{1.18}) described  by Theorem 1.1.
Hence with the aid of (\ref{5.5}), Propositions 5.2 and 6.1
we  are led to find a critical point of
the reduced energy $\widehat{F}_\lambda$, or
equivalently, a critical point of
\begin{eqnarray}\label{7.2}
\widetilde{F}_\varepsilon\big(\mathbf{s},\mathbf{t}\big)
:=\frac{\,2p^2\gamma^{2(p-1)}\,}{4|\log\varepsilon|}
\widehat{F}_\lambda\big(\mathbf{s},\mathbf{t}\big)
\qquad\qquad\qquad\qquad\qquad\qquad\qquad
\qquad\qquad\qquad\qquad\qquad\qquad
\qquad\qquad\qquad\quad\ \,
\nonumber\\
=
\sum_{i=1}^m
c_ia(\xi_i)
\left\{
1
-\frac{1}{4|\log\varepsilon|}\left[
c_i H_a(\xi_i,\xi_i)+\sum_{k=1,\,k\neq i}^mc_k
G_a(\xi_i,\xi_k)
-2\log8+4
\right]
\right\}
+
O\left(\frac{\log^\beta|\log\varepsilon|}{|\log\varepsilon|^2}\right).
\end{eqnarray}
We claim that $\widetilde{F}_\varepsilon\big(\mathbf{s},\mathbf{t}\big)$ can be written as
\begin{eqnarray}\label{7.3}
\widetilde{F}_\varepsilon\big(\mathbf{s},\mathbf{t}\big)=
8\pi\sum_{i=1}^l\left\{a(s_i)
+\frac{1}{|\log\varepsilon|}\big[
a(s_i)\log t_i-t_i\partial_{\nu}a(s_i)
\big]
\right\}+4\pi\sum_{k=l+1}^ma(s_k)
+\Gamma_\varepsilon(\mathbf{s})
+\frac{1}
{|\log\varepsilon|}\Theta_\varepsilon\big(\mathbf{s},\mathbf{t}\big),
\end{eqnarray}
where the smooth functions
$\Theta_\varepsilon\big(\mathbf{s},\mathbf{t}\big)$
depends on $\mathbf{s}$ and $\mathbf{t}$
but $\Gamma_\varepsilon(\mathbf{s})$ only depends on $\mathbf{s}$, and
$\Gamma_\varepsilon(\mathbf{s})$,
$\Theta_\varepsilon\big(\mathbf{s},\mathbf{t}\big)$,
$|\nabla\Gamma_\varepsilon(\mathbf{s})|$
and $|\nabla\Theta_\varepsilon\big(\mathbf{s},\mathbf{t}\big)|$
 uniformly converge to
zero as $\varepsilon\rightarrow0$. In fact,
 using asymptotical properties of the regular part of
the anisotropic Green's function in \cite{AP}, we have that for any
$i=1,\ldots,m$  and $k=1,\ldots,l$,
\begin{equation}\label{7.4}
\aligned
H_a(\xi_i,\xi_k)=\frac1{2\pi}\log\frac{1}{|\xi_i-\hat{\xi}_k|}+
\frac1{2\pi}\nabla\log a(\xi_k)\cdot V(\xi_i-\xi_k)-\frac1{2\pi}\nabla\log a(\hat{\xi}_k)\cdot V(\xi_i-\hat{\xi}_k)+\tilde{\mathrm{z}}(\xi_i,\xi_k),
\endaligned
\end{equation}
where
$\hat{\xi}_k=s_k+\frac{t_k}{|\log\varepsilon|}\nu(s_k)$,
the vector function
$V(\cdot)\in C^{\infty}(\mathbb{R}^2\setminus\{0\})\cap C^{\alpha}(\overline{B_r(0)})$
for any $r>0$ and $0<\alpha<1$,
the mapping
$y\in\Omega_d\mapsto \tilde{\mathrm{z}}(\cdot,y)\in
 C^1\big(\Omega_d,\,C^{1}(\overline{\Omega})\big)$
with $\Omega_d=\big\{
y\in\overline{\Omega}\big|\,
\dist(y,\partial\Omega)<d
\big\}$.
Then
\begin{equation}\label{7.5}
\aligned
H_a(\xi_i,\xi_i)=-\frac{1}{2\pi}\log\left(\frac{2t_i}{|\log\varepsilon|}\right)
+\tilde{\mathrm{z}}(s_i,s_i)+O\left(\frac{t_i^\alpha}{|\log\varepsilon|^\alpha}\right),
\quad\forall\,\,
i=1,\ldots,l.
\endaligned
\end{equation}
Moreover,  if $i, k=1,\ldots,l$ with $i\neq k$,
\begin{eqnarray}\label{7.6}
G_a(\xi_i,\xi_k)=-\frac{1}{2\pi}\log\left|s_i-s_k-\frac{t_i}{|\log\varepsilon|}\nu(s_i)+\frac{t_k}{|\log\varepsilon|}\nu(s_k)
\right|
-\frac{1}{2\pi}\log\left|s_i-s_k-\frac{t_i}{|\log\varepsilon|}\nu(s_i)-\frac{t_k}{|\log\varepsilon|}\nu(s_k)
\right|
\,\,\,\,\,
&&\nonumber\\[0.2mm]
+\frac{t_k}{\pi |\log\varepsilon|}\left\langle \nabla\log a(s_k),\,\nabla V(s_i-s_k)\cdot\nu(s_k)\right\rangle
-\frac{t_k}{\pi |\log\varepsilon|}\left\langle \big(\nabla\times(\nabla\log a)\big)(s_k)\cdot\nu(s_k),\,V(s_i-s_k)\right\rangle
&&\nonumber\\[0.1mm]
+\tilde{\mathrm{z}}(s_i,s_k)
-\frac{1}{|\log\varepsilon|}\left\langle\nabla_{(s_i,s_k)}\tilde{\mathrm{z}}(s_i,s_k),\,\big(t_i\nu(s_i),\,t_k\nu(s_k)\big)\right\rangle
+O\left(\frac{\,t_i^2+t_k^2\,}{|\log\varepsilon|^2}\right),
\qquad\qquad\qquad\qquad\quad\,\,\,
&&
\end{eqnarray}
while if $i=l+1,\ldots,m$  and $k=1,\ldots,l$,
\begin{eqnarray}\label{7.7}
G_a(\xi_i,\xi_k)=-\frac{1}{2\pi}\log\left|s_i-s_k+\frac{t_k}{|\log\varepsilon|}\nu(s_k)
\right|
-\frac{1}{2\pi}\log\left|s_i-s_k-\frac{t_k}{|\log\varepsilon|}\nu(s_k)
\right|
\qquad\qquad\qquad\qquad\quad\qquad\qquad\,\,
&&\nonumber\\[0.2mm]
+\frac{t_k}{\pi |\log\varepsilon|}\left\langle \nabla\log a(s_k),\,\nabla V(s_i-s_k)\cdot\nu(s_k)\right\rangle
-\frac{t_k}{\pi |\log\varepsilon|}\left\langle \big(\nabla\times(\nabla\log a)\big)(s_k)\cdot\nu(s_k),\,V(s_i-s_k)\right\rangle
&&\nonumber\\[0.1mm]
+\tilde{\mathrm{z}}(s_i,s_k)
-\frac{t_k}{|\log\varepsilon|}\left\langle\nabla_{s_k} \tilde{\mathrm{z}}(s_i,s_k),\,\nu(s_k)\right\rangle
+O\left(\frac{t_k^2}{\,|\log\varepsilon|^2\,}\right).
\qquad\qquad\qquad\qquad\qquad\qquad\qquad\qquad\,\,
&&
\end{eqnarray}
On the other hand,
using  the smooth property of $a(x)$ over $\oo$,  we  perform
a Taylor expansion around each boundary point $s_i$ along the inner normal vector
$-\nu(s_i)$ to give
\begin{equation}\label{7.8}
\aligned
a(\xi_i)=a(s_i)-\frac{t_i}{|\log\varepsilon|}\partial_{\nu}a(s_i)+O\left(\frac{t_i^2}{|\log\varepsilon|^2}\right),
\quad\forall\,\,
i=1,\ldots,l.
\endaligned
\end{equation}
Inserting  (\ref{7.5})-(\ref{7.8}) into (\ref{7.2}) and using
(\ref{2.18}) and the fact that
$a(\xi_i)G_a(\xi_i,\xi_k)=a(\xi_k)G_a(\xi_k,\xi_i)$ for all
$i,k=1,\ldots,m$ with $i\neq k$, we conclude that expansion (\ref{7.3}) holds.

We seek a critical point of  $\widetilde{F}_\varepsilon$ by degree theory.
Let $\partial_{T(s_i)}$ be the tangential derivative which is defined on $s_i\in\po$.
Set
$$
\aligned
A(s_i,t_i)=a(s_i)
+\frac{1}{|\log\varepsilon|}\big[
a(s_i)\log t_i-t_i\partial_{\nu}a(s_i)
\big],
\,\,\ \,i=1,\ldots,l.
\endaligned
$$
Then
$$
\aligned
\partial_{T(s_i)}A(s_i,t_i)=\left(1+
\frac{\log t_i}{|\log\varepsilon|}
\right)\partial_{T}a(s_i)
-\frac{t_i}{|\log\varepsilon|}
\partial_{T}\partial_{\nu}a(s_i)
\large,
\,\qquad\,
\partial_{t_i}A(s_i,t_i)=\frac{1}{|\log\varepsilon|}
\left(
\frac{a(s_i)}{t_i}-\partial_{\nu}a(s_i)
\right).
\endaligned
$$
Due to $\partial_{\nu}a(\xi^*_i)>0$ with
$i=1,\ldots,l$,
 we can choose $d$ small enough so that
for
any $s_i\in B_{d}(\xi^*_i)\cap\po$,
there exists a unique positive
$t_i=t_i(s_i)=\frac{a(s_i)}{\partial_{\nu}a(s_i)}$
such that $\partial_{t_i}A(s_i,t_i)=0$ and $\partial^2_{t_it_i}A(s_i,t_i)<0$.
Set $t_i^*=t_i(\xi_i^*)$, $i=1,\ldots,l$.
Since $\xi^*_1,\ldots,\xi^*_m$ are $m$ different strict local maximum or
strict local minimum points of $a(x)$ on $\po$, we have that that for
any sufficiently small $d$, $\varepsilon$ and any $i=1,\ldots,l$,  the
Brouwer degrees
$$
\aligned
\deg\big(\big(\partial_{T(s_i)}A,&\,\partial_{t_i}A\big),\,\big(B_{d}(\xi^*_i)\cap\po\big)\times\big(t_i^*-d,\,t_i^*+d\big),\,0\big)\\[1mm]
=&
\deg\big(\big(\partial_{T}a(s_i),\,\partial_{t_i}A\big),\,\big(B_{d}(\xi^*_i)\cap\po\big)\times\big(t_i^*-d,\,t_i^*+d\big),\,0\big)\\[0.5mm]
=&\sign\,\det\left[
\left(\aligned
    \partial^2_{TT}a(\xi^*_i) \quad &\quad  \frac{1}{|\log\varepsilon|}
\left(
\frac{\partial_{T}a(\xi^*_i)}{\,t^*_i\,}-\partial_{T}\partial_{\nu}a(\xi^*_i)
\right) \\
    0 \quad \quad \quad &\quad \quad \   -\frac{1}{\,|\log\varepsilon|\,}\frac{1}{(t^*_i)^2}a(\xi^*_i) \\
  \endaligned
\right)
\right]=\pm1
\neq0,
\endaligned
$$
and
$$
\aligned
\deg\left(\big(\partial_{T(s_{l+1})}a,\ldots,\partial_{T(s_m)}a\big),
\,\prod_{k=l+1}^m\big(B_{d}(\xi^*_k)\cap\po\big),\,0\right)
=\sign\left(\prod_{k=l+1}^m\partial^2_{TT}a(\xi^*_k)\right)=\pm1
\neq0.
\endaligned
$$
Then by (\ref{7.3}),
$$
\aligned
&\,\deg\left(\nabla_{\left(T(\mathbf{s}),\mathbf{t}\right)}\widetilde{F}_\varepsilon\big(\mathbf{s},\mathbf{t}\big),\,\,\prod_{i=1}^l\big( B_{d}(\xi^*_i)\cap\po\big)\times\prod_{k=l+1}^m\big(B_{d}(\xi^*_k)
\cap\po\big)\times\prod_{i=1}^l\big(t_{i}^*-d,\,t_{i}^*+d\big),\,\,0\right)\\[1mm]
=&
\prod_{i=1}^l\deg\large\big(\big(\partial_{T(s_i)}A,\partial_{t_i}A\big),\big(B_{d}(\xi^*_i)\cap\po\big)
\times\big(t_i^*-d,t_i^*+d\big),0\large\big)\times
\deg\left(\big(\partial_{T(s_{l+1})}a,\ldots,\partial_{T(s_m)}a\big),
\prod_{k=l+1}^m\big(B_{d}(\xi^*_k)\cap\po\big),0\right)\\[1mm]
\neq&\,0.
\endaligned
$$
Hence if  $\varepsilon$ is small enough, there exists
$(\mathbf{s}^\varepsilon,\mathbf{t}^\varepsilon)$ such that
$\nabla_{\left(T(\mathbf{s}),\mathbf{t}\right)}\widetilde{F}_\varepsilon\big(\mathbf{s}^\varepsilon,\mathbf{t}^\varepsilon\big)=0$.
In particular, $\mathbf{s}^\varepsilon=(s^\varepsilon_1,\ldots,s^\varepsilon_m)\rightarrow(\xi^*_1,\ldots,\xi^*_m)$ as $\varepsilon\rightarrow0$.

\vspace{1mm}

\noindent {\bf Proof of (\ref{1.14}).}\,
Let
$\xi^\varepsilon=\xi(\mathbf{s}^\varepsilon,\mathbf{t}^\varepsilon)
=\big(\xi_1(\mathbf{s}^\varepsilon,\mathbf{t}^\varepsilon),\ldots,\xi_m(\mathbf{s}^\varepsilon,\mathbf{t}^\varepsilon)\big)$
and
$\upsilon_\lambda(x)=p\gamma^{p-1}u_\lambda(x)-p\gamma^{p}$.
Then
$\upsilon_\lambda(x)=\big(V_{(\xi^\varepsilon)'}+\phi_{(\xi^\varepsilon)'}\big)
(\frac{x}{\varepsilon} )$
and
\begin{equation}\label{6.9}
\aligned
\lambda u_\lambda^{p-1}e^{u_\lambda^p}
=
\lambda\gamma^{p-1}e^{\gamma^p}
\left(1+\frac{\upsilon_\lambda}{p\gamma^p}\right)^{p-1}
e^{\gamma^p\left[\left(1+\frac{\upsilon_\lambda}{p\gamma^p}\right)^p-1\right]}
=\lambda\gamma^{p-1}e^{\gamma^p}f\big(V_{(\xi^\varepsilon)'}+\phi_{(\xi^\varepsilon)'}\big)
\Big(\frac{x}{\varepsilon}\Big)
\endaligned
\end{equation}
where the function $f(\cdot)$ is given by (\ref{2.27}).
For any $\psi\in C_{c}(\oo)$,
by (\ref{1.5}) and (\ref{6.9}) we get
$$
\aligned
p\gamma^{p-1}\int_{\Omega}
\lambda u_\lambda^{p-1}e^{u_\lambda^p}
\psi
dx
=&
\left[
\,
\sum\limits_{i=1}^m
\left(\int_{\Omega_\varepsilon\cap B_{\mu_i|\log\varepsilon|^\tau}\left((\xi_i^\varepsilon)'\right)}
+\int_{
\Omega_\varepsilon\cap \big( B_{1/(\varepsilon|\log\varepsilon|^{2\kappa})}\left((\xi_i^\varepsilon)'\right)\setminus
B_{\mu_i|\log\varepsilon|^\tau}\left((\xi_i^\varepsilon)'\right)\big)
}
\right)
\right.
\\[1.5mm]
&
\left.
+
\int_{\Omega_\varepsilon\setminus\cup_{i=1}^mB_{1/(\varepsilon|\log\varepsilon|^{2\kappa})}\big((\xi_i^\varepsilon)'\big)}
\right]
f\big(V_{(\xi^\varepsilon)'}+\phi_{(\xi^\varepsilon)'}\big)
(y)\,\psi(\varepsilon y)dy.
\endaligned
$$
Observe  that $\phi_{(\xi^\varepsilon)'}$
is a  higher-order term in $\upsilon_{\lambda}$
because   $\|\phi_{(\xi^\varepsilon)'}\|_{L^{\infty}(\Omega_\varepsilon)}=O\left(1/|\log\varepsilon|^3\right)$.
Then by  (\ref{2.39}),
$$
\aligned
&\,
\int_{\Omega_\varepsilon\cap B_{\mu_i|\log\varepsilon|^\tau}\left((\xi_i^\varepsilon)'\right)}
f\big(V_{(\xi^\varepsilon)'}+\phi_{(\xi^\varepsilon)'}\big)
(y)\,\psi(\varepsilon y)dy\\
=&\,
\int_{\Omega_\varepsilon\cap B_{\mu_i|\log\varepsilon|^\tau}\left((\xi_i^\varepsilon)'\right)}
e^{\omega_{\mu_i}\left(y-(\xi_i^\varepsilon)'\right)}\psi(\varepsilon y)
\left[1+
O
\left(\frac{\log^\beta|\log\varepsilon|}{|\log\varepsilon|}\right)
\right]
(y)dy\\[2mm]
=&\,
c_i
\psi(\xi_i^\varepsilon)
+
o
\left(1\right).
\endaligned
$$
By  (\ref{2.46})-(\ref{2.43}),
$$
\aligned
&\,
\left|
\int_{
\Omega_\varepsilon\cap \big( B_{1/(\varepsilon|\log\varepsilon|^{2\kappa})}\left((\xi_i^\varepsilon)'\right)\setminus
B_{\mu_i|\log\varepsilon|^\tau}\left((\xi_i^\varepsilon)'\right)\big)
}
f\big(V_{(\xi^\varepsilon)'}+\phi_{(\xi^\varepsilon)'}\big)
(y)\,
\psi(\varepsilon y)dy
\right|\\
\leq
&\,D^2
\int_{
\Omega_\varepsilon\cap \big( B_{1/(\varepsilon|\log\varepsilon|^{2\kappa})}\left((\xi_i^\varepsilon)'\right)\setminus
B_{\mu_i|\log\varepsilon|^\tau}\left((\xi_i^\varepsilon)'\right)\big)
}
|\psi(\varepsilon y)|\,
e^{\left[
1-
\frac{1}{4}
\sum_{j=1}^4\left(\frac{p-1}{p}\right)^j\frac{D^j_{\mu_i}}{\gamma^{jp}}
\right]
\omega_{\mu_i}\left(y-\xi'_i\right)
}
\left(1+\frac{\log^{p-1}|\log\varepsilon|}{|\log\varepsilon|^{p-1}}\right)
dy\\[1mm]
=&\,
o
\left(1\right).
\endaligned
$$
By  (\ref{2.36})-(\ref{2.37}),
$$
\aligned
&\,
\left|
\int_{\Omega_\varepsilon\setminus\cup_{i=1}^mB_{1/(\varepsilon|\log\varepsilon|^{2\kappa})}\big((\xi_i^\varepsilon)'\big)}
f\big(V_{(\xi^\varepsilon)'}+\phi_{(\xi^\varepsilon)'}\big)
(y)\,
\psi(\varepsilon y)dy
\right|\\
\leq
&\,C
\int_{\Omega_\varepsilon\setminus\cup_{i=1}^mB_{1/(\varepsilon|\log\varepsilon|^{2\kappa})}\big((\xi_i^\varepsilon)'\big)}
|\psi(\varepsilon y)|
\frac{
\varepsilon^{\frac4p}\log^{p-1}|
\log\varepsilon|
}{|\log\varepsilon|^{p-1}}
\exp\left[O\left(\frac{\log^{p}|
\log\varepsilon|}{|\log\varepsilon|^{p-1}}\right)\right]
dy\\[1mm]
=&\,
o
\left(1\right).
\endaligned
$$
Therefore,
$$
\aligned
p\gamma^{p-1}\int_{\Omega}
\lambda u_\lambda^{p-1}e^{u_\lambda^p}
\psi
dx
=
\sum\limits_{i=1}^mc_i
\psi(\xi_i^\varepsilon)
+
o
\left(1\right)
\rightarrow
\sum\limits_{i=1}^mc_i
\psi(\xi_i^*)
\quad\quad
\textrm{as}
\,\,\,\,
\lambda\rightarrow0.
\endaligned
$$

\vspace{1mm}

\noindent {\bf Proof of (\ref{1.15})-(\ref{1.17}).}\,
Since $\phi_{(\xi^\varepsilon)'}$
is a  correction of third order for $V_{(\xi^\varepsilon)'}$,
by (\ref{2.381}) and  (\ref{6.10}) we can compute
$$
\aligned
\frac{\lambda p}{2}\int_{\Omega}
a(x)
\big(e^{u_\lambda^p}-1\big)dx=&\frac{1}{2\gamma^{2(p-1)}}
\left\{\sum_{i=1}^m
a(\xi_i^\varepsilon)
\int_{\Omega_\varepsilon \cap B_{\mu_i|\log\varepsilon|^\tau}((\xi_i^\varepsilon)')}
e^{\omega_{\mu_i}\left(y-(\xi_i^\varepsilon)'\right)}
\left[1+\frac{p-1}{p}\frac1{\gamma^p}
B_1
\right.
\right.
\\[2mm]
&
+\left.\left.
\left(\frac{p-1}{p}\right)^2\frac1{\gamma^{2p}}
\left(B_2+\frac12(B_1)^2\right)
\right]
dy
+O\left(\frac{\log^\beta|\log\varepsilon|}{|\log\varepsilon|^3}\right)
\right\}.
\endaligned
$$
Taking into account  $\xi^\varepsilon=(\xi^\varepsilon_1,\ldots,\xi^\varepsilon_{m})\in\Omega^l\times(\po)^{m-l}$,
by (\ref{2.2})-(\ref{2.3}) we give
$$
\aligned
&\left(\frac{\lambda p}{2}\int_{\Omega}
a(x)
\big(e^{u_\lambda^p}-1\big)dx
\right)^{\frac{2-p}{p}}
=
\left(
\frac{1}{2\gamma^{2(p-1)}}
\sum_{k=1}^m
c_ka(\xi_k^\varepsilon)
\right)^{\frac{2-p}{p}}\left[1+
O\left(\frac{\log^\beta|\log\varepsilon|}{|\log\varepsilon|^3}\right)
\right]
\\
&
\qquad\qquad
\times
\left\{
\sum_{i=1}^m
\underbrace{
\frac{l_i^\varepsilon}{8\pi}
\int_{\mathbb{R}^2}
e^{\omega_{\mu_i}\left(y-(\xi_i^\varepsilon)'\right)}
\left[1+\frac{p-1}{p}\frac1{\gamma^p}
B_1
+\left(\frac{p-1}{p}\right)^2\frac1{\gamma^{2p}}
\left(B_2+\frac12(B_1)^2\right)
\right]
dy}\limits_{a_i}\right\}^{\frac{2-p}{p}},
\endaligned
$$
where
$$
\aligned
l_i^\varepsilon
=
c_ia(\xi_i^\varepsilon)
\left/
\left(\sum\limits_{k=1}^m
c_ka(\xi_k^\varepsilon)\right),
\qquad\quad
i=1,\ldots,m.
\right.
\endaligned
$$
Similarly, by (\ref{2.380})-(\ref{2.381}) we give
$$
\aligned
&\qquad\qquad\qquad\qquad
\left(\frac{\lambda p}{2}\int_{\Omega}a(x)u_\lambda^p e^{u_\lambda^p}dx
\right)^{\frac{2(p-1)}{p}}
=\left(\frac{1}{2\gamma^{p-2}}
\sum_{k=1}^m
c_ka(\xi_k^\varepsilon)
\right)^{\frac{2(p-1)}{p}}
\left[1+
O\left(\frac{\log^\beta|\log\varepsilon|}{|\log\varepsilon|^3}\right)
\right]
\\
&
\times\left\{\sum_{i=1}^m
\underbrace{
\frac{l_i^\varepsilon}{8\pi}
 \int_{\mathbb{R}^2}
e^{\omega_{\mu_i}\left(y-(\xi_i^\varepsilon)'\right)}
\left[1
+\frac{1}{p\gamma^p}
\Big(\,pA_1+(p-1)B_1\Big)
+\frac{p-1}{p\gamma^{2p}}
\big(A_1B_1+B_1\big)
+
\left(\frac{p-1}{p}\right)^2\frac1{\gamma^{2p}}
\left(B_2+\frac12(B_1)^2\right)
\right]
dy}\limits_{b_i}
\right\}^{\frac{2(p-1)}{p}}.
\endaligned
$$
Owing to the vanishing  identity of first order (see  Corollary A.8 of \cite{Z1})
%$$
%\aligned
%\frac{\,2-p\,}{p}
%\int_{\mathbb{R}^2}
%e^{\omega_{\mu_i}\left(y-\xi'_i\right)}
%B_1
%dy
%+
%\frac{2}{\,p\,}\int_{\mathbb{R}^2}
%e^{\omega_{\mu_i}\left(y-\xi'_i\right)}
%\big[\,p A_1+(p-1)B_1\big]
%dy
%\,\equiv\,
%0,
%\endaligned
%$$
\begin{eqnarray}\label{8.18}
\frac{\,2-p\,}{p}
\int_{\mathbb{R}^2}
e^{\omega_{\mu_i}\left(y-(\xi_i^\varepsilon)'\right)}
B_1
dy
+
\frac{2}{\,p\,}\int_{\mathbb{R}^2}
e^{\omega_{\mu_i}\left(y-(\xi_i^\varepsilon)'\right)}
\big[\,p A_1+(p-1)B_1\big]
dy
\,\equiv\,
0,
\end{eqnarray}
by the Taylor expansion we can compute
\begin{eqnarray}\label{6.6}
\beta_\lambda=
\left(\frac{\lambda p}{2}
\int_{\Omega}
a(x)\big(e^{u_\lambda^p}-1\big)dx
\right)^{\frac{2-p}{p}}
\left(\frac{\lambda p}{2}
\int_{\Omega}
a(x)u_\lambda^p
e^{u_\lambda^p}dx
\right)^{\frac{2(p-1)}{p}}
\qquad\qquad\qquad\qquad\qquad\qquad\qquad\qquad\qquad\qquad\qquad\, \,\,\,
&&
\nonumber
\\[1mm]
=
\frac{1}{2}
\left(\sum\limits_{k=1}^m
c_k a(\xi_k^\varepsilon)\right)
\left\{1+
\frac{(p-1)^2}{8\pi  p^2\gamma^{2p}}\sum_{i=1}^m
l_i^\varepsilon \left[\int_{\mathbb{R}^2}
e^{\omega_{\mu_i}\left(y-(\xi_i^\varepsilon)'\right)}
\Big(B_2+\frac12(B_1)^2\Big)
dy
+2
\int_{\mathbb{R}^2}
e^{\omega_{\mu_i}\left(y-(\xi_i^\varepsilon)'\right)}
(A_1B_1+B_1)
dy
\right]
\right.
&&
\nonumber
\\[1mm]
+\left.
\frac{(p-1)(p-2)}{(16\pi )^2 p^2\gamma^{2p}}
\left(\sum_{i=1}^m
l_i^\varepsilon
\int_{\mathbb{R}^2}
e^{\omega_{\mu_i}\left(y-(\xi_i^\varepsilon)'\right)}
B_1
dy
\right)^2\right\}
\left[1+
O\left(\frac{\log^\beta|\log\varepsilon|}{|\log\varepsilon|^3}\right)
\right].
\qquad\qquad\qquad\qquad\qquad\qquad\quad\quad\
\,\,\,
&&
\end{eqnarray}
Since
$$
\aligned
\int_{\mathbb{R}^2}
e^{\omega_{\mu_i}\left(y-(\xi_i^\varepsilon)'\right)}
\Big(B_2+\frac12(B_1)^2\Big)
dy
=\frac{1}{p-1}O\left(
\log^2\mu_i
\right),
\endaligned
$$
and
$$
\aligned
\int_{\mathbb{R}^2}
e^{\omega_{\mu_i}\left(y-(\xi_i^\varepsilon)'\right)}
(A_1B_1+B_1)
dy
=O\left(
\log^2\mu_i
\right),
\qquad\qquad\qquad
\int_{\mathbb{R}^2}
e^{\omega_{\mu_i}\left(y-(\xi_i^\varepsilon)'\right)}
B_1
dy
=O\left(
\log\mu_i
\right)
\endaligned
$$
(see Lemmas A.4, A.7 and B.5 of \cite{Z1}), we find
\begin{equation*}\label{6.16}
\aligned
\beta_\lambda=
\frac{1}{2}
\left(\sum\limits_{k=1}^m
c_k a(\xi_k^\varepsilon)\right)\left[
1
+
O\left(\frac{\log^2\mu_i}{|\log\varepsilon|^2}\right)
\right]\rightarrow \frac{1}{2}
\sum\limits_{k=1}^m
c_k a(\xi_k^*).
\endaligned
\end{equation*}
By  H\"{o}lder's inequality we observe
$$
\aligned
\left(\sum_{i=1}^m l_i^\varepsilon \int_{\mathbb{R}^2}
e^{\omega_{\mu_i}\left(y-(\xi_i^\varepsilon)'\right)}
B_1
dy
\right)^2&\leq
\left(
\sum_{i=1}^m l_i^\varepsilon
\right)\times
\sum_{i=1}^m l_i^\varepsilon \left(\int_{\mathbb{R}^2}
e^{\omega_{\mu_i}\left(y-(\xi_i^\varepsilon)'\right)}
B_1
dy
\right)^2
=\sum_{i=1}^m l_i^\varepsilon \left(\int_{\mathbb{R}^2}
e^{\omega_{\mu_i}\left(y-(\xi_i^\varepsilon)'\right)}
B_1
dy
\right)^2.
\endaligned
$$
%$$
%\aligned
%\left(\sum_{i=1}^m l_i^\varepsilon \int_{\mathbb{R}^2}
%e^{\omega_{\mu_i}\left(y-(\xi_i^\varepsilon)'\right)}
%B_1
%dy
%\right)^2&\leq
%\left(
%\sum_{i=1}^m l_i^\varepsilon
%\right)\times
%\sum_{i=1}^m l_i^\varepsilon \left(\int_{\mathbb{R}^2}
%e^{\omega_{\mu_i}\left(y-(\xi_i^\varepsilon)'\right)}
%B_1
%dy
%\right)^2\\[1mm]
%&=\sum_{i=1}^m l_i^\varepsilon \left(\int_{\mathbb{R}^2}
%e^{\omega_{\mu_i}\left(y-(\xi_i^\varepsilon)'\right)}
%B_1
%dy
%\right)^2,
%\endaligned
%$$
From (\ref{6.6})  we can derive that if  $0<p<1$,
\begin{eqnarray*}\label{6.14}
\beta_\lambda\leq
\frac{1}{2}
\left(\sum\limits_{k=1}^m
c_k a(\xi_k^\varepsilon)\right)
\left\{1+
\frac{p-1}{\,p^2\gamma^{2p}\,}
\sum_{i=1}^m
l_i^\varepsilon
\left[\frac{p-1}{8\pi}
\left(\int_{\mathbb{R}^2}
e^{\omega_{\mu_i}\left(y-(\xi_i^\varepsilon)'\right)}
\Big(B_2+\frac12(B_1)^2\Big)
dy
\right.
\right.
\right.
\qquad\qquad\qquad\qquad\quad\,
&&
\nonumber
\\[1mm]
\left.\left.\left.
+2
\int_{\mathbb{R}^2}
e^{\omega_{\mu_i}\left(y-(\xi_i^\varepsilon)'\right)}
(A_1B_1+B_1)
dy
\right)
+
\frac{p-2}{\,(16\pi)^2\,}
\left(
\int_{\mathbb{R}^2}
e^{\omega_{\mu_i}\left(y-(\xi_i^\varepsilon)'\right)}
B_1
dy
\right)^2\right]
\right\}
\left[1+
O\left(\frac{\log^\beta|\log\varepsilon|}{|\log\varepsilon|^3}\right)
\right].
\end{eqnarray*}
Owing to  the identity  of second order (see  Corollary B.6 of \cite{Z1})
%$$
%\aligned
%\frac{p-1}{8\pi}
%\left[\int_{\mathbb{R}^2}
%e^{\omega_{\mu_i}\left(y-\xi'_i\right)}
%\Big(B_2+\frac12(B_1)^2\Big)
%dy
%+2
%\int_{\mathbb{R}^2}
%e^{\omega_{\mu_i}\left(y-\xi'_i\right)}
%(A_1B_1+B_1)
%dy
%\right]
%+\frac{p-2}{\,(16\pi)^2\,}
%\left(
%\int_{\mathbb{R}^2}
%e^{\omega_{\mu_i}\left(y-\xi'_i\right)}
%B_1
%dy
%\right)^2\,\equiv\,4.
%\endaligned
%$$
\begin{eqnarray}\label{9.10}
\frac{p-1}{8\pi}
\Big[\int_{\mathbb{R}^2}
e^{\omega_{\mu_i}\left(y-(\xi_i^\varepsilon)'\right)}
\Big(B_2+\frac12(B_1)^2\Big)
dy
+2
\int_{\mathbb{R}^2}
e^{\omega_{\mu_i}\left(y-(\xi_i^\varepsilon)'\right)}
(A_1B_1+B_1)
dy
\Big]
+\frac{p-2}{(16\pi)^2}
\Big(
\int_{\mathbb{R}^2}
e^{\omega_{\mu_i}\left(y-(\xi_i^\varepsilon)'\right)}
B_1
dy
\Big)^2\equiv4,
\end{eqnarray}
we have that for $0<p<1$,
\begin{eqnarray*}\label{6.14}
\beta_\lambda\leq
\frac{1}{2}
\left(\sum\limits_{k=1}^m
c_k a(\xi_k^\varepsilon)\right)
\left\{1+
\frac{4(p-1)}{\,p^2\gamma^{2p}\,}
\left[1+
O\left(\frac{\log^\beta|\log\varepsilon|}{|\log\varepsilon|}\right)
\right]\right\}
<\,\frac{1}{2}
\sum\limits_{k=1}^m
c_k a(\xi_k^\varepsilon).
\end{eqnarray*}
But for $1<p<2$,  by  H\"{o}lder's inequality for vectors in $\mathbb{R}_{+}^m$ we  get
\begin{equation*}
\aligned
\beta_\lambda
&=
\frac{1}{2}
\left(\sum\limits_{k=1}^m
c_k a(\xi_k^\varepsilon)\right)
\left(
\sum_{i=1}^{m}a_i
\right)^{\frac{2-p}{p}}
\left(
\sum_{i=1}^{m}b_i
\right)^{\frac{2(p-1)}{p}}
\left[1+
O\left(\frac{\log^\beta|\log\varepsilon|}{|\log\varepsilon|^3}\right)
\right]\\[2mm]
&\geq
\frac{1}{2}
\left(\sum\limits_{k=1}^m
c_k a(\xi_k^\varepsilon)\right)
\left(
\sum_{i=1}^m
a_i^{\frac{2-p}{p}}b_i^{\frac{2(p-1)}{p}}
\right)
\left[1+
O\left(\frac{\log^\beta|\log\varepsilon|}{|\log\varepsilon|^3}\right)
\right].
\endaligned
\end{equation*}
Applying the Taylor expansion and using the identities (\ref{8.18}) and (\ref{9.10}) again, we
can compute
\begin{eqnarray*}
\left(
\frac{a_i}{l_i^\varepsilon}
\right)^{\frac{2-p}{p}}
\left(
\frac{b_i}{l_i^\varepsilon}
\right)^{\frac{2(p-1)}{p}}
=
1+
\frac{(p-1)^2}{ p^2\gamma^{2p}}\frac{1}{8\pi}
\left[\int_{\mathbb{R}^2}
e^{\omega_{\mu_i}\left(y-(\xi_i^\varepsilon)'\right)}
\Big(B_2+\frac12(B_1)^2\Big)
dy
+2
\int_{\mathbb{R}^2}
e^{\omega_{\mu_i}\left(y-(\xi_i^\varepsilon)'\right)}
(A_1B_1+B_1)
dy
\right]
&&
\nonumber
\\[1mm]
+\,
\frac{(p-1)(p-2)}{ p^2\gamma^{2p}}\frac{1}{(16\pi)^2}
\left(\int_{\mathbb{R}^2}
e^{\omega_{\mu_i}\left(y-(\xi_i^\varepsilon)'\right)}
B_1
dy
\right)^2+
O\left(\frac{\log^\beta|\log\varepsilon|}{|\log\varepsilon|^3}\right)
\qquad\qquad\qquad\quad\quad\,
&&\nonumber
\\[1mm]
=1+
\frac{\,4(p-1)\,}{ p^2\gamma^{2p}}
+
O\left(\frac{\log^\beta|\log\varepsilon|}{|\log\varepsilon|^3}\right).
\qquad\qquad\qquad\qquad\qquad\qquad\qquad
\qquad\qquad\quad\qquad\qquad\,
&&
\end{eqnarray*}
Hence for $1<p<2$,
\begin{eqnarray*}
\qquad\qquad\qquad\qquad\quad\quad\quad
\beta_\lambda\geq
\frac{1}{2}
\left(\sum\limits_{k=1}^m
c_k a(\xi_k^\varepsilon)\right)
\left\{1+
\frac{4(p-1)}{\,p^2\gamma^{2p}\,}
\left[1+
O\left(\frac{\log^\beta|\log\varepsilon|}{|\log\varepsilon|}\right)
\right]\right\}
>\frac{1}{2}
\sum\limits_{k=1}^m
c_k a(\xi_k^\varepsilon).
\qquad\qquad\qquad
\quad\,\,\,\,\,
\square
\end{eqnarray*}

\vspace{1mm}
\vspace{1mm}
\vspace{1mm}
\vspace{1mm}

\noindent {\bf Proof of Theorem 1.2.}
Let us just  find
a critical point $\xi^\varepsilon=(\xi^\varepsilon_1,\ldots,\xi^\varepsilon_m)\in\Omega^l\times(\po)^{m-l}$
of  $F_\lambda$ such  that points $\xi^\varepsilon_1,\ldots,\xi^\varepsilon_m$
accumulate to $\xi_*$.
For this aim, we consider the configuration space
$$
\aligned
\mathcal{O}^*_{d,\varepsilon}:=\left\{\xi=(\xi_1,\ldots,\xi_m)\in\big(B_d(\xi_*)\cap\Omega\big)^l\times\big(B_d(\xi_*)\cap\po\big)^{m-l}
\left|
\min\limits_{i\neq j}\big|\xi_i-\xi_j\big|>\frac1{|\log\varepsilon|^{\kappa}},
\quad
\min\limits_{1\leq k\leq l}\dist(\xi_k,\po)>\frac1{|\log\varepsilon|^{\kappa}}
\right.\right\},
\endaligned
$$
where $d>0$ is a sufficiently small but fixed number, independent of $\varepsilon$.
Using   (\ref{5.5}), Propositions 5.2 and 6.1  together with the fact that
$a(\xi_i)G_a(\xi_i,\xi_k)=a(\xi_k)G_a(\xi_k,\xi_i)$ for all
$i,k=1,\ldots,m$ with $i\neq k$,
we obtain that $F_\lambda$  reduces to
\begin{eqnarray}\label{7.9}
F_\lambda(\xi)=\frac{8\pi}{p^2\gamma^{2(p-1)}}
\left\{2\sum_{i=1}^la(\xi_i)\left[
|\log\varepsilon|-2\pi H_a(\xi_i,\xi_i)
-2\pi \sum_{k=1,\,k\neq i}^l
G_a(\xi_i,\xi_k)
\right]-4\pi\sum_{i=1}^l\sum_{k=l+1}^m a(\xi_k)
G_a(\xi_k,\xi_i)
\right.
\nonumber\\
\left.+\sum_{i=l+1}^ma(\xi_i)\left[
|\log\varepsilon|+\sum_{k=l+1,\,k\neq i}^m
\log|\xi_i-\xi_k|
\right]
+O\left(1\right)
\right\}
\qquad\qquad\qquad\qquad
\qquad\qquad\qquad\quad\,\,\,
\qquad\qquad
\end{eqnarray}
$C^0$-uniformly in
$\mathcal{O}^*_{d,\varepsilon}$.
Let us claim that
for any  $m\geq1$, $0\leq l\leq m$ and for any $\varepsilon$ small enough, the maximization problem
$$
\aligned
\max\limits_{(\xi_1,\ldots,\xi_m)\in\overline{\mathcal{O}}^*_{d,\varepsilon}}
F_\lambda(\xi_1,\ldots,\xi_m)
\endaligned
$$
has a solution in the interior of $\mathcal{O}^*_{d,\varepsilon}$.
If this claim  is proven, we can easily give  all the qualitative properties  (\ref{1.20})-(\ref{1.25}) of
clustered bubbling solutions of (\ref{1.1})
through  analogous arguments for those qualitative properties described by Theorem 1.1.

Let
$\xi^\varepsilon=(\xi^\varepsilon_1,\ldots,\xi^\varepsilon_m)\in
\overline{\mathcal{O}}^*_{d,\varepsilon}$ be the maximizer of $F_\lambda$. We are led to prove that $\xi^\varepsilon$
belongs to the interior of $\mathcal{O}^*_{d,\varepsilon}$.
First, we obtain a lower bound for $F_\lambda$ over $\overline{\mathcal{O}}^*_{d,\varepsilon}$.
Around the point $\xi_*\in\po$, we consider a smooth change of variables
$$
\aligned
H_{\xi_*}^\varepsilon(y)=
\varepsilon^{-2}H_{\xi_*}(\varepsilon^2 y),
\endaligned
$$
where $H_{\xi_*}:B_d(\xi_*)\mapsto\mathcal{M}$
is a  diffeomorphism and
$\mathcal{M}$  is an open neighborhood of the origin such that
$H_{\xi_*}(B_d(\xi_*)\cap\Omega)=\mathcal{M}\cap\mathbb{R}_+^2$
and
$H_{\xi_*}(B_d(\xi_*)\cap\partial\Omega)=\mathcal{M}\cap\partial\mathbb{R}_+^2$.
Let
$$
\aligned
\xi^0_i=\xi_*-\frac{t_i}{\sqrt{|\log\varepsilon|}}\nu(\xi_*),
\quad i=1,\ldots,l,
\qquad\quad\textrm{but}\quad\qquad
\xi^0_i=\varepsilon^2(H_{\xi_*}^\varepsilon)^{-1}\left(
\frac{\varepsilon^{-2}}{\sqrt{|\log\varepsilon|}}\hat{\xi}_i^0
\right),
\quad i=l+1,\ldots,m,
\endaligned
$$
where $t_i>0$  and
$\hat{\xi}_i^0\in \mathcal{M}\cap\partial\mathbb{R}_+^2$  satisfy
$t_{i+1}-t_i=\rho$,
$|\hat{\xi}_i^0-\hat{\xi}^0_{i+1}|=\rho$
for all $\rho>0$ sufficiently
small, fixed and independent of $\varepsilon$.
By using the expansion
$(H_{\xi_*}^\varepsilon)^{-1}(z)=\varepsilon^{-2}\xi_*+z+O(\varepsilon^{2}|z|)$
we find
$$
\aligned
\xi^0_i=\xi_*+\frac{1}{\sqrt{|\log\varepsilon|}}\hat{\xi}_i^0
+O\left(\frac{\varepsilon^2}{\sqrt{|\log\varepsilon|}}|\hat{\xi}_i^0|\right),
\quad i=l+1,\ldots,m.
\endaligned
$$
Then it is clear to see
$\xi^0=(\xi^0_1,\ldots,\xi^0_m)\in\mathcal{O}^*_{d,\varepsilon}$
because of $\kappa>1$.
Since $\xi_*\in\po$ is a strict local maximum point of $a(x)$ over $\oo$
and satisfies $\partial_{\nu}a(\xi_*)=\langle\nabla a(\xi_*),\,\nu(\xi_*)\rangle=0$, there exists
a constant $C>0$ independent of $\varepsilon$ such that
$$
\aligned
a(\xi_*)-\frac{C}{\,|\log\varepsilon|\,}\leq a(\xi^0_i)<a(\xi_*),
\qquad i=1,\ldots,m.
\endaligned
$$
From (\ref{7.4})  it follows  that
for any $i=1,\ldots,l$ and $k=1,\ldots,m$ with $i\neq k$,
$$
\aligned
H_a(\xi^0_i,\xi^0_i)=\frac1{4\pi}\log|\log\varepsilon|+O\left(1\right),
\qquad\quad
G_a(\xi^0_k,\xi^0_i)=H_a(\xi^0_k,\xi^0_i)
-\frac1{2\pi}\log|\xi^0_k-\xi^0_i|=\frac1{2\pi}\log|\log\varepsilon|+O\left(1\right).
\endaligned
$$
Moreover, for any $i,k=l+1,\ldots,m$ with $i\neq k$,
$$
\aligned
\log|\xi^0_i-\xi^0_k|=-\frac1{2}\log|\log\varepsilon|+O\left(1\right).
\endaligned
$$
Hence by  (\ref{7.9}),
\begin{eqnarray}\label{7.10}
\max\limits_{\xi\in\overline{\mathcal{O}}^*_{d,\varepsilon}}
F_\lambda(\xi)\geq
F_\lambda(\xi^0)\geq
\frac{8\pi}{p^2\gamma^{2(p-1)}}
\left\{
(m+l)a(\xi_*)|\log\varepsilon|
-\frac12(m+l)(m+l-1)
a(\xi_*)\log|\log\varepsilon|
+O(1)
\right\}.
\end{eqnarray}
Next, we suppose $\xi^\varepsilon=(\xi^\varepsilon_1,\ldots,\xi^\varepsilon_m)\in\partial\mathcal{O}^*_{d,\varepsilon}$.
There are four possibilities:\\
C1. \,\,There exists an $i_0\in\{1,\ldots,l\}$ such that
$\xi^\varepsilon_{i_0}\in\partial B_d(\xi_*)\cap\Omega$, in which case,
$a(\xi^\varepsilon_{i_0})<a(\xi_*)-d_0$ for some\\
\indent\indent $d_0>0$ independent of $\varepsilon$;\\
C2. \,\,There exists an $i_0\in\{l+1,\ldots,m\}$ such that
$\xi^\varepsilon_{i_0}\in\partial B_d(\xi_*)\cap\partial\Omega$, in which case,
$a(\xi^\varepsilon_{i_0})<a(\xi_*)-d_0$ for\\
\indent\indent some $d_0>0$ independent of $\varepsilon$;\\
C3. \,\,There exists an $i_0\in\{1,\ldots,l\}$ such that
$\dist(\xi_{i_0}^\varepsilon,\po)=|\log\varepsilon|^{-\kappa}$;\\
C4. \,\,There exist indices $i_0$, $k_0$, $i_0\neq k_0$ such that
$|\xi_{i_0}^\varepsilon-\xi_{k_0}^\varepsilon|=|\log\varepsilon|^{-\kappa}$.\\
From (\ref{1.7}), (\ref{7.4})  and the maximum principle we have that
 for all $i=1,\ldots,l$ and $k=1,\ldots,m$ with $i\neq k$,
\begin{equation*}\label{7.11}
\aligned
G_a(\xi^\varepsilon_k,\xi^\varepsilon_i)>0,
\qquad\,\,
H_a(\xi^\varepsilon_k,\xi^\varepsilon_i)>0
\qquad\,\,
\textrm{and}
\qquad\,\,
H_a(\xi^\varepsilon_i,\xi^\varepsilon_i)=-\frac1{2\pi}\log\big[\dist(\xi_i^\varepsilon,\po)\big]+O\left(1\right).
\endaligned
\end{equation*}
Thus in the first and second cases,
\begin{eqnarray*}\label{7.12}
\max\limits_{\xi\in\overline{\mathcal{O}}^*_{d,\varepsilon}}
F_\lambda(\xi)=
F_\lambda(\xi^\varepsilon)\leq
\frac{8\pi}{p^2\gamma^{2(p-1)}}
\big\{
\big[(m+l)a(\xi_*)-d_0\big]|\log\varepsilon|
+O\left(\log|\log\varepsilon|\right)
\big\},
\end{eqnarray*}
which contradicts to (\ref{7.10}).
This shows that $a(\xi_i^\varepsilon)\rightarrow a(\xi_*)$. By the condition of $a(x)$ over $\oo$,
we deduce $\xi_i^\varepsilon\rightarrow \xi_*$ for all $i=1,\ldots,m$.\\
In the third case,
\begin{eqnarray}\label{7.13}
\max\limits_{\xi\in\overline{\mathcal{O}}^*_{d,\varepsilon}}
F_\lambda(\xi)=
F_\lambda(\xi^\varepsilon)
\leq
\frac{8\pi}{p^2\gamma^{2(p-1)}}
\left\{
(m+l)a(\xi_*)|\log\varepsilon|
-4\pi a(\xi^\varepsilon_{i_0})H_a(\xi^\varepsilon_{i_0},\xi^\varepsilon_{i_0})
+O(1)
\right\}
\nonumber\\
\leq
\frac{8\pi}{p^2\gamma^{2(p-1)}}
\left\{
(m+l)a(\xi_*)|\log\varepsilon|
-2\kappa a(\xi^\varepsilon_{i_0})\log|\log\varepsilon|
+O(1)
\right\}.
\,\,
\end{eqnarray}
In the last case,  if
$i_0\in\{1,\ldots,m\}$ and $k_0\in\{1,\ldots,l\}$,
\begin{eqnarray}\label{7.14}
\max\limits_{\xi\in\overline{\mathcal{O}}^*_{d,\varepsilon}}
F_\lambda(\xi)=
F_\lambda(\xi^\varepsilon)
\leq
\frac{8\pi}{p^2\gamma^{2(p-1)}}
\left\{
(m+l)a(\xi_*)|\log\varepsilon|
+
2a(\xi^\varepsilon_{i_0})
\log|\xi^\varepsilon_{i_0}-\xi^\varepsilon_{k_0}|
+O(1)
\right\}
\nonumber\\
\leq
\frac{8\pi}{p^2\gamma^{2(p-1)}}
\left\{
(m+l)a(\xi_*)|\log\varepsilon|
-2\kappa a(\xi^\varepsilon_{i_0})\log|\log\varepsilon|
+O(1)
\right\},
\,\ \,
\end{eqnarray}
while if
$i_0\in\{l+1,\ldots,m\}$ and $k_0\in\{l+1,\ldots,m\}$,
\begin{eqnarray}\label{7.15}
\max\limits_{\xi\in\overline{\mathcal{O}}^*_{d,\varepsilon}}
F_\lambda(\xi)=
F_\lambda(\xi^\varepsilon)
\leq
\frac{8\pi}{p^2\gamma^{2(p-1)}}
\left\{
(m+l)a(\xi_*)|\log\varepsilon|
+
a(\xi^\varepsilon_{i_0})
\log|\xi^\varepsilon_{i_0}-\xi^\varepsilon_{k_0}|
+O(1)
\right\}
\nonumber\\
\leq
\frac{8\pi}{p^2\gamma^{2(p-1)}}
\left\{
(m+l)a(\xi_*)|\log\varepsilon|
-\kappa a(\xi^\varepsilon_{i_0})\log|\log\varepsilon|
+O(1)
\right\}.
\,\,\,\,
\end{eqnarray}
Comparing (\ref{7.13})-(\ref{7.15}) with  (\ref{7.10}), we obtain
\begin{equation}\label{7.14}
\aligned
2\kappa a(\xi^\varepsilon_{i_0})\log|\log\varepsilon|
+O(1)\leq
\frac12(m+l)(m+l-1)
a(\xi_*)\log|\log\varepsilon|
+O(1),
\endaligned
\end{equation}
which  is impossible by the choice of $\kappa$ in (\ref{2.5}).
\,\qquad\qquad\qquad\qquad\qquad\qquad\qquad\qquad\qquad\qquad
\qquad\qquad\qquad\qquad\qquad\qquad\quad$\square$

\end{document}